\documentclass[12pt]{amsart}
\usepackage{amssymb}
\usepackage{amsmath}
\usepackage{longtable}
\newcommand\g{{\mathfrak g}}
\newcommand\h{{\mathfrak h}}

\renewcommand{\a}{\mathfrak{a}}
\renewcommand{\c}{\mathfrak{c}}
\newcommand\m{\mathfrak m}
\newcommand\lfr{\mathfrak l}
\newcommand\n{\mathfrak n}

\newcommand\z{\mathfrak z}
\renewcommand\k{\mathfrak k}
\newcommand\s{\mathfrak s}

\newcommand\skewperp{\angle}
\renewcommand{\t}{\mathfrak{t}}
\newcommand\codim{\operatorname{codim}}

\newcommand\tr{\operatorname{tr}}
\newcommand\td{\operatorname{tr.deg}}
\newcommand\im{\operatorname{im}}
\newcommand\Spec{\operatorname{Spec}}
\newcommand\Quot{\operatorname{Quot}}

\newcommand\Aut{\operatorname{Aut}}

\newcommand\C{\mathbb C}

\newcommand\R{\mathbb R}

\newcommand\Z{\mathbb Z}

\renewcommand\sl{\mathfrak{sl}}
\newcommand\so{\mathfrak{so}}
\renewcommand\sp{\mathfrak{sp}}

\newcommand\GL{\mathop{\rm GL}\nolimits}

\newcommand\SL{\mathop{\rm SL}\nolimits}
\newcommand\Sp{\mathop{\rm Sp}\nolimits}
\newcommand\SO{\mathop{\rm SO}\nolimits}

\newcommand\Span{\operatorname{Span}}
\newcommand{\ad}{\mathop{\rm ad}\nolimits}
\newcommand{\Ad}{\mathop{\rm Ad}\nolimits}

\newcommand{\rank}{\mathop{\rm rk}\nolimits}
\newcommand{\cork}{\mathop{\rm cork}\nolimits}
\newcommand{\defe}{\mathop{\rm def}\nolimits}

\renewcommand{\Ad}{\mathop{\rm Ad}\nolimits}

\newcommand\quo{/\!/}
\newtheorem{Thm}{Theorem}[section]
\newtheorem{Prop}[Thm]{Proposition}
\newtheorem{Cor}[Thm]{Corollary}
\newtheorem{Lem}[Thm]{Lemma}
\theoremstyle{definition}
\newtheorem{Ex}[Thm]{Example}
\newtheorem{Conj}[Thm]{Conjecture}
\newtheorem{defi}[Thm]{Definition}
\newtheorem{Rem}[Thm]{Remark}

\unitlength=1mm   \numberwithin{equation}{section}
\numberwithin{table}{section} \oddsidemargin=0cm
\author{Ivan V. Losev}
\title{On fibers of algebraic invariant moment maps}
\begin{document}
%\fontsize{12pt}{20pt}\selectfont
\begin{abstract}
In this paper we study some properties of fibers of the invariant
moment map for a Hamiltonian action of a reductive group on an
affine symplectic variety. We prove that all fibers have equal
dimension. Further, under some additional restrictions, we show that
the quotients of fibers are irreducible normal schemes.
\end{abstract}

\maketitle \tableofcontents
\section{Introduction}
Let $K$ be a connected compact Lie group acting on a symplectic real
manifold $M$ by symplectomorphisms. Suppose there exists a moment
map $\mu:M\rightarrow \k^*$ (see, for instance, \cite{GS} for the
definition of moment maps). It is an important problem in symplectic
geometry to study properties of $\mu$. In fact, usually one studies
not the map $\mu$ itself, but some coarser map, which we call the
{\it invariant moment map}. It is constructed as follows.  One
chooses a Weyl chamber $C\subset \k^*$. The inclusion
$C\hookrightarrow \k^*$ induces a homeomorphism $C\cong \k^*/K$ of
topological spaces. By definition, the invariant moment map $\psi$
is the composition of $\mu:M\rightarrow \k^*$ and the quotient map
$\k^*\rightarrow C$. It turns out that the map $\psi$ has the
following amazing properties provided $M$ is compact:

\begin{itemize}
\item[(a)] The image of $\psi$ is a  convex polytope in $C$.
\item[(b)] All fibers of $\psi$ are connected.
\item[(c)] $\psi$ is an open map onto its image.
\end{itemize}

(a) and (b) were proved by Kirwan in \cite{Kirwan}, (c) is due to
Knop, \cite{Knop_convex}. Since $\mu$ is $K$-equivariant, one can
extract some information about the image of $\mu$ from (a). From (b)
one derives that all fibers of $\mu$ are connected. Hamiltonian
$K$-manifolds satisfying (a)-(c) were called {\it convex} in
\cite{Knop_convex}. In fact,  all interesting classes of Hamiltonian
manifolds (compact manifolds, Stein complex manifolds, cotangent
bundles) are convex, see \cite{Knop_convex} for details.

An algebraic analog of the category of smooth manifolds with an
action of a compact Lie group is the category of smooth {\it affine}
varieties acted on by a reductive algebraic group. Similarly to the
case of compact groups one can define the notion of a Hamiltonian
action of a reductive group, see Subsection \ref{SUBSECTION_Ham1}.
It is an interesting problem to understand:
\begin{enumerate}
\item what are algebraic analogs of properties (a)-(c)?
\item what varieties satisfy these properties?
\end{enumerate}

The study of these two questions was initiated by Knop in the early
90's (see the details below).

In the sequel all groups and varieties are defined over $\C$. First
of all, we need to define the invariant moment map in the algebraic
category. Let $X$ be a symplectic algebraic variety and $G$ a
reductive algebraic group acting on $X$ in a Hamiltonian way. Fix a
moment map $\mu_{G,X}:X\rightarrow \g^*$ for this action. In the
sequel it will be convenient to identify $\g$ and $\g^*$ by means of
a nondegenerate invariant symmetric form of $\g$ and consider
$\mu_{G,X}$ as a morphism $X\rightarrow \g$. By the invariant moment
map for $X$ we mean the morphism
$\psi_{G,X}:=\pi_{G,\g}\circ\mu_{G,X}$, where $\pi_{G,\g}$ denotes
the quotient morphism $\g\rightarrow \g\quo G$ for the adjoint
action $G:\g$. Note that the relation between $\mu_{G,X}$ and
$\psi_{G,X}$ is more loose than in the case of compact groups. For
example, one cannot determine $\im\mu_{G,X}$ by $\im\psi_{G,X}$.

It turns out that the morphism $\psi_{G,X}$ does have some good
properties.

\begin{Thm}\label{Thm:1}
The morphism $\psi_{G,X}$ is equidimensional (i.e., all irreducible
components of nonempty fibers have the same dimensions equal,
obviously, to $\dim X-\dim\overline{\im\psi_{G,X}}$).
\end{Thm}

In fact, a more precise result holds, see Theorem \ref{Thm:4.0.1}.

However, $\psi_{G,X}$ does not seem to have other good properties.
For example, even its general fiber may be disconnected, see
\cite{Knop12}, Introduction. Therefore one needs to modify the
morphism $\psi_{G,X}$.

To this end we introduce a kind of  Stein factorization of
$\psi_{G,X}$. Namely, let $A$  denote  the integral closure of the
subalgebra $\psi_{G,X}^*(\C[\g]^G)$ in $\C[X]^G$. Set
$C_{G,X}:=\Spec(A)$. There are a natural $G$-invariant morphism
$\widetilde{\psi}_{G,X}: X\rightarrow C_{G,X}$ and a finite morphism
$\tau_{G,X}:C_{G,X}\rightarrow \g\quo G$ such that
$\tau_{G,X}\circ\widetilde{\psi}_{G,X}=\psi_{G,X}$. Note that at
least the general fibers of $\widetilde{\psi}_{G,X}:X\rightarrow
C_{G,X}$ are connected whenever $G$ is connected. The idea to
replace $\psi_{G,X}$ with $\widetilde{\psi}_{G,X}$ is due to F.
Knop, see \cite{Knop1}.

In \cite{Knop12} Knop proved that any fiber of
$\widetilde{\psi}_{G,X}$ is connected provided $X$ is the cotangent
bundle of some smooth irreducible (not necessarily affine)
$G$-variety. On the other hand, he constructed an example of a
four-dimensional affine Hamiltonian $\C^\times$-variety $X$ such
that $\widetilde{\psi}_{G,X}$ has a disconnected fiber.

On the other hand, Theorems 1.2.5,1.2.7 from \cite{alg_hamil}
describe the image of $\widetilde{\psi}_{G,X}$. This description is
particularly easy when $X$ satisfies some additional conditions that
can be described as a presence of a grading on $\C[X]$ compatible
with the structure of a Hamiltonian variety.

\begin{defi}\label{defi:2.2.1}
An affine Hamiltonian $G$-variety $X$ equipped with an action
$\C^\times:X$ commuting with the action of $G$ is said to be {\it
conical} if the following two conditions are fulfilled
\begin{itemize}
\item[(Con1)] The morphism $\C^\times\times X\quo G\rightarrow X\quo
G, (t,\pi_{G,X}(x))\mapsto \pi_{G,X}(tx),$ can be extended to a
morphism $\C\times X\quo G\rightarrow X\quo G$.
\item[(Con2)] There exists a positive integer $k$ (called the {\it degree} of $X$) such that $t_*\omega=t^{-k}\omega$ and
$\mu_{G,X}(tx)=t^k\mu_{G,X}(x)$  for all $t\in \C^\times, x\in X$. Here $\omega$ denotes the symplectic form
on $X$ and $t_*\omega$ is the push-forward of $\omega$ under the automorphism of $X$ induced by $t$.
\end{itemize}
\end{defi}

For example, a symplectic $G$-module and the cotangent bundle of a
smooth affine $G$-variety are conical.

If $X$ is conical, then $C_{G,X}$ is a quotient of a vector space by
a finite group and $\widetilde{\psi}_{G,X}$ is surjective, see
\cite{alg_hamil}, Theorem 1.2.7. More precisely, there is a subspace
$\a\subset \g$ (called the Cartan space of $X$) and a subgroup
$W\subset N_G(\a)/Z_G(\a)$ (the Weyl group) such that $C_{G,X}\cong
\a/W$ and the finite morphism $\tau_{G,X}:C_{G,X}\rightarrow \g\quo
G$ is induced by the embedding $\a\hookrightarrow \g$. So the
subspace $\a\subset \g$ and the group $W$ encode the difference
between $\widetilde{\psi}_{G,X}$ and $\psi_{G,X}$. This description
partially generalizes Knop's results for cotangent bundles and
symplectic vector spaces (\cite{Knop1},\cite{Knop6}).

We have no examples of conical Hamiltonian $G$-varieties, where
$\widetilde{\psi}_{G,X}$ has a disconnected fiber. We conjecture
that in this case all fibers of $\widetilde{\psi}_{G,X}$ are
connected and, more precisely, that $X$ enjoys the following
property:

\begin{itemize}
\item[(Irr)] Any fiber of $\widetilde{\psi}_{G,X}\quo G:X\quo G\rightarrow C_{G,X}$ is irreducible.
\end{itemize}

We are able to prove (Irr) only under another restriction on $X$.

\begin{defi}\label{Def:0.4}
An affine Hamiltonian $G$-variety $X$ is said to be {\it untwisted}
if
\begin{itemize}
\item[(Utw1)] $C_{G,X}$ is smooth.
\item[(Utw2)] The morphism $\widetilde{\psi}_{G,X}$ is smooth in
codimension 1 (that is, the complement to the set of smooth points
of $\widetilde{\psi}_{G,X}$ in $X$ has codimension at least 2).
\end{itemize}
\end{defi}

\begin{Thm}\label{Thm:0.5}
Let $G$ be connected and $X$ a conical Hamiltonian $G$-variety.
\begin{enumerate}
\item If $X$ is untwisted, then any fiber of $\widetilde{\psi}_{G,X}\quo G$ is a normal Cohen-Macaulay scheme.
\item If $X$ satisfies (Utw1) and all fibers  of $\widetilde{\psi}_{G,X}\quo G$ are normal (as schemes), then
$X$ satisfies (Irr).
\item Suppose $X$ is algebraically simply connected. If $X$ satisfies (Irr), then $X$ is untwisted.
\end{enumerate}
\end{Thm}

The term "untwisted" is partially justified by Remark \ref{Rem:0.6}.

We recall that a smooth irreducible variety $X$ is called {\it
algebraically simply connected} if a finite \'{e}tale morphism
$\varphi:Y\rightarrow X$ is an isomorphism whenever $Y$ is
irreducible.

Note that a fiber of $\widetilde{\psi}_{G,X}\quo G$ can be thought
as an algebraic analog of a Marsden-Weinstein reduction, \cite{MW}.

Now let us describe some classes of conical untwisted Hamiltonian
$G$-varieties. Knop showed in \cite{Knop2} that the cotangent bundle
of any smooth irreducible affine variety is untwisted. In the
present paper we give  alternative proofs of this result and prove
that a symplectic $G$-module is untwisted.

Let us briefly describe the content of the paper. In Section
\ref{SECTION_prelim} we recall some known results concerning
Hamiltonian actions in the algebraic setting. Section
\ref{SECTION_dimension} is devoted to the proof of  Theorem
\ref{Thm:1} (in fact, of a more precise statement). In Section
\ref{SECTION_Weyl} we prove some results concerning the Weyl groups
of Hamiltonian actions (see above). These results are used in the
proof of Theorem \ref{Thm:0.5}. Besides, they play a crucial role in
the computation of Weyl groups and  root lattices of affine
$G$-varieties, the former is done in the preprint \cite{Weyl}.
Section \ref{SECTION_untwisted} is devoted to the proof of Theorem
\ref{Thm:0.5}. We also present there some classes of untwisted
varieties. In Section \ref{SECTION_open} we discuss some open
problems related to the subject of the paper. Finally, Section
\ref{SECTION_Notation} contains conventions and the list of notation
we use. In the beginning of  Sections 2-5 their content is
described in more detail.

{\bf Acknowledgements.} Part of the work on this paper was done
during my visit  to Ruhr University, Bochum, in July, 2005, in the
framework  of Euler program. I would like to thank this institution
and especially Professor H. Flenner for hospitality. I also express
my gratitude to Professor F. Knop for his kind permission to use his
counterexample in Subsection \ref{SUBSECTION_counterexamples}.
Finally, I wish to thank the referees for useful remarks on an
earlier version of this text.

\section{Preliminaries}\label{SECTION_prelim}
In this section $G$ is a reductive algebraic group and $X$ is a
smooth variety equipped with a regular symplectic form $\omega$ and
an action of $G$ by symplectomorphisms.

In Subsection \ref{SUBSECTION_Ham1} we recall the definition of a
Hamiltonian action and give some examples.  Subsection
\ref{SUBSECTION_Ham2} is devoted to conical Hamiltonian varieties
introduced in \cite{alg_hamil}. In Subsection
~\ref{SUBSECTION_local} we study a local structure of Hamiltonian
actions. At first, we recall the theory of cross-sections of
Hamiltonian actions (Proposition \ref{Prop:1.1}) tracing back to
Guillemin-Sternberg, \cite{GS}. Next, in this subsection we recall
the symplectic slice theorem from \cite{slice}. These two results
are  key ingredients of most proofs in this paper. Finally, in
Subsection \ref{SUBSECTION_recall} we recall some results from
\cite{alg_hamil}, \cite{Comb_Ham}. The most important ones are
Propositions \ref{Lem:4.4.1}, \ref{Thm:2.2}.

\subsection{Hamiltonian actions}\label{SUBSECTION_Ham1}
 Let $U$ be an open subset of $X$ and $f$ a regular
function on $U$.
 The {\it skew-gradient} $v(f)$ of $f$ is, by definition, the regular vector field on $U$
 given by the equality
\begin{equation*}%\label{eq:2.1:1}
\omega_x(v(f),\eta)=\langle d_xf, \eta\rangle, x\in U, \eta\in T_xX.
\end{equation*}
 For $f,g\in \C[U]$ one  defines their Poisson bracket
$\{f,g\}\in \C[U]$ by
\begin{equation*}%\label{eq:2.1.2}
\{f,g\}=\omega(v(f),v(g)).
\end{equation*}
Clearly, $\{f,g\}=L_{v(f)}g$, where $L$ denotes the Lie derivative.

To any element $\xi\in\g$ one associates the velocity vector field
$\xi_*$. Suppose  there is a linear map $\g\rightarrow \C[X],
\xi\mapsto H_\xi,$ satisfying the following two conditions:
\begin{itemize}
\item[(H1)] The map $\xi\mapsto H_\xi$ is $G$-equivariant.
\item[(H2)] $v(H_\xi)=\xi_*$.
\end{itemize}

\begin{defi}\label{Def:2.1.1}
The action  $G:X$ equipped with a linear map $\xi\mapsto H_\xi$
satisfying (H1),(H2) is said to be {\it Hamiltonian} and $X$ is
called a Hamiltonian $G$-variety. %The functions $H_\xi$ are said to
%be the hamiltonians of $X$.
\end{defi}

\begin{Rem}\label{Rem:2.1.2}
Very often the definition of a Hamiltonian action is given in a
slightly different way. Namely, for a connected group $G$  condition
(H1) is replaced by the condition $\{H_\xi,H_\eta\}=H_{[\xi,\eta]}$.
However, these two conditions are equivalent provided (H2) is
fulfilled. Note also that one can consider Hamiltonian actions on arbitrary Poisson varieties,
see, for example, \cite{alg_hamil}.
\end{Rem}

For a Hamiltonian action $G:X$ we define the morphism
$\mu_{G,X}:X\rightarrow \g^*$ by the formula
\begin{equation*}%\label{eq:2.1:3}
\langle \mu_{G,X}(x),\xi\rangle= H_{\xi}(x),\xi\in\g,x\in X.
\end{equation*}
This morphism is called the {\it moment map} of the Hamiltonian
$G$-variety $X$.

Conditions (H1),(H2) are equivalent, respectively, to

\begin{itemize}
\item[(M1)] $\mu_{G,X}$ is  $G$-equivariant.
\item[(M2)] $\langle d_x\mu_{G,X}(v),\xi\rangle=\omega_x(\xi_x,v),$ for all $ x\in X,v\in
T_xX,\xi\in\g$.
\end{itemize}

Here and below we write $\xi_x$ instead of $\xi_{*x}$.

Any two maps $\mu_{G,X}:X\rightarrow \g^*$ satisfying  conditions
(M1),(M2)  differ by an element of $\g^{*G}$. Moreover,
$H_{[\xi,\eta]}=\{H_\xi, H_\eta\} =\omega(\xi_*,\eta_*)$ (see, for
example, \cite{GS},\cite{Vinberg}). Conversely, for any
$\eta\in\g^{*G}$ there exists the unique Hamiltonian $G$-variety
$X_{\eta}$ coinciding with $X$ as a symplectic $G$-variety and such
that $\mu_{G,X_\eta}=\mu_{G,X}+\eta$.

Let us choose some effective $G$-module $V$ and put
$(\xi,\eta)=\tr_V(\xi\eta)$ for $\xi,\eta\in\g$. The form $(\cdot,\cdot)$ is
$G$-invariant, symmetric and its restriction to the Lie algebra of
any reductive subgroup of $G$ is nondegenerate. Using this form, we
identify $\g$ and $\g^*$. In particular, we may consider $\mu_{G,X}$
as a morphism from $X$ to $\g$.

Let us now give some examples of Hamiltonian $G$-varieties.

\begin{Ex}[Cotangent bundles]\label{Ex:2.1.4} Let $X_0$ be a smooth $G$-variety, $X:=T^*X_0$ the cotangent
bundle of $X_0$. $X$ is a symplectic algebraic variety (the
symplectic form is presented, for example, in
\cite{GS},\cite{Vinberg}). The action of $G$ on $X$ is Hamiltonian.
The moment map is given by $\langle\mu_{G,X}((y,\alpha)),
\xi\rangle=\langle \alpha, \xi_{y}\rangle$. Here $y\in X_0,
\alpha\in T^*_yX_0,\xi\in\g$.
\end{Ex}

\begin{Ex}[Symplectic vector spaces]\label{Ex:2.1.5}  Let $V$ be a vector space equipped with a non-degenerate
skew-symmetric bilinear form $\omega$. Then $V$ is a symplectic
variety. Let $G$   act on $V$ by linear symplectomorphisms. Then the
action $G:V$ is Hamiltonian. The moment map $\mu_{G,V}$ is given by
$\langle\mu_{G,V}(v), \xi\rangle = \frac{1}{2}\omega(\xi v,v),
\xi\in\g, v\in V$.
\end{Ex}

\begin{Ex}[Model varieties]\label{Ex:2.1.6}
This example generalizes the previous one.
 Let $H$  be a reductive subgroup of $G$, $\eta\in \g^H$, $V$ a symplectic $H$-module. Put $U=(\z_\g(\eta)/\h)^*$. Let us
equip the homogeneous vector bundle $X=G*_H(U\oplus V)$ with a
certain closed 2-form. Let $\eta_n,\eta_s$ denote nilpotent and
semisimple parts of $\eta$, respectively. If $\eta_n\neq 0$, choose
an $\sl_2$-triple $(\eta_n,h,f)$ in $\z_\g(\eta_s)^H$ (where $h$ is
semisimple and $f$ is nilpotent). If $\eta_n=0$, we set $h=f=0$. The
$H$-module $U$ can be identified with $\z_\g(f)\cap \h^\perp$. Fix a
point $x=[1,(u,v)]\in X$. The tangent space $T_xX$ is naturally
identified with $\h^\perp\oplus U\oplus V$, where $U\oplus V$ is the
tangent space to the fiber of the projection $G*_H(U\oplus
V)\rightarrow G/H$ and the embedding $\h^\perp\hookrightarrow T_xX$
is given by $\xi\mapsto \xi_x$. Put
\begin{equation*}%\label{eq:model_form}
\begin{split}
\omega_x(u_1+v_1+\xi_1,&u_2+v_2+\xi_2)=\omega_V(v_1,v_2)+(\xi_1,u_2)-(\xi_2,u_1)+(\eta+u+\mu_{H,V}(v),[\xi_1,\xi_2]),\\
& u_1,u_2\in U, v_1,v_2\in V, \xi_1,\xi_2\in \h^\perp.
\end{split}
\end{equation*}

The corresponding map $\omega: U\oplus V\rightarrow \bigwedge^2
(\h^\perp\oplus U\oplus V)^*$ is $H$-equivariant. Thus $\omega$ can
be extended to the unique $G$-invariant 2-form on $X$, which is
denoted also by $\omega$. It turns out that $\omega$ is closed and
nondegenerate in any point of the zero section $G/H$, \cite{slice},
assertion 1 of Proposition 1. If $\eta$ is nilpotent, then $\omega$
is nondegenerate on the whole variety $X$. In the general case the
subset $X_r=\{x\in G*_H(U\oplus V)| \omega_x \text{ is nondegenerate
in }x\}$ is affine. The action $G:X_r$ is Hamiltonian. The moment
map is given by  (see \cite{slice}, assertion 3 of Proposition 1)
$$\mu_{G,X_r}([g,(u,v)])=\Ad(g)(\eta+u+\mu_{H,V}(v)).$$
 We denote the Hamiltonian variety $X_r$ by
$M_{G}(H,\eta,V)\index{mghev@$M_G(H,\eta,V)$}$ and call it a {\it
model variety}.
\end{Ex}

\begin{Rem}\label{Rem:2.1.7}
The Hamiltonian structure on $M_G(H,\eta,V)$ depends on the choice
of an $\sl_2$-triple $(\eta_n,h,f)$ in $\z_\g(\eta_s)^H$ (if
$\eta_n\neq 0$). However, Hamiltonian varieties corresponding to
different choices of $h,f$ are isomorphic (see Remark 1
from~\cite{slice}). In the sequel we say that $(\eta_n,h,f)$ is an
$\sl_2$-triple {\it generating} $M_G(H,\eta,V)$.
\end{Rem}

\begin{Rem}\label{Rem:2.1.8}
For $\eta_0\in \g^G$ the Hamiltonian  $G$-varieties
$M_G(H,\eta+\eta_0,V), M_G(H,\eta,V)_{\eta_0}$ are naturally
identified. They even coincide as  subsets in $G*_H(U\oplus V)$.
\end{Rem}

Now we consider two constructions with Hamiltonian varieties.

\begin{Ex}[Restriction to a subgroup]\label{Ex:2.1.9} Let $H$ be a reductive subgroup
of $G$ and $X$ a Hamiltonian $G$-variety. Then $X$ is a Hamiltonian
$H$-variety with the moment map $\mu_{H,X}=p\circ\mu_{G,X}$. Here
$p$ denotes the restriction map $\g^*\twoheadrightarrow \h^*$.
\end{Ex}

\begin{Ex}[Products]\label{Ex:2.1.10}
Suppose $X_1,X_2$ are Hamiltonian $G$-varieties. Being the product
of symplectic varieties,  the variety $X_1\times X_2$ has a natural
symplectic structure. The action $G:X_1\times X_2$ is Hamiltonian.
The moment map is given by the formula $\mu_{G, X_1\times
X_2}(x_1,x_2)=\mu_{G,X_1}(x_1)+\mu_{G,X_2}(x_2)$ for $x_1\in
X_1,x_2\in X_2$.
\end{Ex}

\begin{Rem}\label{Rem:2.1.11}
It follows directly from the construction of a model variety that if
$(H,\eta,V)$ is the same as in Example~\ref{Ex:2.1.6} and $V_0$ is a
trivial symplectic $H$-module, then the Hamiltonian $G$-varieties
$M_G(H,\eta,V\oplus V_0)\cong M_G(H,\eta,V)\times V_0$ are
isomorphic (the action $G:V_0$ is assumed to be trivial).
\end{Rem}

Now we define some important numerical invariants of an irreducible
Hamiltonian $G$-variety $X$. For an action of $G$ on an algebraic
variety $Y$ we denote by $m_G(Y)\index{mgy@$m_G(Y)$}$ the maximal
dimension of a $G$-orbit on $Y$.  The number
$m_G(X)-m_G(\overline{\im\mu_{G,X}})$ is called the {\it defect} of
$X$ and is denoted by $\defe_G(X)\index{defgx@$\defe_G(X)$}$. The
number $\dim X-\defe_G(X)-m_G(X)$ is called the {\it corank} of $X$
and is denoted by $\cork_G(X)\index{corkgx@$\cork_G(X)$}$.
Equivalently, $\cork_G(X)=\td \C(X)^G-\defe_G(X)$. An irreducible
Hamiltonian $G$-variety $X$ such that $\cork_G(X)=0$ is called {\it
coisotropic}.

It follows from the standard properties of the moment map (see, for
example, \cite{GS},\cite{Vinberg}) that the  defect and the corank
of $X$ coincide, respectively, with  $\dim\ker\omega|_{\g_*x}$,
$\rank\omega|_{(\g_*x)^\skewperp}$ for a point $x\in X$ in general
position. Further, the following statement holds, see
\cite{alg_hamil}, Proposition 3.1.7.
\begin{Lem}\label{Lem:2.1.12}
$\dim C_{G,X}=\dim\overline{\im\psi_{G,X}}=\defe_{G}(X)$.
\end{Lem}

\begin{defi}\label{defi:3.5} Let $X_1,X_2$ be Hamiltonian
$G$-varieties. A morphism $\varphi:X_1\rightarrow X_2$ is called
{\it Hamiltonian} if it is an \'{e}tale $G$-equivariant
symplectomorphism intertwining the moment maps.
\end{defi}

Note that a Hamiltonian morphism $\varphi:X_1\rightarrow X_2$
induces the unique morphism $\varphi_0:C_{G,X_1}\rightarrow C_{G,X_2}$
such that
$\widetilde{\psi}_{G,X_2}\circ\varphi=\varphi_0\circ\widetilde{\psi}_{G,X_1}$.

\begin{Rem}\label{Rem:2.1.13}
One can similarly define Hamiltonian actions on complex analytic
manifolds. The definitions of the corank and the defect can be
extended to this case without any  noticeable modifications.
\end{Rem}

\subsection{Conical Hamiltonian varieties}\label{SUBSECTION_Ham2}
The definition of a conical Hamiltonian variety was given in
Introduction, Definition \ref{defi:2.2.1}.

\begin{Ex}[Cotangent bundles]\label{Ex:2.2.2}
Let $X_0,X$ be as in Example \ref{Ex:2.1.4}. The variety $X$ is a
vector bundle over $X_0$. The action $\C^\times:X$ by  fiberwise
multiplication turns $X$ into a conical variety of degree 1.
\end{Ex}

\begin{Ex}[Symplectic vector spaces]\label{Ex:2.2.3}
The symplectic $G$-module $V$ equipped with the action $\C^\times:V$
given by $(t,v)\mapsto tv$ is conical of degree 2.
\end{Ex}

\begin{Ex}[Model varieties]\label{Ex:2.2.4}
This example generalizes the previous one. Let $H,\eta,V$ be  as in
Example \ref{Ex:2.1.6} and $X=M_G(H,\eta,V)$. Suppose that $\eta$ is
nilpotent. Here we define an action $\C^\times:X$ turning $X$ into a
conical Hamiltonian variety of degree 2. Let $(\eta,h,f)$ be the
$\sl_2$-triple in $\g^H$ generating $X$. As a $G$-variety,
$X=G*_H(U\oplus V)$, where $U=\z_\g(f)\cap \h^\perp$. Note that $h$
is an image of a coroot under an embedding of Lie algebras.
 In particular, there
exists a one-parameter subgroup $\gamma:\C^\times\rightarrow G$ with
$\frac{d}{dt}|_{t=0}\gamma=h$. Since $[h,\h]=0, [h,f]=-2f$, we see
that  $\gamma(t)(\h^\perp)=\h^\perp,\gamma(t)(U)=U$. Define a
morphism $\C^\times\times X\rightarrow X$ by formula
\begin{equation}\label{eq:2.4}
(t, [g,(u,v)])\mapsto [g\gamma(t),t^2\gamma(t)^{-1}u,tv], t\in
\C^\times, g\in G, u\in U, v\in V.
\end{equation}
One checks directly that the  morphism (\ref{eq:2.4}) is
well-defined and determines an action of $\C^\times$ on $X$
commuting with the action of $G$. Let us check that $X$ with this
action is a conical Hamiltonian variety. The action of $\C^\times$
on $X\quo G$ coincides with that induced by the action $\C^\times:X$
given by
\begin{equation}\label{eq:2.2:1} (t,[g,(u,v)])\mapsto
[g,t^2\gamma(t)^{-1}u,tv].\end{equation} The eigenvalues of $\ad(h)$
on $\z_\g(f)$ are not positive. Thus the  morphism (\ref{eq:2.2:1}) can
be extended to a morphism $\C\times X\rightarrow X$. This yields
(Con1). (Con2) for $k=2$ is verified directly using the construction
of Example~\ref{Ex:2.1.6}.
\end{Ex}

\begin{Rem}\label{Rem:2.2.5}
Let $X$ be as in the previous example. The action $\C^\times:X$
induces a non-negative grading on $\C[X]^G$. In the notation of the
previous example $\C[X]^G\cong \C[U\oplus V]^H$. The grading on
$\C[U\oplus V]^H$ is induced from the following grading on
$\C[U\oplus V]$:

all elements of $V^*\subset \C[U\oplus V]$ have degree 1. The
$H$-module
 $U^*$ is naturally identified with $\z_\g(\eta)\cap
 \h^\perp$. Put $\g_i=\{\xi\in \g| [h,\xi]=i \xi\}$. All elements of $\z_\g(\eta)\cap\h^\perp\cap\g_i$
 have degree $i+2$.\end{Rem}

\begin{Lem}[\cite{alg_hamil}, Lemma 3.3.6]\label{Lem:2.2.6}
Let $X$ be a conical Hamiltonian $G$-variety of degree $k$. Then
\begin{enumerate}
\item $0\in \im\psi_{G,X}$.
\item Assume that $X$ is irreducible and normal. Then the subalgebra $\C[C_{G,X}]\subset
\C[X]^G$ is $\C^\times$-stable. The morphisms
$\widetilde{\psi}_{G,X}:X\rightarrow C_{G,X},
\tau_{G,X}:C_{G,X}\rightarrow \g\quo G$ are $\C^\times$-equivariant,
where the action $\C^\times:\g\quo G$ is induced from the action
$\C^\times:\g$ given by $(t,x)\mapsto t^kx, t\in \C^\times, x\in\g$.
\item Under the assumptions of  assertion 2, there is the unique point $\lambda_0\in C_{G,X}$ such that
$\tau_{G,X}(\lambda_0)=0$. For any point $\lambda\in C_{G,X}$ the
limit $\lim_{t\rightarrow 0}t\lambda$ exists and is equal to
$\lambda_0$.
\end{enumerate}
\end{Lem}

%\begin{Rem}\label{Rem:2.2.7}
%It follows easily from Lemma~\ref{Lem:2.2.6} that
%$0\in\a_{G,X}^{(X_L)}$ provided $X$ is conical.
%\end{Rem}

\subsection{Local structure of Hamiltonian actions}\label{SUBSECTION_local}
Firstly, we review the algebraic variant of  the Guillemin-Sternberg
local cross-section theory, see \cite{Knop2}, Section 5,
\cite{alg_hamil}, Subsection 5.1. Let $L$ be a Levi subgroup of $G$
and $\lfr$ the corresponding Lie algebra. Put
$\lfr^{pr}=\{\xi\in\lfr| \z_\g(\xi_s)\subset \lfr\}$.

\begin{Prop}[\cite{Knop2}, Theorem 5.4 and \cite{alg_hamil},
Corollary 5.1.3, Propositions 5.1.2, 5.1.4, 5.1.7]\label{Prop:1.1}
Let $x\in X,
\lfr=\z_\g(\mu_{G,X}(x)_s),Y=\mu_{G,X}^{-1}(\lfr^{pr})$. Then
\begin{enumerate}
\item
 $T_yX=\lfr^\perp_*y\oplus T_yY$
is a skew-orthogonal direct sum  for any $y\in Y$. In particular,
$Y$ is a smooth subvariety of $X$ and the restriction of $\omega$ to
$Y$ is nondegenerate. Thus $Y$ is equipped with a symplectic
structure.
\item The action $N_G(L):Y$ is Hamiltonian with the
moment map $\mu_{G,X}|_Y$. \item The natural morphism
$G*_{N_G(L)}Y\rightarrow X$ is \'{e}tale. Its image is saturated.
\item If $x$ is in general position, then the natural morphism $G*_{N_G(L)}Y\rightarrow
X$ is an open embedding and $N_G(L)$ permutes the connected
components of $Y$ transitively.
\end{enumerate}
\end{Prop}

A subset $Z^0$ of a $G$-variety $Z$ is said to be {\it saturated} if
there exist a $G$-invariant morphism $\varphi:Z\rightarrow Z_0$ and
a subset $Z_0^0\subset Z_0$ such that $Z^0=\varphi^{-1}(Z_0^0)$.

\begin{defi}\label{defi:1.2}
An irreducible (=connected) component of $\mu_{G,X}^{-1}(\lfr^{pr})$
equipped with the structure of a Hamiltonian $L$-variety obtained by
restriction of the Hamiltonian structure from
$\mu_{G,X}^{-1}(\lfr^{pr})$ is called an {\it $L$-cross-section} of
$X$.
\end{defi}

\begin{defi}\label{defi:1.3}
The Levi subgroup $L=Z_G(\mu_{G,X}(x)_s)$, where $x\in X$ is in
general position, is said to be the {\it principal centralizer} of
$X$.
\end{defi}

Note that the principal centralizer is determined uniquely up to
$G$-conjugacy.

\begin{Lem}\label{Lem:2.3.5}
Let $L$ be the principal centralizer and $X_L$ an $L$-cross-section
of $X$. Then the following conditions are equivalent:
\begin{enumerate}
\item $m_G(X)=\dim G$.
\item $\defe_G(X)=\rank G$.
\item $\overline{\im \mu_{G,X}}=\g$.
\item $L$ is a maximal torus in $G$ and $m_L(X_L)=\defe_L(X_L)=\rank G$.
\item The stabilizer in general position for the action $G:X$ is finite.
\end{enumerate}
Under these conditions, $\cork_G(X)=\dim X-\dim G-\rank G$.
\end{Lem}
\begin{proof}
The equivalence of conditions (1)-(4) was proved in \cite{Comb_Ham},
Lemma 4.5. The  equality for $\cork_G(X)$ follows from (1) and (2). It is well-known that (5) is equivalent to
(1).
\end{proof}

%REFEREE'S REMARK: more detailed proof
\begin{Lem}\label{Lem:2.3.11}
Let $L$ be the principal centralizer and $X_L$ an $L$-cross-section
of $X$. Suppose that the stabilizer in general position $L_0$ for
the action $L:X_L$ is reductive and that $0\in
\overline{\im\psi_{G,X}}$. Then
$\overline{\im\mu_{G,X}}=\overline{G(\lfr\cap\lfr_0^\perp)}$.
\end{Lem}
\begin{proof}
From Proposition \ref{Prop:1.1} it follows that
$\overline{\im\mu_{G,X}}=\overline{G\im\mu_{L,X_L}}$. Since $0\in
\overline{\im\psi_{G,X}}$, we see that
$0\in\overline{\im\mu_{L,X_L}}$. By Theorem 4.1.1 from
\cite{alg_hamil}, $(L,L)\subset L_0$. Therefore
$\overline{\im\mu_{L,X_L}}=\lfr\cap\lfr_0^\perp$.
\end{proof}

Now we turn to the problem of describing the structure of an {\it
affine} Hamiltonian $G$-variety in some neighborhood of a point with
closed $G$-orbit. A neighborhood is taken with respect to the
complex topology (in the sequel we call such neighborhoods {\it
analytical}).

At first, we define some invariants of the triple $(G,X,x)$. Put
$H=G_x, \eta=\mu_{G,X}(x)$. The subgroup $H\subset G$ is reductive
and $\eta\in \g^H$. Put $V=(\g_*x)^\skewperp/(\g_*x\cap
\g_*x^\skewperp)$. This is a symplectic $H$-module. We say that
$(H,\eta,V)$  is the {\it determining triple} of  $X$ at $x$. For
example, the determining triple of $X=M_G(H,\eta,V)$ at
$x=[1,(0,0)]$ is $(H,\eta,V)$, see \cite{slice}, assertion 4 of
Proposition 1.

As the name suggests, a determining triple should determine the
structure of the Hamiltonian  $G$-variety  $X$ near $x$. In fact, a
slightly stronger claim holds.

\begin{defi}\label{defi:4.3.1}
Let $X_1,X_2$ be affine Hamiltonian $G$-varieties, $x_1\in X_1,
x_2\in X_2$ be points with closed  $G$-orbits. The pairs
$(X_1,x_1),(X_2,x_2)$ are called {\it analytically equivalent}, if
there are saturated open analytical neighborhoods $O_1,O_2$ of
$x_1\in X_1, x_2\in X_2$, respectively, that are isomorphic as
complex-analytical Hamiltonian $G$-manifolds.
\end{defi}

\begin{Rem}\label{Rem:4.3.2}
An open saturated analytical neighborhood in  $X$ is the inverse
image of an {\it open} analytical neighborhood in $X\quo G$ under
$\pi_{G,X}$. See, for example, \cite{slice}, Lemma 5.
\end{Rem}

\begin{Prop}[Symplectic slice theorem, \cite{slice}]\label{Prop:4.3.3}
Let $X$ be an affine Hamiltonian  $G$-variety,  $x\in X$  a point
with closed $G$-orbit, $(H,\eta,V)$  the determining triple of $X$
at $x$. Then the pair $(X,x)$ is analytically equivalent to the pair
$(M_{G}(H,\eta,V), [1,(0,0)])$.
\end{Prop}

Now we prove two lemmas, which will be used in Subsection
\ref{SUBSECTION_affham6}.

We have two approaches to the local study of affine Hamiltonian
varieties: the cross-sections theory  and the symplectic slice
theorem. Let us establish a connection between them.

\begin{Lem}\label{Lem:4.3.4} Let $x\in X$ be a point with closed   $G$-orbit
and $(H,\eta,V)$ the determining triple of $X$ at $x$. Put
$M=Z_{G}(\eta_s)$. Denote by $X_M$ the unique $M$-cross-section of $X$
containing $x$.  Then the following assertions hold
\begin{enumerate}\item  $Mx$ is closed in $X_M$ and $(H,\eta,V)$
is the determining triple of  $X_M$ at $x$.
\item There exists an affine saturated open (with respect to Zariski topology) neighborhood
$X_M^0\subset X_M$ of $x$ such that the following conditions are
satisfied:
\begin{enumerate}
\item the natural morphism $X_M^0\quo M\rightarrow X\quo G,
\pi_{M,X_M}(z)\mapsto \pi_{G,X}(z)$ is \'{e}tale;
\item
for any $z\in X^0_M$ the orbit
 $Mz$ is closed in $X_M^0$ (equivalently, in $X_M$) iff $Gz$ is closed in $X$.
\end{enumerate}
\end{enumerate}
\end{Lem}
\begin{proof}
The morphism $\varphi:G*_MX_M\rightarrow X, [g,x]\mapsto gx,$ is
\'{e}tale (assertion 3 of Proposition~\ref{Prop:1.1}). Since $Gx$ is
closed in $X$, we see that $G[1,x]$ is closed in $G*_MX_M$,
equivalently,   $Mx$ is closed in $X_M$. Since $G_z\subset
Z_G(\mu_{G,X}(z))\subset Z_G(\mu_{G,X}(z)_s)=M$, we have $G_z=M_z$
for $z\in X_M$. By construction of $\mu_{M,X_M}$,
$\mu_{M,X_M}(z)=\mu_{G,X}(z)$. Assertion 1 will follow if we check
that the $H$-modules $\g_*x^\skewperp/(\g_*x^\skewperp\cap \g_*x)$
and $\m_*x^\skewperp/(\m_*x^\skewperp\cap \m_*x)$ are isomorphic.
Here the skew-orthogonal complement to $\g_*x$ (resp., to $\m_*x$)
is taken in $T_xX$ (resp., in $T_xX_M$)). The existence of an
isomorphism stems from $\g_*x=\m^{\perp}_*x\oplus \m_*x$ and
assertion 1 of Proposition~\ref{Prop:1.1}.

By the above, the orbits $G[1,x],Gx$ are closed and isomorphic via
$\varphi$. It follows from Luna's fundamental lemma, \cite{Luna1},
that for some open affine neighborhood $U$ of the point
$\pi_{M,X_M}(x)$ in $X_M\quo M\cong (G*_MX_M)\quo G$ the morphism
$\varphi\quo G:U\rightarrow X\quo G$ is \'{e}tale and
\begin{equation}\label{eq:4.3:1}\pi_{G,G*_MX_M}^{-1}(U)\cong
U\times_{X\quo G}X.\end{equation}
 Clearly, $\pi_{G,G*_MX_M}^{-1}(U)\cong G*_M\pi_{M,X_M}^{-1}(U)$.
 Thanks to (\ref{eq:4.3:1}), we see that for all $z\in
X_M^0:=\pi^{-1}_{M,X_M}(U)$ the orbit  $G[1,z]$ is closed in $
G*_M\pi_{M,X_M}^{-1}(U)$ iff $Gz$ is closed in $X$.
\end{proof}

The next lemma studies the behavior of determining triples under
replacing $G$ with some connected subgroup $G^1\subset G$ containing
$(G,G)$.

\begin{Lem}\label{Lem:4.3.5}
Let $x\in X$ be a point with closed $G$-orbit and $(H,\eta,V)$ the
determining triple of $X$ at $x$. Then $G^1x$ is closed in $X$ and
the determining triple of the Hamiltonian  $G^1$-variety $X$ at $x$
has the form $(H\cap G^1, \eta_0,V\oplus V_0)$, where $V_0$ is a
trivial  $H\cap G^1$-module and $\eta_0$ is the projection of $\eta$
to $\g^1$.
\end{Lem}
\begin{proof}
Since $G^1$ is a normal subgroup of $G$, we see that all
$G^1$-orbits in $Gx$ have the same dimension whence closed.
Obviously, $G^1_x=G^1\cap H, \mu_{G^1,X}(x)=\eta_0$. Clearly,
 $\g^1_*x\subset \g_*x$ and $\g_*x^\skewperp\subset
\g^1_*x^\skewperp$. Therefore we have a natural  embedding
$\g_*x^\skewperp/(\g_*x^\skewperp\cap \g^1_*x)\hookrightarrow
\g^1_*x^\skewperp/(\g^1_*x^\skewperp\cap \g^1_*x)$ and a natural
projection $\g_*x^\skewperp/(\g_*x^\skewperp\cap
\g^1_*x)\twoheadrightarrow \g_*x^\skewperp/(\g_*x^\skewperp\cap
\g_*x)$. The cokernel of the former is a quotient of the $H\cap
G^1$-module $\g^1_*x^\skewperp/\g_*x^\skewperp\cong
(\g_*x/\g^1_*x)^*$, while the kernel of the latter is a submodule in
$\g_*x/\g^1_*x$. Since $\g_*x/\g_*^1x$ is a trivial $H\cap
G^1$-module, we are done.
\end{proof}

\subsection{Some results concerning
$\widetilde{\psi}_{G,X},C_{G,X}$}\label{SUBSECTION_recall} Let us,
at first, define two important invariants of a Hamiltonian variety:
its Cartan space and Weyl group. The proofs of the facts below
concerning these invariants can be found in \cite{alg_hamil},
Subsection 5.2.

Let $L$ be the principal centralizer and $X_L$ an $L$-cross-section
of a Hamiltonian $G$-variety $X$. It turns out that
$\overline{\im\mu_{Z(L)^\circ,X_L}}$ is an affine subspace in
$\z(\lfr)$.  We denote this affine subspace by $\a_{G,X}^{(X_L)}$
and call it the {\it Cartan space} of $X$. It intersects the Lie
algebra of the inefficiency kernel for the action $Z(L)^\circ:X_L$
in the unique point (by the inefficiency kernel of a group action
$\Gamma:Y$ we mean the kernel of the corresponding homomorphism
$\Gamma\rightarrow \Aut(Y)$). Taking this point as the origin in
$\a_{G,X}^{(X_L)}$ we may (and will) consider $\a_{G,X}^{(X_L)}$ as
a vector space.

The group $N_G(L,X_L)$ acts linearly on $\a_{G,X}^{(X_L)}$. We
denote the image of $N_G(L,X_L)$ in $\GL(\a_{G,X}^{(X_L)})$ by
$W_{G,X}^{(X_L)}$ and call it the {\it Weyl group} of $X$. If $G$ is
connected, then $W_{G,X}^{(X_L)}$ is naturally identified with
$N_G(L,X_L)/L$.

Note that, in a suitable sense, the pair
$(\a_{G,X}^{(X_L)},W_{G,X}^{(X_L)})$ does not depend up to
$G$-conjugacy from the choice of $L,X_L$. When a particular choice
of $L,X_L$ does not matter, we write $\a^{(\cdot)}_{G,X}$ for
$\a_{G,X}^{(X_L)}$ and $W^{(\cdot)}_{G,X}$ for $W_{G,X}^{(X_L)}$.

Note that $\im\psi_{L,X_L}\subset\a_{G,X}^{(X_L)}\hookrightarrow
\lfr\quo L$. There is the unique $G$-invariant morphism
$\widehat{\psi}_{G,X}:X\rightarrow
\a_{G,X}^{(\cdot)}/W_{G,X}^{(\cdot)}$ coinciding with
$\psi_{N_G(L,X_L),X_L}$ on $X_L$. The morphism
$\psi_{G,X}:X\rightarrow \g\quo G$ is the composition of
$\widehat{\psi}_{G,X}$ and the finite morphism
$\tau_{G,X}^1:\a_{G,X}^{(\cdot)}/W^{(\cdot)}_{G,X}\rightarrow \g\quo G$ induced
by the embedding $\a_{G,X}^{(X_L)}\hookrightarrow \g$. So
$\widehat{\psi}_{G,X}$ factors through $\widetilde{\psi}_{G,X}$ and
the respective morphism $\tau_{G,X}^2:C_{G,X}\rightarrow
\a_{G,X}^{(\cdot)}/W_{G,X}^{(\cdot)}$ is finite and dominant.

%The following lemma follows directly from the definitions of $\a_{G,X}^{(X_L)}$ and the morphisms
%$\tau_{G,X}^1,\tau_{G,X}^2$.

\begin{Lem}\label{Lem:2.51}
Assume, in addition, that $X$ is conical of degree $k$. Then
$\a_{G,X}^{(X_L)}$ is a vector subspace of $\g$ so one can equip
$\a_{G,X}^{(X_L)}$ with the action of $\C^\times$ given by
$(t,\xi)\mapsto t^k\xi$. Let us equip
$\a_{G,X}^{(\cdot)}/W_{G,X}^{(\cdot)}$ with the induced action. Then
the morphisms $\tau_{G,X}^1,\tau_{G,X}^2$ are
$\C^\times$-equivariant.
\end{Lem}
\begin{proof}
Note that $X_L$ is $\C^\times$-stable and the morphism $\mu_{L,X_L}:X_L\rightarrow \lfr$ is $\C^\times$-equivariant
(here $\C^\times$ acts on $\lfr$ by $(t,\xi)\mapsto t^k\xi$). Now everything follows
directly from the definitions of $\a_{G,X}^{(X_L)}$ and the morphisms
$\tau_{G,X}^1,\tau_{G,X}^2$.
\end{proof}

Now we want to describe the behavior of $\widehat{\psi}_{G,X}$ under
some simple modifications of the pair $(G,X)$. To do this we need to
recall some  results obtained in \cite{Comb_Ham}. The proofs of
these results are mostly straightforward.

Let $X,L,X_L$ be such as above. Let $M$ be a Levi subgroup of $G$
containing $L$, $G^1$ a connected subgroup of $G$ containing $(G,G),
L^1:=G^1\cap L$, $G_1,\ldots,G_k$ be all simple normal subgroups of
$G$, so that $G=Z(G)^\circ G_1\ldots G_k$ is the decomposition into
the locally direct product. Finally, let $X'$ be another affine
irreducible Hamiltonian $G$-variety and $\varphi:X\rightarrow X'$ a
generically finite dominant $G$-equivariant morphism such that
$\mu_{G,X'}\circ\varphi=\mu_{G,X}$.

By Lemma 6.9 from \cite{Comb_Ham},
$\a_{G,X}^{(X_L)}=\a_{G^\circ,X}^{(X_L)}, W_{G^\circ,X}^{(X_L)}$ is
a normal subgroup of $W_{G,X}^{(X_L)}$.

Suppose $G$ is connected. Recall, \cite{Comb_Ham}, Lemmas 4.6,6.10,
that there exists the unique $M$-cross-section $X_M$ of $X$ containing
$X_L$ and $\a_{M,X_M}^{(X_L)}=\a_{G,X}^{(X_L)}, W_{M,X_M}^{(X_L)}=
W_{G,X}^{(X_L)}\cap M/L$.

 By Lemma 4.6 from
\cite{Comb_Ham}, $L$ is the principal centralizer  of $X'$ and there
exists the unique $L$-cross-section $X_L'$ of $X'$ such that
$\varphi(X_L)\subset X_L'$. Further, by Lemma 6.11 from
\cite{Comb_Ham}, $\a_{G,X}^{(X_L)}=\a_{G,X'}^{(X_L')}$,
$W_{G,X}^{(X_L)}\subset W_{G,X'}^{(X_L')}$.

Suppose, in addition, that $0\in \overline{\im\psi_{G,X}}$.  Recall,
\cite{Comb_Ham}, Lemma 4.6, that $L^1$ is the principal centralizer
and $X_L$ is an $L^1$-cross-section of the Hamiltonian $G^1$-variety
$X$. Further, by \cite{Comb_Ham}, Lemma 6.13,
$\a_{G,X}^{(X_L)}\cap\g^1\subset \a_{G^1,X}^{(X_L)}$, the groups
$W_{G^1,X}^{(X_L)},W_{G,X}^{(X_L)}$ are naturally identified, and
the orthogonal projection $\g\twoheadrightarrow \g^1$ induces the
$W_{G,X}^{(X_L)}$-equivariant epimorphism
$\a_{G,X}^{(X_L)}\twoheadrightarrow \a_{G^1,X}^{(X_L)}$.

Finally, suppose $X$ satisfies the equivalent conditions of Lemma
\ref{Lem:2.3.5}. Put $T=L, T_i=L\cap G_i$. Recall, \cite{Comb_Ham},
Lemma 4.6, that $T_i$ is the principal centralizer of the
Hamiltonian $G_i$-variety $X$ and there is the unique
$T_i$-cross-section $X_{T_i}$ of $X$ containing $(\prod_{j\neq
i}G_j)X_T$. Further, Lemma 6.14 from \cite{Comb_Ham} implies that
$\a_{G_i,X}^{(X_{T_i})}=\t_i$, $W_{G,X}^{(X_T)}\subset \prod_{i=1}^k
W_{G_i,X}^{(X_{T_i})}$ and the projection of $W_{G,X}^{(X_T)}$ to
$\GL(\t_i)$ coincides with $W_{G_i,X}^{(X_{T_i})}$.

\begin{Lem}\label{Lem:1.9}
Let $G,X,X_L,M,G^1,L^1,G_1,\ldots,G_k,
X',\varphi,X_M,X_L',T,T_i,X_{T_i}$ be as above.
\begin{enumerate}
\item  $\widehat{\psi}_{G,X}$ is the composition of
$\widehat{\psi}_{G^\circ,X}$ and the natural morphism of quotients
$\a^{(X_L)}_{G^\circ,X}/W^{(X_L)}_{G^\circ,X}\rightarrow
\a_{G,X}^{(X_L)}/W_{G,X}^{(X_L)}$ induced by the inclusion
$W_{G^\circ,X}^{(X_L)}\subset W_{G,X}^{(X_L)}$.
\item Suppose $G$ is connected.   Then the following diagram is
commutative.

\begin{picture}(80,30)
\put(2,2){$X$}\put(2,22){$X_M$}\put(27,2){$\a_{G,X}^{(\cdot)}/W_{G,X}^{(\cdot)}$}
\put(27,22){$\a_{M,X_M}^{(\cdot)}/W_{M,X_M}^{(\cdot)}$}
\put(65,2){$\g\quo G$}\put(65,22){$\m\quo M$}
\put(4,20){\vector(0,-1){14}} \put(34,20){\vector(0,-1){14}}
\put(68,20){\vector(0,-1){14}} \put(6,4){\vector(1,0){20}}
\put(6,24){\vector(1,0){20}} \put(45,4){\vector(1,0){19}}
\put(48,24){\vector(1,0){16}} \put(12,5){\tiny
$\widehat{\psi}_{G,X}$} \put(12,25){\tiny $\widehat{\psi}_{M,X_M}$}
\put(52,6){\tiny $\tau^1_{G,X}$} \put(52,26){\tiny $\tau^1_{M,X_M}$}
\end{picture}

Here the morphism $X_M\rightarrow X$ is the inclusion, the morphism
$\a_{M,X_M}^{(\cdot)}/W_{M,X_M}^{(\cdot)}\rightarrow
\a_{G,X}^{(\cdot)}/W_{G,X}^{(\cdot)}$ is  given by $
W_{M,X_M}^{(X_L)}\xi\mapsto W_{G,X}^{(X_L)}\xi$, and the morphism
$\m\quo M\rightarrow \g\quo G$ is induced by the restriction of
functions from $\g$ to $\m$.
\item The following
diagram is commutative.

\begin{picture}(60,30)
\put(2,2){$\a_{G,X}^{(\cdot)}/W_{G,X}^{(\cdot)}$}\put(9,22){$X$}
\put(42,2){$\a_{G,X'}^{(\cdot)}/W_{G,X'}^{(\cdot)}$}\put(49,22){$X'$}
\put(11,20){\vector(0,-1){13}}\put(51,20){\vector(0,-1){13}}
\put(25,5){\vector(1,0){15}}  \put(15,24){\vector(1,0){32}}
\end{picture}
\item Suppose $G$ is connected and $0\in \overline{\im\psi_{G,X}}$.
 Then the
following
 diagram is commutative.

 \begin{picture}(50,30)
\put(30,20){$X$} \put(5,2){$\a_{G,X}^{(\cdot)}/W_{G,X}^{(\cdot)}$}
\put(42,2){$\a_{G^1,X}^{(\cdot)}/W_{G^1,X}^{(\cdot)}$}
\put(25,4){\vector(1,0){15}}\put(29,19){\vector(-1,-1){12}}
\put(33,19){\vector(1,-1){12}}
\end{picture}
\item Suppose $G$ is connected and $X$ satisfies the equivalent
conditions of Lemma \ref{Lem:2.3.5}.  Then the following diagram,
where the map $\a_{G,X}^{(\cdot)}/W^{(\cdot)}_{G,X}\rightarrow
\a_{G_i,X}^{(\cdot)}/W_{G_i,X}^{(\cdot)}$ is induced by the natural
epimorphism $\g\rightarrow \g_i$, is commutative.

 \begin{picture}(50,30)
\put(30,20){$X$} \put(5,2){$\a_{G,X}^{(\cdot)}/W_{G,X}^{(\cdot)}$}
\put(42,2){$\a_{G_i,X}^{(\cdot)}/W_{G_i,X}$}
\put(25,4){\vector(1,0){15}}\put(29,19){\vector(-1,-1){12}}
\put(33,19){\vector(1,-1){12}}
\end{picture}
\end{enumerate}
\end{Lem}
\begin{proof}
The proofs of assertions 1,3,4 follow directly from the definition
of $\widehat{\psi}_{\bullet,\bullet}$.

Let us prove assertion 2. The commutativity of the right square of
the diagram follows directly from the definition of
$\tau^1_{\bullet,\bullet}$. To prove the commutativity of the left
square we note that  both morphisms $X_M\rightarrow
\a_{G,X}^{(\cdot)}/W_{G,X}^{(\cdot)}$ from the diagram are
$M$-invariant and their restrictions to $X_L$ coincide with
$\widehat{\psi}_{N_G(L,X_L),X_L}$. To complete the proof  it remains
to recall that $MX_L$ is dense in $X_M$.

We proceed to assertion 5. The morphism
$\widehat{\psi}_{G_i,X}|_{X_i}$ is $Z(G)^\circ\prod_{j\neq
i}G_j$-invariant. It follows that $\widehat{\psi}_{G_i,X}$ is
$G$-invariant. It remains to note that the restrictions of both
morphisms $X\rightarrow \a_{G_i,X}^{(\cdot)}/W_{G_i,X}^{(\cdot)}$
coincide on $X_T$.
\end{proof}

Now we are going to quote some properties of
$C_{G,X},\widetilde{\psi}_{G,X},\widehat{\psi}_{G,X}$ proved in
\cite{alg_hamil}.

\begin{Prop}\label{Lem:4.4.1}
The morphism $\widehat{\psi}_{G,X}\quo G:X\quo G\rightarrow
\a_{G,X}^{(\cdot)}/W_{G,X}^{(\cdot)}$ is equidimensional and open.
Further, for any closed subvariety $Y\subset
\im\widehat{\psi}_{G,X}$ and any irreducible component $Z$ of
$(\widehat{\psi}_{G,X}\quo G)^{-1}(Y)$ the subset
$(\widehat{\psi}_{G,X}\quo G)(Z)$ is dense in $Y$.
\end{Prop}
\begin{proof}
Note that $\a_{G,X}^{(\cdot)}/W_{G,X}^{(\cdot)}$ is a normal variety
of dimension $\defe_G(X)$. Thanks to Theorem 1.2.3 from
\cite{alg_hamil}, $\widehat{\psi}_{G,X}\quo G$ is equidimensional.
The openness stems from \cite{Chevalley}, Proposition 3 in Section
5.5. The last assertion of the proposition is an easy corollary of
the fact that $\widehat{\psi}_{G,X}\quo G$ is equidimensional.
\end{proof}

\begin{Prop}[\cite{alg_hamil}, Theorem 1.2.7]\label{Thm:2.2}
 Suppose $X$ is conical. Then $C_{G,X}\cong \a_{G,X}^{(\cdot)}/W_{G,X}^{(\cdot)}$ and
$\widetilde{\psi}_{G,X}=\widehat{\psi}_{G,X}$. Further, the algebra $\C[C_{G,X}]$
coincides with the intersection of $\C[X]$ and the Poisson center of $\C(X)^G$. 
\end{Prop}

\section{Dimensions of fibers}\label{SECTION_dimension}
Throughout the section $G$ is a connected reductive group and $X$ is
a Hamiltonian $G$-variety with  symplectic form $\omega$.

In Subsection \ref{SUBSECTION_affham2} we prove a variant of the
Luna-Richardson restriction theorem (\cite{LR}) for Hamiltonian
varieties.  This allows us to reduce a general affine Hamiltonian
$G$-variety to one satisfying the equivalent conditions of
Lemma~\ref{Lem:2.3.5}.

Subsection~\ref{SUBSECTION_affham4} deals with a stratification of
fibers of the morphism $\psi_{G,X}\quo G:X\quo G\rightarrow \g\quo
G$. A stratum consists of the images of all points with closed
$G$-orbit and the same  determining triple. The main results of the
subsection are the proof that any stratum is smooth and the formula
for the dimensions of the strata (Proposition~\ref{Prop:4.4.2}).

The main part of this section is Subsection
\ref{SUBSECTION_affham5}. There we prove the following result that
strengthens Theorem \ref{Thm:1}.

\begin{Thm}\label{Thm:4.0.1}
The morphisms $\psi_{G,X},\widetilde{\psi}_{G,X},$
$\widehat{\psi}_{G,X}$ are equidimensional. The morphisms
$\widehat{\psi}_{G,X},\widetilde{\psi}_{G,X}$ are open. For any
closed irreducible subvariety  $Y\subset \im\widehat{\psi}_{G,X}$
and any irreducible component $\widetilde{Y}\subset
\widehat{\psi}_{G,X}^{-1}(Y)$ the subvariety
$\pi_{G,X}(\widetilde{Y})\subset X\quo G$ is an irreducible
component of $(\widehat{\psi}_{G,X}\quo G)^{-1}(Y)$.
\end{Thm}

The proof uses the stratification introduced in
Subsection~\ref{SUBSECTION_affham4} and the estimate on  dimensions
of fibers of $\pi_{G,X}$ obtained in Proposition~\ref{Prop:4.5.1}.

\subsection{A Hamiltonian version of the Luna-Richardson theorem}\label{SUBSECTION_affham2}
Let $H$ be a reductive subgroup of $G$. The subvariety $X^H\subset
X$ is smooth (see \cite{VP}, Subsection 6.5) and $N_G(H)$-stable.
Let us equip $X^H$ with a structure of a Hamiltonian
$N_G(H)$-variety.

\begin{Prop}\label{Lem:4.2.1}
\begin{enumerate}
\item $\omega|_{X^H}$ is nondegenerate, thus
$X^H$ is equipped with the symplectic structure.
\item The action
$N_G(H):X^H$ is Hamiltonian with the moment map
$\mu_{N_G(H),X^H}=\mu_{G,X}|_{X^H}$.
\end{enumerate}
\end{Prop}
\begin{proof}
For a symplectic vector space $V$ and a reductive subgroup
$H\subset\Sp(V)$ the $H$-modules $V^{H}$ and $V/(V^H)^{\skewperp}$
are isomorphic. Thus  $\omega|_{V^H}$ is nondegenerate. Since
$T_x(X^H)=(T_xX)^H$, see \cite{VP}, Subsection 6.5, we see that
$\omega|_{X^H}$ is nondegenerate.

Note that the Lie algebra of $N_G(H)$ coincides with $\g^H+\h$.
Since $\mu_{G,X}$ is $G$-equivariant, we have $\mu_{G,X}(X^H)\subset
\g^H$. Clearly, $\mu_{G,X}|_{X^H}$ is  $N_G(H)$-equivariant. It
remains to check that
\begin{equation}\label{eq:4.2:1}v(H_{\xi}|_{X^H})_x=\xi_x\end{equation} for all $\xi\in
\g^H+\h,x\in X^H$. Obviously, $v(H_\xi)_x=\xi_x=0$ for all
$\xi\in\h,x\in X^H$. Thus (\ref{eq:4.2:1}) holds for $\xi\in\h$. Now
let $\xi\in \g^H$. Then $H_\xi\in\C[X]^H$, and $v(H_\xi)_x$ is an
$H$-invariant vector for $x\in X^H$. It follows from the
construction of the symplectic form on $X^H$ that
$v(H_\xi)_x=v(H_\xi|_{X^H})_x$.
\end{proof}

%\begin{Rem}\label{Rem:4.2.2}
%Choose an irreducible (=connected) component $X_0$ of $X^H$ and
%denote by $\widetilde{G}_0$ the stabilizer of $X_0$ in $N_G(H)$.
%Clearly, $H\subset \widetilde{G}_0$. Put $G_0:=\widetilde{G}_0/H$.
%The Lie algebra of $G_0$  is naturally identified with
%$\widetilde{\g}_0\cap\h^{\perp}$.  The action $G_0:X_0$ is
%Hamiltonian with $\mu_{G_0,X_0}=p\circ\mu_{\widetilde{G}_0,X_0}$,
%where $p:\widetilde{\g}_0\twoheadrightarrow \g_0$ is the natural
%epimorphism.
%
% Note that
%\begin{equation}\label{eq:4.43}\defe_{\widetilde{G}_0}(X_0)=\defe_{G_0}(X_0)=\dim\overline{\psi_{G,X}(X_0)}.\end{equation}
%Indeed,  $\defe_{G_0}(X_0)=\dim\overline{\im\psi_{G_0,X_0}}$. Let
%$\tau$ denote the morphism $(\g_0+\xi)\quo G_0\rightarrow \g\quo G$
%induced by the restriction of functions. This morphism is finite
%because it is the composition of the embedding
%$\g_0+\xi_0\hookrightarrow \g^H$ and the morphism $(\g^H)\quo
%\widetilde{G_0}\rightarrow \g\quo G$. The last morphism is finite by
%Luna's theorem~\cite{Luna2}. (\ref{eq:4.43}) follows now from the
%equality $\psi_{G,X}|_{X_0}=\tau\circ\psi_{\widetilde{G}_0,X_0}$.
%\end{Rem}

Now we will apply the previous construction to a special choice of
$H$.

 Let $L$ be the principal
centralizer of $X$ and $X_L$ an $L$-cross-section. By Corollary
4.2.3 from~\cite{alg_hamil}, the restriction of
$\pi_{(L,L),X_L}:X_L\rightarrow X_L\quo (L,L)$ to
$X_L^{(L,L)}\subset X_L$ is an isomorphism. Denote by $L_0$ the unit
component of the inefficiency kernel of the action $L:X_L\quo
(L,L)\cong X_L^{(L,L)}$.
 It follows from Theorem 4.2.1, \cite{alg_hamil}, that
$L_0=(L,L)T_0$, where $T_0$ is the unit component of the
inefficiency kernel for the action $Z(L):X_L$.   Let $X_0$ be the unique connected component of $X^{L_0}$ containing $X_L^{(L,L)}$.
Put $\widetilde{G}_0=N_G(L_0,X_0)$ (the stabilizer of $X_0$ under
the action of $N_G(L_0)$), $ G_0=\widetilde{G}_0/L_0$.  We identify
$\g_0$ with $\g^{L_0}\cap\lfr_0^\perp$. It follows from
Proposition~\ref{Lem:4.2.1} that the action $\widetilde{G}_0:X_0$ is
Hamiltonian with  moment map $\mu_{G,X}|_{X_0}$. By Remark 3.1.2
from \cite{alg_hamil}, the action $G_0:X_0$ is Hamiltonian with the
moment map $\mu_{G_0,X_0}:=p\circ\mu_{\widetilde{G}_0,X_0}$, where
$p$ denotes the natural projection
$\widetilde{\g}_0\rightarrow\g_0$.

The following proposition is what we mean by a "Hamiltonian version
of the Luna-Richardson theorem".

\begin{Prop}\label{Prop:4.2.2}
In the notation introduced above the following statements hold.
\begin{enumerate}
\item
The morphism $X_0\quo G_0\rightarrow X\quo G$ induced by the
restriction of functions is an isomorphism. \item
 $m_{G_0}(X_0)=\dim G_0$,
 $\defe_{G_0}(X_0)=\defe_G(X),\cork_G(X)=\cork_{G_0}(X_0)$.
 \item $L/L_0$ is the principal centralizer of $X_0$. The subvariety $X_L^{L_0}$
 is dense in the unique
 $L/L_0$-cross-section $X_{0L}$ of $X_0$,
$\a_{G_0,X_0}^{(X_{L})}=\a_{G,X}^{(X_{0L})}-\xi_0$, where $\xi_0\in
\lfr_0\cap\a_{G,X}^{(X_{L})},$ and
$W_{G_0,X_0}^{(X_{0L})}=W_{G,X}^{(X_{L})}$.
\item  $\widehat{\psi}_{G,X}|_{X_0}=\widehat{\psi}_{G_0,X_0}$.
\end{enumerate}  \end{Prop}

In the proof we will use some notions of the theory of algebraic
transformation groups.  Let $Y$ be an irreducible affine variety
acted on by a reductive group $H$. It is known, see \cite{VP},
Theorem 7.12, that there  exists an open subset $Y_0\subset Y\quo H$
such that for any $y\in Y_0$ the closed orbit in $\pi_{H,Y}^{-1}(y)$
is isomorphic to $H/C$, where $C$ is a reductive subgroup of $H$.

\begin{defi}\label{defi:4.2.3} Such a subgroup $C$ (determined uniquely up to $H$-conjugacy)
is called the {\it principal isotropy subgroup} for the action
$H:Y$.
\end{defi}

The action $H:Y$ is called {\it stable} if its general orbit is
closed and {\it locally free} if $m_H(Y)=\dim H$.

\begin{proof}[Proof of Proposition \ref{Prop:4.2.2}]
The action $Z(L)^\circ:X_L^{(L,L)}\cong X_L\quo (L,L)$ is stable
(\cite{alg_hamil}, Proposition 4.5.1). Thus $L_0$ is the unit
component of the principal isotropy subgroup for the action $L:X_L$.
Since the natural morphism $G*_LX_L\rightarrow X$ is \'{e}tale and
its image is saturated, we see that the group $L_0$ is the unit
component of the principal isotropy subgroup for the action $G:X$
and that the morphism $X_0\quo G_0\rightarrow X\quo G$ is dominant.
By the Luna-Richardson theorem (\cite{LR}), the morphism $X_0\quo
G_0\rightarrow X\quo G$ is an isomorphism and the action of $G_0$ on
$X_0$ is locally free.
 The latter yields $\defe_{G_0}(X_0)=\rank
G_0=\rank G-\rank L_0=\defe_G(X)$. By Theorem 1.2.9 from
\cite{alg_hamil}, $\C(X)^G=\Quot(\C[X]^G),
\C(X_0)^{G_0}=\Quot(\C[X_0]^{G_0})$. So
$$\cork_G(X)=\td \C(X)^G-\defe_G(X)=\td \C(X_0)^{G_0}-\defe_{G_0}(X_0)=\cork_{G_0}(X_0).$$

%REFEREE'S REMARK: more detailed proof
We proceed to assertion 3. Since $m_{G_0}(X_0)=\dim G_0$, the
maximal torus $L/L_0\subset G_0$ is the principal centralizer of
$X_0$ (see Lemma \ref{Lem:2.3.5}) and
$\a_{G_0,X_0}^{(X_{0L})}=\lfr\cap\lfr_0^{\perp}=\a_{G,X}^{(X_L)}-\xi_0$
for any $L/L_0$-cross-section $X_{0L}$ of $X_0$. The natural
morphism $X_L\quo L\rightarrow X\quo G$ is dominant and quasifinite,
therefore so is the natural morphism $(X_L^{L_0})\quo
(L/L_0)\rightarrow X_0\quo G_0$.  It follows from \cite{alg_hamil},
Theorem 1.2.9, that the actions $L/L_0:X_L^{L_0},G_0:X_0$ are
stable. It follows that $\dim X_L^{L_0}=\dim (X_L^{L_0})\quo
(L/L_0)+\dim L/L_0=\dim X_0\quo G_0+\dim L/L_0=\dim X_0-\dim
G_0+\dim L/L_0$. Since
$\mu_{G,X}(X_L^{L_0})\subset\mu_{G,X}(X_L)\cap \mathfrak{g}_0\subset \mathfrak{l}^{pr}\cap\mathfrak{g}_0
\subset (\mathfrak{l}/\mathfrak{l}_0)^{pr}$ (the last subset is taken w.r.t. the Lie algebra
$\g_0$), we see that $X_L^{L_0}$ lies in the unique
$L/L_0$-cross-section $X_{0L}$ of $X_0$. Comparing the dimensions,
we see that $X_L^{L_0}$ is dense in $X_{0L}$. The equality for the
Weyl groups stems from $N_G(L,X_L)/L_0\subset G_0,
N_G(L,X_L^{L_0})=N_G(L,X_L)$.

Finally,  both morphisms in assertion 4 are $G_0$-invariant and
their restrictions to $X_L^{L_0}$ are equal to the restriction of
$\widehat{\psi}_{N_G(L,X_L),X_L}$.
\end{proof}

\subsection{A stratification of a fiber of $\psi_{G,X}\quo G$}\label{SUBSECTION_affham4}

In this subsection we introduce a stratification of  fibers of the
morphism $\psi_{G,X}\quo G:X\quo G\rightarrow \g\quo G$. We consider
fibers of $\psi_{G,X}\quo G$ as algebraic varieties. Namely, let
$\eta\in \g$, $H$ be a reductive subgroup of $G_\eta$ and $V$ a
symplectic $H$-module. We put
\begin{equation*}S_{G,X}(H,\eta,V)\index{sgxfev@$S_{G,X}(H,\eta,V)$}=\{\pi_{G,X}(x)| Gx \text{ is
closed},(H,\eta,V) \text{ is the determining triple of }X\text{ at
}x\}.\end{equation*} Clearly,
$S_{G,X}(H_1,\eta_1,V_1)=S_{G,X}(H_2,\eta_2,V_2)$ iff  there is
$g\in G$ and a linear isomorphism $\iota:V_1\rightarrow V_2$ such
that $\Ad(g)\eta_1=\eta_2$, $gH_1g^{-1}=H_2$ and
$(ghg^{-1})\iota(v)=\iota(hv)$ for all $h\in H_1$.

The main result of this subsection is the following
\begin{Prop}\label{Prop:4.4.2}
Let $X,G, H,\eta,V$ be as above, $\lambda=\pi_{G,\g}(\eta)$. Then
$S_{G,X}(H,\eta,V)$ is a locally-closed smooth  subvariety  of pure
codimension $\cork_G(X)-\dim V^H$ in $(\psi_{G,X}\quo
G)^{-1}(\lambda)$.
\end{Prop}
\begin{proof}
Firstly, we show that $S_{G,X}(H,\eta,V)$ is a locally-closed
subvariety of $X\quo G$. Denote by $Y$ the set of all points $x\in
X$ such that $Gx$ is closed, $G_x=H$, and  $T_xX/\g_*x\cong V\oplus
(\g_\eta/\h)^*$. It follows from the Luna slice theorem applied to
any point of $Y$ that $Y$ is a locally-closed subvariety in $X$.
Therefore $Y_\eta=Y\cap \mu_{G,X}^{-1}(\Ad(G)\eta)$ is a locally
closed subvariety of $X$. Since all orbits in $Y_\eta$ are closed in
$X$, we see that $Y_{\eta}$ is an open saturated subvariety of
$\overline{Y_\eta}$. Thus $S_{G,X}(H,\eta,V)=\pi_{G,X}(Y_\eta)$ is
open in $\overline{Y_\eta}\quo G$.

Applying Proposition~\ref{Prop:4.3.3}, we reduce the codimension and
smoothness claims to the case $X=M_{G}(H,\eta,V)$. Put
$\s=\z_\g(\eta_s)$. Choose an $\sl_2$-triple $(\eta_n,h,f)$ in
$\s^H$ generating $M_G(H,\eta,V)$. Denote by $U$ the $H$-module
$\z_\s(f)\cap\h^\perp$.

\begin{Lem}\label{Lem:4.4.3}
In the above notation  $\eta$ is an isolated point of
$(\eta+\z_\s(f))\cap \overline{\Ad(G)\eta}$.
\end{Lem}
\begin{proof}[Proof of Lemma~\ref{Lem:4.4.3}]
Note that $T_\eta(\eta+\z_\s(f))=\z_\s(f),
T_\eta\overline{\Ad(G)\eta}=[\g,\eta]$. It is enough to show
$\z_\s(f)\cap [\g,\eta]=\{0\}$. The equality $\s=\z_\g(\eta_s)$
yields $[\g,\eta]=[\s^{\perp},\eta]+ [\s,\eta]=\s^{\perp}\oplus
[\s,\eta_n]$. Thanks to the representation theory of $\sl_2$,
$[\s,\eta_n]\cap \z_\s(f)=0$ whence the required equality.
\end{proof}

In virtue of Remark~\ref{Rem:2.1.11}, it is enough to assume that
$V^H=\{0\}$. Put $x:=[1,(0,0)]$. Everything will follow if we check
that $\pi_{G,X}(x)$ is an isolated point in $S_{G,X}(H,\eta,V)$.
Indeed, by Proposition~\ref{Lem:4.4.1}, $\cork_G(X)=\dim X\quo
G-\defe_G(X)=\dim_{\pi_{G,X}(x)}(\psi_{G,X}\quo G)^{-1}(\lambda)$.

 There exists a neighborhood
$O'$ of $\eta$ in $\eta+\z_\s(f)$ such that $O'\cap
\overline{\Ad(G)\eta}=\eta$. Replacing $O'$ with $HO'$, if
necessary, we may assume that $O'$ is $H$-stable. Set
$O:=\{[g,(u,v)]\in M_G(H,\eta,V)| \eta+u+\mu_{H,V}(v)\in O'\}$. By
definition, $O$ is an open $G$-subvariety of  $X$ containing $x$. It
is enough to show that any point $x_1\in O$ with closed $G$-orbit
and the determining triple
 $(H,\eta,V)$ is $G$-conjugate to $x$. Assume the
converse. Put $x_1=[g,(u,v)]$, $u\in U, v\in V, (u,v)\neq 0$. Recall
that $\mu_{G,X}(x_1)=\Ad(g)(\eta+u+\mu_{H,V}(v))$. Since
$\mu_{G,X}(x_1)=\eta$, Lemma~\ref{Lem:4.4.3} implies that
$u+\mu_{H,V}(v)=0$. Since $U\cap \h=\{0\}$, we have $u=0$. The
subgroup $H_v\subset H$ is conjugate to $H$ in $G$. Thus $v\in
V^H=\{0\}$. Contradiction.
\end{proof}

\subsection{The proof of Theorem~\ref{Thm:4.0.1}}\label{SUBSECTION_affham5}
At first, we obtain an estimate for the dimension of a fiber of
$\pi_{G,X}$.

\begin{Prop}\label{Prop:4.5.1}
The dimension of any fiber of $\pi_{G,X}:X\rightarrow X\quo G$ does
not exceed $\dim X-\defe_G(X)-\frac{\cork_G(X)}{2}$.
\end{Prop}
\begin{proof}
The proof is carried out in two steps. Firstly, we consider the case
when $X$ satisfies the equivalent conditions of
Lemma~\ref{Lem:2.3.5} and then deduce the general case from this
one.

{\it Step 1. } Suppose $X$ satisfies the equivalent conditions of
Lemma~\ref{Lem:2.3.5}. Then
$$\defe_G(X)+\frac{\cork_G(X)}{2}=\frac{\dim
X-\dim G+\rank G}{2}.$$

Let $y\in X\quo G$,  $x$ be a point from the unique closed $G$-orbit
in $\pi_{G,X}^{-1}(y)$, $H=G_x$, $\eta=\mu_{G,X}(x)$,
$U=(\z_\g(\eta)/\h)^*, V=(\g_*x)^{\skewperp}/(\g_*x\cap
(\g_*x)^{\skewperp})$. The $H$-modules $U\oplus V$ and $T_xX/\g_*x$
are isomorphic.

Using the Luna slice theorem, we see that it is enough to check
\begin{equation}\label{eq:4.5:0}\dim\pi_{H,U\oplus V}^{-1}(0)\leqslant \dim U+\dim V-\frac{\dim
X-\dim G+\rank G}{2}\end{equation}

\begin{Lem}[\cite{Schwarz}, Proposition 2.10]\label{Lem:4.5.2}
Let  $H$ be a reductive group, $T_H$ a maximal torus of $H$, and
$V$ a self-dual $H$-module. Then
\begin{equation*}%\label{eq_dim_null}
\dim\pi_{H,V}^{-1}(0)\leqslant \frac{1}{2}(\dim V-\dim V^{T_H}+\dim
H-\dim T_H).
\end{equation*}
\end{Lem}

\begin{Lem}\label{Lem:4.5.3}
 $U\oplus V$ is a self-dual $H$-module.
\end{Lem}
\begin{proof}[Proof of Lemma~\ref{Lem:4.5.3}]
Note that the $H$-modules $U\oplus V$ and $T_xX/\g_*x$ are
isomorphic. The module $T_xX$ is symplectic, while the module
$\g_*x\cong \g/\h\cong\h^\perp$ is orthogonal. Hence both these
modules are self-dual. Therefore the quotient module $U\oplus V$ is
self-dual too.
%Replacing $\g$ with $\z_\g(\eta_s)$ if necessary, we may  assume
%that $\eta_s\in\z(\g)$ and easily reduce the proof to the case when
%$\eta$ is nilpotent.  The representation $H:V$ is symplectic and
%thus self-dual.  Choose an $\sl_2$-triple $(\eta,h,f)\subset \g^H$
%and consider the grading $\g=\bigoplus_{i\in \Z}\g_i$, where
%$\g_i:=\{\xi\in\g| [h,\xi]=i \xi \}$. Put
%$\z_\g(\eta)_i:=\z_\g(\eta)\cap \g_i$. It is well known that
%$\z_\g(\eta)=\bigoplus_{i\geqslant 0} \z_\g(\eta)_i$ and that
%$\z_\g(\eta)_0$ is a reductive subalgebra of $\g$. The latter
%implies that  $\z_\g(\eta)_0/\h\subset U$ is an orthogonal
%$H$-module. Now it is enough to show that the $H$-module
%$\z_\g(\eta)_i$ is self-dual.
%
%Choose a connected subgroup $S\subset G$ corresponding to the
%$\sl_2$-triple $\langle \eta_n,h,f\rangle$, and let $w$ be an
%element of $N_S(\langle h\rangle)$ acting on $\langle h\rangle$ by
%$-\id$. For $x,y\in \z_\g(\eta)_i$ put $(x,y)_w=(x,wy)$. Since $w\in
%Z_G(H)$, the form $(\cdot,\cdot)_w$ is $H$-invariant. It follows
%from the representation theory of $\sl_2$ that
%$$ [\g,\eta]\oplus \z_\g(f)=\g,w(\z_{\g}(\eta)_i)=\z_\g(w(\eta))\cap \g_{-i}=\z_\g(f)\cap
%\g_{-i}.$$
%
%Note that $\z_\g(\eta)^\perp= [\g,\eta], (\g_i,\g_j)=0 \text{ for
%}i+j\neq 0$. Thence the pairing of $\z_\g(\eta)_i$ and $\z_\g(f)\cap
%\g_{-i}$ induced by $(\cdot,\cdot)$ is nondegenerate. Equivalently,
% $(\cdot,\cdot)_w$ is nondegenerate on $\z_\g(\eta)_i$.
\end{proof}

We see that the $H$-module $U\oplus V$ satisfies the assumptions of
Lemma~\ref{Lem:4.5.2}. Let $T_H$ be a maximal torus of $H$. Let us
show that $\dim U^{T_H}\geqslant \rank\g-\rank\h$. Since $\dim
\h^{T_H}=\rank\h$, it is enough to show that
$\dim\z_\g(\xi)^{T_H}\geqslant \rank\g$ for any $\xi\in\g^H$. It is
enough to check the last inequality for $\xi\in \g^H$ in general
position. But in this case $\xi$ is semisimple. Thence $\z_\g(\xi)$
is a Levi subalgebra of $\g$ and everything is clear.

By Lemma~\ref{Lem:4.5.2}, we have the following inequalities
\begin{equation}\label{eq:4.5:1}
\begin{split}
&\dim\pi_{H,U\oplus V}^{-1}(0)\leqslant \frac{1}{2}(\dim U+\dim
V-\dim U^{T_H}-\dim V^{T_H}+\dim \h-\rank\h)\\ &\leqslant
\frac{1}{2}(\dim U+\dim V-(\rank\g-\rank\h)+\dim\h-\rank\h).
\end{split}
\end{equation}

One may check directly that the last expression in (\ref{eq:4.5:1})
coincides with the r.h.s of (\ref{eq:4.5:0}).

{\it Step 2.} Now we consider the general case. Let $X_0,G_0$ be as
in Subsection~\ref{SUBSECTION_affham2}.

By  Proposition~\ref{Prop:4.2.2}, $\cork_G(X)=\cork_{G_0}(X_0),
\defe_G(X)=\defe_{G_0}(X_0)$. The proposition will follow if we show
that
\begin{equation}\label{eq:4.5:11}\codim_X\pi_{G,X}^{-1}(y)\geqslant
\codim_{X_0}\pi_{G_0,X_0}^{-1}(y), \end{equation} for any $y\in
X\quo G$. It follows from Proposition~\ref{Prop:4.2.2} that $X_0\quo
G_0\cong X\quo G$, $\pi_{G_0,X_0}^{-1}(y)=\pi_{G,X}^{-1}(y)\cap
X_0$. Now (\ref{eq:4.5:11}) stems from the following general fact of
Algebraic geometry:

 $\dim_xY\cap Z\geqslant \dim_x Y+\dim_x Z-\dim X$ for any  subvarieties $Y,Z$ of an irreducible  variety $X$
 and $x\in Y\cap Z$ provided $X$ is smooth.
\end{proof}

\begin{proof}[Proof of Theorem~\ref{Thm:4.0.1}]
Clearly, $\widetilde{\psi}_{G,X},\psi_{G,X}$ are equidimensional
provided $\widehat{\psi}_{G,X}$ is. As we mentioned above, any
equidimensional morphism to a normal variety is open.

To prove the theorem it  remains to check that for all
$\lambda\in\g\quo G$ and any irreducible component $Z$ of
$\psi_{G,X}^{-1}(\lambda)$ the equality $\dim\pi_{G,X}(Z)=\dim X\quo
G-\defe_G(X)$ and the inequality $\dim Z\leqslant \dim X-\defe_G(X)$
take place (the opposite inequality holds automatically, since
$\defe_G(X)=\dim\overline{\im\psi_{G,X}}$). The former equality will
imply
\begin{equation}\label{eq:4.5:2}\dim\pi_{G,X}(Z)=\dim
X\quo G-\defe_G(X)+\dim Y\end{equation}
 for an irreducible component $Z$ of
 $\widehat{\psi}_{G,X}^{-1}(Y)$, where $Y\subset
\im\widehat{\psi}_{G,X}$ is an arbitrary closed irreducible
subvariety  (recall that, by Proposition~\ref{Lem:4.4.1},
$\im\widehat{\psi}_{G,X}=\im(\widehat{\psi}_{G,X}\quo G)$ is an open
subvariety in $\a_{G,X}^{(\cdot)}/W_{G,X}^{(\cdot)}$). Thanks to
Proposition \ref{Lem:4.4.1},  (\ref{eq:4.5:2}) holds iff
$\pi_{G,X}(Z)$ is an irreducible component in
$(\widehat{\psi}_{G,X}\quo G)^{-1}(Y)$.

Choose a subvariety $S_{G,X}(H,\eta,V)\subset (\psi_{G,X}\quo
G)^{-1}(\lambda)$ (see Subsection~\ref{SUBSECTION_affham4}) such
that $\pi_{G,X}(Z)\cap S_{G,X}(H,\eta,V)$ is dense (and so, in
virtue of Proposition~\ref{Prop:4.4.2}, open) in $\pi_{G,X}(Z)$.
Further, choose a point  $x\in Z\cap
\pi_{G,X}^{-1}(S_{G,X}(H,\eta,V))$ with closed $G$-orbit. Applying
Proposition~\ref{Prop:4.3.3} to $x$, we may replace $X$ with
$M_G(H,\eta,V)$. Thanks to Remark~\ref{Rem:2.1.11}, we may assume
that $V^H=0$. From Proposition~\ref{Prop:4.4.2} it follows that
$\pi_{G,X}(Z)$ is a point. By Proposition~\ref{Prop:4.5.1}, $\dim
Z\leqslant \dim X-\defe_G(X)-\frac{1}{2}\cork_G(X)$. It follows that
$\cork_G(X)=0, \dim (\psi_{G,X}\quo G)^{-1}(\lambda)=0, \dim Z=\dim
X-\defe_G(X)$. This verifies the claim in the beginning of the
previous paragraph and completes the proof.
\end{proof}

\begin{Cor}\label{Cor:4.5.2}
For any  $\lambda\in \im\psi_{G,X}$ and any irreducible component
$Z$ of $\psi_{G,X}^{-1}(\lambda)$ there exists an open subset
$Z_0\subset Z\quo G$ such that $Z_0$ is smooth (as a variety),
$\codim_{Z\quo G }(Z\quo G)\setminus Z_0\geqslant 2$, and for any
$z\in Z_0$ and any point $x\in\pi_{G,X}^{-1}(z)$ with closed
$G$-orbit the following condition holds:
\begin{enumerate}
\item[(*)]  $M_G(H,\eta,V/V^H)$ is coisotropic, where $(H,\eta,V)$ is
the determining triple of $X$ at $x$.
\end{enumerate}
Moreover, $M_G(H,\eta,V/V^H)$ does not depend (up to an isomorphism)
on the choice of $z$.
\end{Cor}
\begin{proof}
 (*) is equivalent to $\cork_G(X)=\cork_G(M_G(H,\eta,V))=\dim
V^H$. It follows from Theorem~\ref{Thm:4.0.1} that  $Z$ maps
dominantly  whence, by the standard properties of quotient
morphisms, surjectively onto some irreducible component of
$(\psi_{G,X}\quo G)^{-1}(\lambda)$. The required claims follow now
from Proposition~\ref{Prop:4.4.2}.
\end{proof}

\begin{Cor}\label{Cor:4.5.3}
Let $Y$ be a closed irreducible subvariety in
$\im\widehat{\psi}_{G,X}$. Then
$\overline{\widehat{\psi}_{G,X}(\widetilde{Y})}=Y$ for any
irreducible component $\widetilde{Y}$ of
$\widehat{\psi}_{G,X}^{-1}(Y)$.
\end{Cor}
\begin{proof}
According to Theorem~\ref{Thm:4.0.1}, $\pi_{G,X}(\widetilde{Y})$ is
an irreducible component of $(\widehat{\psi}_{G,X}\quo
G)^{-1}(Y)\subset X\quo G$. It remains to apply
Proposition~\ref{Lem:4.4.1}.
\end{proof}

\begin{Cor}\label{Cor:4.5.4}
A simply connected affine conical Hamiltonian $G$-variety satisfies
(Utw1).
\end{Cor}
\begin{proof}
Thanks to Proposition \ref{Thm:2.2},
$\tau^2_{G,X}:C_{G,X}\rightarrow
\a_{G,X}^{(\cdot)}/W_{G,X}^{(\cdot)}$ is an isomorphism. By
Theorem~\ref{Thm:4.0.1}, the morphism
$\widetilde{\psi}_{G,X}:X\rightarrow C_{G,X}$ is equidimensional.
Since $G$ is connected, the subalgebra $\C[X]^G$ is integrally
closed in $\C[X]$. Thus $\C[C_{G,X}]$ is integrally closed in
$\C[X]$. In other words, a general fiber of $\widetilde{\psi}_{G,X}$
is connected. Summarizing, we see that $\widetilde{\psi}_{G,X}$ is
an equidimensional morphism with a connected general fiber from a
simply connected variety $X$ to $C_{G,X}\cong
\a_{G,X}^{(\cdot)}/W_{G,X}^{(\cdot)}$.  The proof of the proposition
is based on an idea of Panyushev~\cite{Panyushev} and is completely
analogous to that given in~\cite{Knop6}, Theorem 7.2.
\end{proof}

\section{Some results concerning Weyl groups}\label{SECTION_Weyl}
Throughout the section $G,X,\omega$ have the same meaning as in the
previous section.

In this section we study the structure of the Weyl group
$W_{G,X}^{(\cdot)}$. Subsection \ref{SUBSECTION_affham6} contains
three technical propositions, which play a crucial role in the
subsequent exposition. Propositions \ref{Prop:4.6.1},
\ref{Prop:4.6.5} allow one to  reduce the study of  an arbitrary
affine Hamiltonian $G$-variety to the study of a coisotropic conical
model variety. Proposition \ref{Prop:4.6.3} describes the behavior
of Weyl groups under this reduction.

Using results of Subsection \ref{SUBSECTION_affham6}, in Subsection
\ref{SUBSECTION_Weyl_aff1} we establish some properties of Weyl
groups of varieties satisfying the equivalent conditions of Lemma
\ref{Lem:2.3.5}. In particular, we get some restrictions on
varieties with a "small" Weyl group (Proposition \ref{Prop:5.2.1},
Corollary \ref{Cor:5.2.3}) and show that a Weyl group cannot be "too
small" (Corollary \ref{Cor:5.2.5}). As a consequence of Corollary
\ref{Cor:5.2.5} we get  some explicit restrictions on Weyl groups
for simple $G$ of types $A-E$ in Proposition \ref{Prop:5.2.6},
Corollary \ref{Cor:5.2.7}.

Finally, in Subsection \ref{SUBSECTION_Weyl_computation} we compute
the Weyl groups of  linear actions of simple groups satisfying some
additional restrictions. This computation will be used in Subsection
\ref{SUBSECTION_Utw3} to check that any symplectic $G$-module is an
untwisted Hamiltonian variety.

\subsection{Some technical propositions}\label{SUBSECTION_affham6}
\begin{Prop}\label{Prop:4.6.1}
 Let $L$ be the principal centralizer and $X_L$ an
$L$-cross-section of $X$, $\xi\in\a_{G,X}^{(X_L)}$,
$\alpha=\pi_{W_{G,X}^{(X_L)},\a_{G,X}^{(X_L)}}(\xi)$, $M=Z_G(\xi)$.
Suppose $\alpha\in \im\widehat{\psi}_{G,X}$. Choose an irreducible
component $Z$ of $\widehat{\psi}_{G,X}^{-1}(\alpha)$. Then there
exists $x\in X$ possessing the following properties:
\begin{itemize}
\item[(a)] $x\in Z$.
\item[(b)] $\mu_{G,X}(x)_s\in \z(\m)\cap\m^{pr}$.
\item[(c)] A unique $M$-cross-section $X_M$ of $X$
containing  $x$ contains  $X_L$ and
$\widehat{\psi}_{M,X_M}(x)=\pi_{W_{M,X_M}^{(X_L)},\a_{M,X_M}^{(X_L)}}(\xi)$.
\item[(d)] $Gx$ is closed in $X$.
\item[(e)]  Let $(H,\eta,V)$ be the determining triple of $X$ (or, equivalently, of $X_M$) at $x$ and
$\widehat{G}$ be a connected subgroup of $M$ containing
$(M,M)H^\circ$. The orbit $\widehat{G}x$ is closed in $X_M$ and the
Hamiltonian $\widehat{G}$-variety
$\widehat{X}:=M_{\widehat{G}}(H\cap \widehat{G},\eta_n,V/V^H)$
 is coisotropic.
\end{itemize}
\end{Prop}
\begin{Rem}\label{Rem:4.6.2}
If $X$ satisfies the equivalent conditions of Lemma \ref{Lem:2.3.5},
then so does the Hamiltonian $\widehat{G}$-variety $\widehat{X}$.
This stems easily from Proposition \ref{Prop:4.3.3}.
\end{Rem}
\begin{proof}[Proof of Proposition \ref{Prop:4.6.1}]
Choose a point $z\in Z$ with closed $G$-orbit. Let us show that $gz$
satisfies (b),(c) for some $g\in G$.  Put $M_1=Z_G(\mu_{G,X}(z)_s)$.
Since $\pi_{G,\g}(\mu_{G,X}(z)_s)= \pi_{G,\g}(\xi)$, we have
$M_1\sim_G M$. Let $X_{M_1}$ be an $M_1$-cross-section of $X$
containing $z$, $L_1$ be the principal centralizer and $X_{L_1}$ an
$L_1$-cross-section of $X_{M_1}$. Replacing $z$ with $gz$ for an
appropriate element $g\in G$, we may assume that $L_1=L,
X_{L_1}=X_L$. Next, replacing $z$ with $mz$ for some $m\in M_1$, one
obtains
$\mu_{G,X}(z)_s\in\a_{M_1,X_{M_1}}^{(X_L)}=\a_{G,X}^{(X_L)}$. By the
commutative diagram of assertion 2 of Lemma \ref{Lem:1.9}, for some
$n\in N_G(L,X_L)$ the following equality holds
\begin{equation}\label{eq:4.6:2}
\widehat{\psi}_{M_1,X_{M_1}}(z)=\pi_{W_{M_1,X_{M_1}}^{(X_L)},\a_{G,X}^{(X_L)}}(n\xi).
\end{equation}
 Note that $\psi_{M_1,X_{M_1}}(z)\in
\z(\m_1)\hookrightarrow \m_1\quo M_1$. From (\ref{eq:4.6:2}) it
follows that $\pi_{M_1,\m_1}(n\xi)\in\z(\m_1)\hookrightarrow
\m_1\quo M_1$ whence $n\xi\in\z(\m_1)$. On the other hand,  $n\xi\in
\z(\Ad(n)\m)\cap (\Ad(n)\m)^{pr}$ and so $\m_1\subset \Ad(n)\m$. We
have seen above that $M_1\sim_G M$ whence $M_1=nMn^{-1}$. Replacing
$z$ with $n^{-1}z$, we get the point $z$ satisfying (a)-(c). Put
$\alpha'=\pi_{W_{M,X_M}^{(X_L)},\a_{M,X_M}^{(X_L)}}(\xi)$.

 According to Lemma~\ref{Lem:4.3.4}, there  exists an
open affine $M$-saturated subvariety $X_M^0\subset X_M$ containing
$z$ such that  for any $x\in X_M^0$ the orbit $Gx$ is closed in $X$
iff $Mx$ is closed in $X_M$. Further, by Lemma~\ref{Lem:4.3.5},
$\widehat{G}x\subset X_M$ is closed whenever  $Mx$ is closed.

From assertion 2 of Lemma \ref{Lem:1.9}, Theorem \ref{Thm:4.0.1} and
the fact that the natural morphism $G*_{M}X_M\rightarrow X$ is
\'{e}tale we get $\dim Z\cap
X_M=\dim\widehat{\psi}_{M,X_M}^{-1}(\alpha')$. Hence there is an
irreducible component $Z'$ of $\psi_{M,X_M}^{-1}(\alpha')$
containing $z$ and contained in $Z\cap X_M$. By Corollary
\ref{Cor:4.5.2}, there is an open subset $Y^0\subset
\pi_{M,X_M}(Z')$ such that any point $x\in \pi_{M,X_M}^{-1}(Y^0)$
with closed $M$-orbit satisfies (a)-(d) and (e) for $\widehat{G}=M$.
When $\widehat{G}\neq M$, there is a covering
$T_0\times\widehat{G}\twoheadrightarrow M$ and a finite Hamiltonian
morphism $T^*(T_0)\times \widehat{X}\rightarrow M_M(H\cap
\widehat{G},\eta_n,V/V^H)$,  where $T_0$ is a torus. Since
$H^\circ\subset \widehat{G}$, we are done.
\end{proof}

\begin{Prop}\label{Prop:4.6.5}
Let $X,L,X_L$ be as in Proposition \ref{Prop:4.6.1}, $T_0$ denote
the unit component of the inefficiency kernel of the action
$Z(L)^\circ:X_L$, $\xi_0\in\a_{G,X}^{(X_L)}$, $M=Z_G(\xi_0)$.
Suppose $0\in\im\widehat{\psi}_{G,X}$. Put
$\z:=\z(\m)\cap\a_{G,X}^{(X_L)},
\underline{Z}:=\pi_{W_{G,X}^{(X_L)},\a_{G,X}^{(X_L)}}(\z)$. Choose
an irreducible component $\widetilde{Z}$ of
$\widehat{\psi}_{G,X}^{-1}(\underline{Z})$. Let $\xi\in \z$ be a
point in general position. Then there is a component $Z$ of
$\widehat{\psi}_{G,X}^{-1}(\pi_{W_{G,X}^{(X_L)},\a_{G,X}^{(X_L)}}(\xi))$
lying in $\widetilde{Z}$ and a point $x\in Z$ satisfying the
conditions (b)-(e) of Proposition \ref{Prop:4.6.1} and
\begin{itemize}
\item[(f)] $G_x^\circ\subset (M,M)T_0$.
\end{itemize}
\end{Prop}

\begin{Rem}\label{Rem:4.6.4}
Under the assumptions of Proposition \ref{Prop:4.6.5}  one may
assume that  $\widehat{G}$ defined in (d) coincides with $(M,M)T_0$.
If $X$ satisfies the equivalent conditions of Lemma~\ref{Lem:2.3.5},
then one can take $(M,M)$ for $\widehat{G}$.
\end{Rem}

\begin{proof}[Proof of Proposition \ref{Prop:4.6.5}]
The morphism $\widehat{\psi}_{G,X}$ is open, Theorem
\ref{Thm:4.0.1}.  So $Z,\widetilde{Z}$ do exist. Choose a point
$z\in Z$ satisfying conditions (a)-(e) and such that
$\pi_{G,X}(\widetilde{Z})$ is the only component of
$(\widehat{\psi}_{G,X}\quo G)^{-1}(\underline{Z})$ (see Theorem
\ref{Thm:4.0.1}) containing $\pi_{G,X}(z)$.

Let $X_M$ be as in (c).  By the choice of $z$, any irreducible
component $\widetilde{Z}'$ of $\psi_{M,X_M}^{-1}(\z)$ containing $z$
is contained in $\widetilde{Z}\cap X_M$, compare with the proof of
Proposition \ref{Prop:4.6.1}. As in that proof, there is an open
subset $Y^0\subset\pi_{M,X_M}(\widetilde{Z}')$ such that any $x\in
\pi_{M,X_M}^{-1}(Y^0)$ with closed $M$-orbit satisfies conditions
(a)-(e) (for
appropriate $\xi$). %Since $M\mu_{M,X_M}(x)_s\cap
%\a_{G,X}^{(X_L)}\neq\varnothing$ for any $x\in X_M$, we see that
%$\widetilde{Z}'$ is automatically an irreducible component of
%$\widehat{\psi}_{M,X_M}^{-1}(\z(\m))$.
%The following lemma completes the proof.

 It remains to prove that
$M_x^\circ\subset (M,M)T_0$ for a general point $x\in
\widetilde{Z}'$ with closed $M$-orbit. Recall (see the discussion
preceding Proposition \ref{Prop:4.2.2}) that $L_0:=(L,L)T_0$ is the
unit component of the principal isotropy group for the action
$M:X_M$.
 Let $C$ denote the principal isotropy
subgroup for the action $M:\widetilde{Z}'$, so $L_0\subset C$.  By
the definition of $C$, there exists an irreducible component $X_1$
of $X_M^C$ such that $\pi_{M,X_M}(X_1\cap \widetilde{Z}')$ is dense
in $\pi_{M,X_M}(\widetilde{Z}')$.

By Lemma~\ref{Lem:4.2.1}, the action $N_M(C,X_1):X_1$ is Hamiltonian
with  moment map $\mu_{N_M(C,X_1),X_1}=\mu_{M,X_M}|_{X_1}$. Since
$0\in \overline{\im\psi_{G,X}}$, we get $0\in
\overline{\psi_{M,X_M}(\widetilde{Z}')}$, equivalently,
$\overline{\mu_{M,X_M}(X_1)}$ contains a nilpotent element. Since
$C$ acts trivially on $X_1$, we get
\begin{equation}\label{eq:4.6:3}\mu_{M,X_M}(X_1)\subset \m^C\cap(\xi+\c^{\perp})\end{equation}
for any $\xi\in \im\mu_{M,X_M}(X_1)$. Since there is a nilpotent
element in $\overline{\mu_{M,X_M}(X_1)}$, we see that the r.h.s. of
(\ref{eq:4.6:3}) coincides with $\m^C\cap\c^{\perp}$.  For brevity,
put $\s=\m^C\cap\c^\perp$. This is an ideal in $\m^C$.

Choose $x\in \widetilde{Z}'\cap X_1$ and put $\eta=\mu_{M,X_M}(x)$.
Then $\eta_s\in\z$ and $(\eta_s-\xi)+\eta_n\in \s$. Clearly,
$\c^C\subset \z(\m^C)$. Thus $[\eta_s-\xi,\eta_n]=0$ whence
$\eta_s-\xi=(\eta-\xi)_s\in \s$ and
\begin{equation}\label{eq:4.6:4}\mu_{M,X_M}(x)_s\in
\z\cap \c^\perp,\forall x\in \widetilde{Z}'\cap X_1.\end{equation}

\begin{Lem}\label{Lem:4.6.6}
$\m=\z+\t_0+[\m,\m]$.
\end{Lem}
\begin{proof}
It is enough to check that
\begin{equation}\label{eq:4.6:5}\t=\z+\t_0+\t_1,\t_1:=\t\cap [\m,\m],\end{equation} where $\t$ denotes a
Cartan subalgebra of $\lfr$. Recall that
\begin{align*}%\label{eq:4.6:6}
&\t=\z(\m)\oplus \t_1,\\%\label{eq:4.6:7}
&\z=\z(\m)\cap\a_{G,X}^{(X_L)}=\z(\m)\cap
(\z(\lfr)\cap\t_0^\perp)=\z(\m)\cap\t_0^\perp.
\end{align*}
Since $\z(\m),\t_1,\t_0$ are the Lie algebras of algebraic groups,
we see that $(\cdot,\cdot)$ is nondegenerate on
$\z(\m),\t_1,\t_0,\z$. To prove (\ref{eq:4.6:5}) it is enough to
note that $\t_0+\t_1=\z^\perp$.
\end{proof}

 If $\c\not\subset[\m,\m]+\t_0$, then, thanks to Lemma
\ref{Lem:4.6.6}, the r.h.s. of (\ref{eq:4.6:4}) is a proper subspace
in $\z$. Hence
$\psi_{M,X_M}(\widetilde{Z}')=\psi_{M,X_M}(\widetilde{Z'}\cap X_1)$
is not dense in $\z$. Since $\z\cap\im\psi_{M,X_M}$ is an open
subset in $\z$, we get a contradiction with
Corollary~\ref{Cor:4.5.3}.
% \end{proof}
\end{proof}

\begin{Prop}\label{Prop:4.6.3}
Let $X,L,X_L,M,X_M,\widehat{G}$ be as in
Proposition~\ref{Prop:4.6.1}, $\widehat{L}=L\cap \widehat{G}$. Let
$x\in X$ satisfy  conditions (a)-(d) of Proposition~\ref{Prop:4.6.1}
(for some $Z$) and $\widehat{X}$ be the model variety constructed by
$x$ as in (e). Then $\widehat{L}$ is the principal centralizer of
$\widehat{X}$ and there is an $\widehat{L}$-cross-section
$\widehat{X}_{\widehat{L}}$ of $\widehat{X}$ such that
$\a_{\widehat{G},\widehat{X}}^{(\widehat{X}_{\widehat{L}})}$ is a
$W_{G,X}^{(X_L)}\cap M/L$-stable subspace of $\a_{G,X}^{(X_L)}$ and
$W_{\widehat{G},\widehat{X}}^{(\widehat{X}_{\widehat{L}})}$ lies in
the image of $ W_{G,X}^{(X_L)}\cap M/L$ in
$\GL(\a_{\widehat{G},\widehat{X}}^{(\widehat{X}_{\widehat{L}})})$.
\end{Prop}
\begin{proof}
Recall, see Lemma \ref{Lem:1.9}, that $\widehat{L}$ is the principal
centralizer and $X_L$ is an $\widehat{L}$-cross-section of the
Hamiltonian $\widehat{G}$-variety $X_M$. Let $(H,\eta,V)$ denote the
determining triple of $X$ at $x$. Thanks to
Lemmas~\ref{Lem:4.3.4},\ref{Lem:4.3.5}, $(H\cap
\widehat{G},\eta_n,V/V^H\oplus V_0)$ is the determining triple of
the Hamiltonian $\widehat{G}$-variety $X_M$ at $x$, where $V_0$ is a
trivial  $H\cap\widehat{G}$-module. Put
$\widehat{X}':=M_{\widehat{G}}(H\cap \widehat{G},\eta_n,V/V^H\oplus
V_0)\cong \widehat{X}\times V_0$. It is enough to prove the analogue
of the assertion of the proposition for $\widehat{X}'$.

By Proposition~\ref{Prop:4.3.3}, there is a $\widehat{G}$-saturated
analytical open neighborhood  $O$ of $[1,(0,0)]$ in $\widehat{X}'$,
that is isomorphic (as a Hamiltonian $\widehat{G}$-manifold)  to a
saturated analytical neighborhood of $x$ in $(X_M)_{-\eta_s}$. One
may assume additionally that $O$ is connected.  By \cite{slice},
Lemma 5, $O_1:=\pi_{\widehat{G},\widehat{X}'}(O)$ is an open
neighborhood of $\pi_{\widehat{G},\widehat{X}'}([1,(0,0)])$ in
$\widehat{X}'\quo \widehat{G}$. Further, according to Example
\ref{Ex:2.2.4}, $\widehat{X}'$ is a conical Hamiltonian variety.
Replacing $O$ with a smaller neighborhood, we may assume that
$t.O\subset O$ for $0\leqslant t\leqslant 1$. Note that
$\widehat{L}$ is the principal centralizer of the Hamiltonian
$\widehat{G}$-variety $\widehat{X}'$. Since $\widehat{G}X_L$ is an
open subvariety of $X_M$ (in Zariski topology), we have $X_L\cap
O\neq \varnothing$. Choose an $\widehat{L}$-cross-section
$\widehat{X}'_{\widehat{L}}$
 of $\widehat{X}'$ such that some connected component of $X_L\cap O$
is contained in $\widehat{X}'_{\widehat{L}}\cap O$.

\begin{Lem}\label{Lem:4.6.4}
The manifold $\widehat{X}'_{\widehat{L}}\cap O$ is connected.
\end{Lem}
\begin{proof}[Proof of Lemma~\ref{Lem:4.6.4}]
Let $(\eta_n,h,f)$ be an  $\sl_2$-triple in $\widehat{\g}^{H\cap
\widehat{G}}$ generating the model variety $\widehat{X}'$. Note that
the action  $\C^\times:\widehat{X}'$ preserves
$\widehat{X}'_{\widehat{L}}$. Let $Y^0,Y^1$ be two distinct
connected components of $\widehat{X}'_{\widehat{L}}\cap O$, $y^i\in
Y^i, i=0,1,$ and $y^t, 0\leqslant t\leqslant 1,$ a continuous curve
connecting $y^0,y^1$ in $\widehat{X}'_{\widehat{L}}$. There is a
positive real $\tau<1$ such that $\tau y^t\in O$ for all
$t,0\leqslant t\leqslant 1$. Finally, note that $\tau_1y^i\in Y^i$
for all real $\tau_1$ such that $\tau\leqslant \tau_1\leqslant 1$
and $i=0,1$. Therefore $t\mapsto \tau y^t$ is a continuous curve in
$\widehat{X}'_{\widehat{L}}\cap O$ connecting points from $Y^0,Y^1$.
Contradiction.
\end{proof}

Now we can complete the proof of the proposition. One easily deduces
from Proposition \ref{Prop:4.3.3} that
$\a_{\widehat{G},\widehat{X}'}^{(\widehat{X}'_{\widehat{L}})}=\a_{\widehat{G},X_M}^{(X_L)}$.
The equalities
$W_{\widehat{G},X_M}^{(X_L)}=W_{M,X_M}^{(X_L)}=W_{G,X}^{(X_L)}\cap
M/L$ hold, see the discussion preceding Lemma \ref{Lem:1.9}.  By
Lemma~\ref{Lem:4.6.4},
$N_{\widehat{G}}(\widehat{L},X'_{\widehat{L}})=N_{\widehat{G}}(\widehat{L},X'_{\widehat{L}}\cap
O)$. It remains to recall that $X'_{\widehat{L}}\cap
O\hookrightarrow X_L$ whence
$N_{\widehat{G}}(\widehat{L},X'_{\widehat{L}}\cap O)\subset
N_{\widehat{G}}(\widehat{L},X_L)$.
\end{proof}

\begin{Rem}\label{Rem:2.8}
We use the notation of Proposition \ref{Prop:4.6.1}. Let
$\widehat{X}', \widehat{X}'_{\widehat{L}},O$ be as in the proof of
Proposition \ref{Prop:4.6.3}.  It can be checked using the
definitions of the morphisms $\widehat{\psi}_{\bullet,\bullet}$ that
the following diagram is commutative

\begin{picture}(165,40)
\put(2,32){$\widehat{X}$}\put(30,32){$\widehat{X}'$}\put(46,32){$O$}\put(60,32){$X_M$}
\put(95,32){$G*_MX_M$}\put(138,32){$X$}
\put(7,4){$\a_{\widehat{G},\widehat{X}}^{(\widehat{X}_{\widehat{L}})}/
W_{\widehat{G},\widehat{X}}^{(\widehat{X}_{\widehat{L}})}$}
\put(50,4){$\a_{\widehat{G},X_M}^{(X_L)}/
W_{\widehat{G},X_M}^{(X_L)}$} \put(90,4){$\a_{M,X_M}^{(X_L)}/
W_{M,X_M}^{(X_L)}$} \put(130,4){$\a_{G,X}^{(X_L)}/ W_{G,X}^{(X_L)}$}
\put(28,34){\vector(-1,0){20}}\put(44,34){\vector(-1,0){8}}
\put(50,34){\vector(1,0){9}}\put(67,34){\vector(1,0){25}}
\put(113,34){\vector(1,0){25}}\put(30,6){\vector(1,0){18}}
\put(89,6){\vector(-1,0){14}}\put(116,6){\vector(1,0){13}}
\put(5,31){\vector(1,-3){6.5}}\put(30,31){\vector(-1,-3){6.5}}
\put(63,31){\vector(0,-1){21}}\put(67,32){\vector(1,-1){23}}
\put(105,31){\vector(0,-1){21}}\put(140,31){\vector(0,-1){21}}
\put(141,22){\footnotesize $\widehat{\psi}_{G,X}$}
\put(79,22){\footnotesize $\widehat{\psi}_{M,X_M}$}
\put(64,17){\footnotesize $\widehat{\psi}_{\widehat{G},X_M}$}
\put(29,22){\footnotesize
$\widehat{\psi}_{\widehat{G},\widehat{X}'}$}\put(9,22){\footnotesize
$\widehat{\psi}_{\widehat{G},\widehat{X}}$}
\end{picture}

Let us explain the meaning of unmarked arrows. The morphism
$\widehat{X}'\cong \widehat{X}\times V^H \rightarrow \widehat{X}$ is
the projection along $V^H$. The maps $O\rightarrow \widehat{X}',X_M$
are open embeddings of complex analytical manifolds. The morphism
$X_M\rightarrow G*_MX_M$ is the embedding $x\mapsto [1,x]$ and the
morphism $G*_MX_M\rightarrow X$ is given by $[g,x]\mapsto gx$, it is
\'{e}tale by Proposition \ref{Prop:1.1}. One easily sees that
$\a_{\widehat{G},\widehat{X}}^{(\widehat{X}_{\widehat{L}})}=\a_{\widehat{G},X_M}^{(X_L)}$,
$\a_{M,X_M}^{(X_L)}=\a_{G,X}^{(X_L)}$, and
$\a_{\widehat{G},X_M}^{(X_L)}$ is the image of $\a_{M,X_M}^{(X_L)}$
under the orthogonal projection $\m\twoheadrightarrow\widehat{\g}$.
Moreover, it follows from Lemma \ref{Lem:1.9} and Proposition
\ref{Prop:4.6.3} that
$W_{\widehat{G},\widehat{X}}^{(\widehat{X}_{\widehat{L}})}\subset
W_{\widehat{G},X_M}^{(X_L)}$, $W_{M,X_M}^{(X_L)}=
W_{G,X}^{(X_L)}\cap M/L$, $W_{\widehat{G},X_M}^{(X_L)}\cong
W_{M,X_M}^{(X_L)}$. The morphism
$\a_{M,X_M}^{(X_L)}/W_{M,X_M}^{(X_L)}\rightarrow\a_{G,X}^{(X_L)}/W_{G,X}^{(X_L)}$
is the natural morphisms of quotients. The morphism
$\a_{\widehat{G},\widehat{X}}^{(\widehat{X}_{\widehat{L}})}/W_{\widehat{G},\widehat{X}}^{(\widehat{X}_{\widehat{L}})}\rightarrow
\a_{\widehat{G},X_M}^{(X_L)}/W_{\widehat{G},X_M}^{(X_L)}$ is the
composition of the natural morphism of quotients and the translation
by the projection of $\eta_s$ to $\a_{\widehat{G},X_M}^{(X_L)}$,
where $\eta$ is as in (e) of Proposition \ref{Prop:4.6.1}. The
morphism $\a_{M,X_M}^{(X_L)}/W_{M,X_M}^{(X_L)}\rightarrow
\a_{\widehat{G},X_M}^{(X_L)}/W_{\widehat{G},X_M}^{(X_L)}$ is induced
by the  the projection $\a_{M,X_M}^{(X_L)}\twoheadrightarrow
\a_{\widehat{G},X_M}^{(X_L)}$.
\end{Rem}

\subsection{The structure of Weyl groups of affine Hamiltonian varieties}\label{SUBSECTION_Weyl_aff1}
In this subsection $G$ is a connected reductive group, $T$ is a
maximal torus of $G$, $X$ is a conical affine Hamiltonian
$G$-variety satisfying the equivalent conditions of
Lemma~\ref{Lem:2.3.5}, and $X_T$ is a $T$-cross-section of $X$. The
goal of this subsection is to obtain some information about
$W_{G,X}^{(\cdot)}$ and some restrictions on a Hamiltonian
$G$-variety $X$ with a given Weyl group. All results are based on
Propositions \ref{Prop:4.6.1},\ref{Prop:4.6.5},\ref{Prop:4.6.3}.
These propositions allow one to reduce the study of
$W_{G,X}^{(\cdot)}$ to the case when  $G$ is semisimple and  $X$ is
a  model Hamiltonian variety $M_G(H,\eta,V)$ such that
$\cork_G(X)=0$ and $\eta$ is nilpotent. First of all, we need to
find out when the Weyl group of the last variety is trivial.

\begin{Prop}\label{Prop:5.2.1}
Let $G$ be a connected reductive  group, $H$ a reductive subgroup,
$\eta$ an $H$-invariant nilpotent element of $\g$, and $V$ a
symplectic $H$-module. Suppose $X:=M_G(H,\eta,V)$ satisfies the
equivalent conditions of Lemma \ref{Lem:2.3.5} and $\cork_G(X)=0$.
\begin{enumerate}
\item
If $W_{G,X}^{(\cdot)}=\{1\}$, then
\begin{itemize}
\item[(*)] $\eta=0, (G,G)\subset H$,
 $(G,G)\cong G_1\times \ldots\times G_k$ for some $k$, where $G_i\cong
\SL_2$. Moreover, the  $G$-modules $V/V^{(G,G)}$ and $ V_1\oplus
V_2\oplus\ldots\oplus V_k$ are isomorphic, where $V_i$ is the direct
sum of two copies of the two-dimensional irreducible $G_i$-module.
\end{itemize}
\item Conversely, if $G$ is semisimple, and $X$ satisfies (*), then
$W_{G,X}^{(\cdot)}=\{1\}$.
\end{enumerate}
\end{Prop}
\begin{proof}
Suppose, at first, that $G$ is semisimple. Let us prove the first
assertion.
 %REFEREE'S REMARK: explain why $cork_G(X)=0$ implies $X\quo G\cong C_{G,X}$ and why $\C[X]^G$ is generated
%elements of degree 2.

Since $\cork_G(X)=0$, we see that the field $\C(X)^G$ is Poisson commutative, compare with \cite{Vinberg}, Section 2.3.
It follows from Proposition \ref{Thm:2.2} that
$\widehat{\psi}_{G,X}\quo G:X\quo G\rightarrow
\a_{G,X}^{(\cdot)}/W_{G,X}^{(\cdot)}$ is an isomorphism. By Lemma \ref{Lem:2.51}, $\C[X]^G$
and $\C[\a_{G,X}^{(\cdot)}/W_{G,X}^{(\cdot)}]$ are isomorphic as graded algebras (where all elements
of $\a_{G,X}^{(X_L)*}$ are supposed to have degree 2, the grading on
$\C[X]^G$ is described in Remark \ref{Rem:2.2.5}).  Thus the
equality $W_{G,X}^{(\cdot)}=\{1\}$ is equivalent to the condition
that $\C[X]^G$ is generated by elements of degree  2.

The morphism $M_G(H^\circ,\eta,V)\rightarrow M_G(H,\eta,V),
[g,(u,v)]\mapsto [g,(u,v)],$ satisfies the assumptions of assertion
4 of Lemma \ref{Lem:1.9}. Thus
$W_{G,M_G(H^\circ,\eta,V)}^{(\cdot)}=\{1\}$ and we may assume that
$H$ is connected.

Put $U=(\z_\g(\eta)\cap \h^\perp)^*$. Let us equip the algebra
$\C[U\oplus V]$ with the grading described in
Remark~\ref{Rem:2.2.5}. The algebra $\C[U\oplus V]^H$ is generated
by elements of degree 2. Recall that $\eta\in U^*$ has degree 4.
Therefore $\eta=0$. Any element from $U^*\cong \h^\perp$ has degree
2. Therefore $U\quo H\cong U^H$, equivalently, $\C[U/U^H]^H=\C$. But
the $H$-module $U/U^H$ is orthogonal (that is, possesses a
nondegenerate $H$-invariant symmetric form) because  $U$ is
orthogonal. Hence $U=U^H$. Equivalently, $\g=\h+\g^\h$. In
particular, $\h$ is an ideal of $\g$. By Lemma \ref{Lem:2.3.5},
\begin{equation}\label{eq:5.2:1} \dim X= \dim\g+\rank \g.\end{equation}
But $\dim X=2\dim \g-2\dim\h+\dim V$. Note that $m_H(V)=\dim H$, for
$m_G(X)=\dim G$ and $\g/\h$ is a trivial $\h$-module.  Since $V$ is
a symplectic $H$-module, we have
\begin{equation}\label{eq:5.2:2}\dim V=
\dim\h+\rank\h+\cork_H(V).\end{equation} From (\ref{eq:5.2:1}),
(\ref{eq:5.2:2}) it follows that
\begin{equation}\label{eq:5.2:3}\dim \g+\rank\g= \dim X=2\dim G/H+\dim
V\geqslant 2\dim\g-\dim\h+\rank\h.\end{equation} We deduce from
(\ref{eq:5.2:3}) that $\dim \g-\rank\g\leqslant \dim\h-\rank\h$.
Since  $\h$ is an ideal of $\g$, the last inequality is equivalent
to $\g=\h$.

Let $G= G_1\ldots G_k$ be the decomposition into the locally direct
product of simple subgroups. By the discussion before Lemma
\ref{Lem:1.9}, $W_{G_i,V}^{(\cdot)}=\{1\}$. By
Propositions~\ref{Prop:4.6.1}, \ref{Prop:4.6.5}, there exists a
point $x\in \psi_{G_i,V}^{-1}(0)$ satisfying the conditions (a)-(f)
of those propositions (with $G,X$ replaced with $G_i,V$). Let
$(H_0,\eta_0,V_0)$ be the determining triple of the $G_i$-variety
$V$ at $x$. Thanks to Proposition~\ref{Prop:4.6.3},
$W_{G_i,M_{G_i}(H_0,\eta_0,V_0)}^{(\cdot)}=\{1\}$. By  assertion 1,
$\eta_0=0, H_0=G_i$ whence $V_0=V$. Further, $V/V^{G_i}$ is a
coisotropic  $G_i$-module of dimension $\dim\g_i+\rank\g_i$. Using
the classification of coisotropic modules obtained
in~\cite{cois},\cite{Knop7}, we get $G_i=\SL_2$. Since
$W_{G_i,V/V^{G_i}}^{(\cdot)}=\{1\}$, we see that $\C[V/V^{G_i}]$ is
generated by an element of degree 2. One easily deduces from this
that $V/V^{G_i}$ is isomorphic to the direct sum of two copies of
the irreducible 2-dimensional $\SL_2$-module.

Since $V/V^{G_i}\cong V^{G_i\skewperp}$ is a symplectic $G$-module
and $\Sp(V/V^{G_i})^{G_i}$ is a torus, the group $\prod_{j\neq
i}G_j$ acts trivially on $V/V^{G_i}$.    Note that
$(\bigoplus_{i=1}^k V^{G_i\skewperp})^{\skewperp}=\bigcap_{i=1}^k
V^{G_i}=V^{G}=0$. The last equality holds because $\cork_G(V)=0$. To
complete the proof of assertion 1 note that $\prod_{i=1}^k G_i$ acts
on $V$ effectively. It follows that the natural epimorphism
$\prod_{i=1}^k G_i\rightarrow G$ is an isomorphism.

Now suppose that $X$ is of the form indicated in (*). It is enough
to check  the equality $W_{G,X}^{(\cdot)}=\{1\}$ for $k=1$. Here the
equality follows from the observation that $\C[X]^G$ is generated by
an element of degree 2.

We proceed to the case when $G$ is not necessarily semisimple.
 Let $x$ be a point of $X$ satisfying conditions
(a)-(f) of Propositions \ref{Prop:4.6.1},\ref{Prop:4.6.3} for $M=G$,
$\widehat{G}=(G,G)$ and $\widehat{X}$ be the model variety
constructed by $\widehat{G},x$ in (e). By Proposition
\ref{Prop:4.6.3}, $W_{\widehat{G},\widehat{X}}^{(\cdot)}=\{1\}$.
Therefore $(G,G)=\prod_{i=1}^k G_i, \widehat{X}=\bigoplus_{i=1}^k
V_i$. Since any stabilizer of a point with closed $G$-orbit is
conjugate to a subgroup of $H$, we have $(G,G)\subset H,\eta=0$. By the above, there is a point $x\in V^{(G,G)}\hookrightarrow X$ such that
$V/\g_*x\cong \bigoplus_{i=1}^k V_i$. This observation completes the
proof.
\end{proof}

Now we are going to obtain a sufficient condition  for
$W_{G,X}^{(\cdot)}$  to intersect any subgroup of $W(\g)$ conjugate
to a certain fixed subgroup. To state the corresponding assertion we
need some definitions.

\begin{defi}\label{defi:5.2.2}
A subset $A\subset \Delta(\g) $ is called  {\it completely
perpendicular} if the following two conditions take place:
\begin{enumerate}
\item $(\alpha,\beta)=0$ for any $\alpha,\beta\in A$.
\item $\Span_{\R}(A)\cap \Delta(\g) =A\cup -A$.
\end{enumerate}
\end{defi}

For example, any one-element subset of $\Delta(\g) $ is completely
perpendicular.

\begin{defi}\label{Def:1.4.1}
A pair $(\h,V)$, where $\h$ is a reductive subalgebra of $\g$ and
$V$ is an $\h$-module, is said to be a $\g$-stratum. Two $\g$-strata
$(\h_1,V_1)$, $(\h_2,V_2)$ are called  {\it equivalent} if there
exist $g\in G$ and a linear isomorphism
$\varphi:V_1/V_1^{\h_1}\rightarrow V_2/V_2^{\h_2}$ such that
$\Ad(g)\h_1=\h_2$ and $ (Ad(g)\xi)\varphi(v_1)=\varphi(\xi v_1)$ for
all $\xi\in\h_1, v_1\in V_1/V_1^{\h_1}$.
\end{defi}

\begin{defi}\label{Def:1.4.2} Let $Y$ be a smooth affine variety and
$y\in Y$ a point with closed $G$-orbit. The pair $(\g_y,
T_yY/\g_*y)$ is called the $\g$-stratum of $y$. We say that $(\h,V)$
is a $\g$-stratum of $Y$ if $(\h,V)$ is equivalent to a $\g$-stratum of
a point of $Y$. In this case we write $(\h,V)\rightsquigarrow_\g Y$.
\end{defi}

\begin{Rem}\label{Rem:1.4.3} Let us justify the terminology. Pairs $(\h,V)$
do define some stratification of $Y\quo G$ by varieties with
quotient singularities. Besides, analogous objects were called
"strata" in \cite{Schwarz3}, where the term is borrowed from.
\end{Rem}

Let $A$ be a nonempty completely perpendicular subset of $\Delta(\g)
$. By $S^{(A)}$ we denote the $\g$-stratum $(\g^{(A)},
\sum_{\alpha\in A}V^{\alpha})$, where $V^\alpha$ is, by definition,
the direct sum of two copies of the two-dimensional irreducible
$\g^{(A)}/\g^{(A\setminus\{\alpha\})}$-module.

\begin{Cor}\label{Cor:5.2.3}
 If $W_{G,X}^{(X_T)}\cap W(\g^{(A)})=\{1\}$, then
$S^{(A)}\rightsquigarrow_\g X$.
\end{Cor}
\begin{proof}
Put $M=Z_G(\bigcap_{\alpha\in A}\ker\alpha)$. We remark  that
$G^{(A)}=(M,M)$. Choose a point $x\in X$ satisfying conditions
(a)-(f) of Propositions~\ref{Prop:4.6.1},\ref{Prop:4.6.5} for
general $\xi\in\z(\m)$. Let $(H,\eta,V)$ be the determining triple
of $X$ at $x$ and $\widehat{X}=M_{G^{(A)}}(H\cap
G^{(A)},\eta_n,V/V^H)$. By Proposition~\ref{Prop:4.6.3},
$W_{G^{(A)},\widehat{X}}^{(\cdot)}=\{1\}$. Using
Proposition~\ref{Prop:5.2.1}, we see that  $S^{(A)}$ is equivalent
to the $\g$-stratum of $x$.
\end{proof}

Now we  obtain some restriction on $W_{G,X}^{(\cdot)}$, namely, we
check that $W_{G,X}^{(\cdot)}$ is large in the sense of the
following definition.

\begin{defi}\label{defi:5.2.4}
A subgroup $\Gamma\subset W(\g)$ is said to be {\it large} if for
any two roots $\alpha,\beta\in \Delta(\g) $ such that $\beta\neq
\pm\alpha, (\alpha,\beta)\neq 0$ there exists  $\gamma\in
\R\alpha+\R\beta$ with $s_\gamma\in \Gamma$.
\end{defi}

\begin{Cor}\label{Cor:5.2.5}
The subgroup $W_{G,X}^{(X_T)}\subset W(\g)$ is large.
\end{Cor}
\begin{proof}
Assume the converse. Choose $\alpha,\beta\in \Delta(\g) $  such that
$\beta\neq \pm\alpha, (\alpha,\beta)\neq 0$ but  $s_\gamma\not\in
W_{G,X}^{(X_T)}$ for all  $\gamma\in \Delta(\g)\cap
(\R\alpha+\R\beta)$. Put $M=Z_G(\ker\alpha\cap\ker\beta)$. Note that
$(M,M)=G^{(\alpha,\beta)}$.  Let $x\in X$ satisfy conditions (a)-(f)
of Propositions \ref{Prop:4.6.1}, \ref{Prop:4.6.5} for general
$\xi\in\z(\m)$. Let $(H,\eta,V)$ be the determining triple of $X$ at
$x$. Put $\widehat{X}:=M_{G^{(\alpha,\beta)}}(H\cap
G^{(\alpha,\beta)},\eta_n,V/V^H)$. It follows from
Proposition~\ref{Prop:4.6.3} that
$W_{G^{(\alpha,\beta)},\widehat{X}}^{(\cdot)}$  contains no
reflection. Let $\widetilde{G}$ denote the simply connected covering
of $G^{(\alpha,\beta)}$. It is a simple simply connected group of
rank 2. Further, denote by $\widetilde{H}$ the connected normal
subgroup of $\widetilde{G}$ corresponding to $\h$. Put
$\widetilde{X}=M_{\widetilde{G}}(\widetilde{H},\eta_n,V/V^H)$. It is
a coisotropic variety. There is a natural morphism
$\widetilde{X}\rightarrow \widehat{X}$ satisfying the assumptions of
the fourth assertion of Lemma \ref{Lem:1.9}. Therefore the group
$W_{\widetilde{G},\widetilde{X}}^{(\cdot)}$ contains no reflection.
On the other hand, by Corollary~\ref{Cor:4.5.4}, the group
$W_{\widetilde{G},\widetilde{X}}^{(\cdot)}$ is generated by
reflections. Therefore
$W_{\widetilde{G},\widetilde{X}}^{(\cdot)}=\{1\}$. By
Proposition~\ref{Prop:5.2.1}, $\widetilde{G}$ is isomorphic to the
direct product of several copies of $\SL_2$. Since $\widetilde{G}$
is simple and of rank 2, this is absurd.
\end{proof}

Now let us describe large subgroups of $W(\g)$ for simple groups $G$
of types $A-E$.

%For a subgroup $\Gamma\subset W(\g)$ we denote by $\Gamma_r$ the
%subgroup of $\Gamma$ generated by all reflections from $\Gamma$. It
%follows from  Definition~\ref{defi:5.2.4} that the subgroup
%$\Gamma\subset W(\g)$ is large iff so is $\Gamma_r\subset W(\g)$. So
%to describe all large subgroups in $W(\g)$ it is enough to describe
%those generated by reflections.

Firstly, we consider the situation when  $\g$ is simple and has type
$A,D,E$, in other words, when all elements of  $\Delta(\g) $ are of
the same length.

Recall the classification of maximal proper root subsystems in
$\Delta(\g) $ (see~\cite{Dynkin}). We fix a system
$\alpha_1,\ldots,\alpha_r\in\Delta(\g) $ of simple roots. Let
$\alpha_0$ be the minimal (=lowest) root and $n_1,\ldots, n_r$
(uniquely determined) nonnegative integers satisfying
$\alpha_0+n_1\alpha_1+\ldots+n_r\alpha_r=0$. A proper root subsystem
$\Delta_0\subset\Delta(\g) $ is maximal iff it is $W(\g)$-conjugate
to one of the following root subsystems.

\begin{itemize}
\item[(a)] $\Span_\Z(\alpha_1,\ldots,\alpha_{i-1},\alpha_{i+1},\ldots,
\alpha_r)\cap\Delta(\g) $ for $n_i=1$.
\item[(b)]
$\Span_\Z
(\alpha_0,\alpha_1,\ldots,\alpha_{i-1},\alpha_{i+1},\ldots,
\alpha_r)\cap\Delta(\g) $ for prime $n_i$.
\end{itemize}

The number $n_i$ depends only on $\Delta_0$. We will call this
number the {\it characteristic} of $\Delta_0$.

For a proper subgroup  $\Gamma\subset W(\g)$ let $\Delta_\Gamma$
denote the set of all  $\alpha\in \Delta(\g) $ such that
$s_\alpha\in \Gamma$.

\begin{Prop}\label{Prop:5.2.6}
Let $\g$ be a simple Lie algebra of type $A,D,E$, $\rank\g>1$, and
$\Gamma$ a proper subgroup in  $W(\g)$.  Then $\Gamma$ is large iff
 $\Delta_\Gamma$ is a maximal proper root subsystem in
$\Delta(\g) $ of characteristic   1 or 2.
\end{Prop}
\begin{Lem}\label{Lem:5.2.9}
Let $\g$ be a simple Lie algebra of type $A,D,E$. Then
$\Delta_\Gamma$ is a root subsystem in $\Delta(\g)$ for any subgroup
$\Gamma\subset W(\g)$.
\end{Lem}
\begin{proof}
Let $\alpha,\beta\in \Delta(\g)$. Since all roots of $\Delta(\g) $
are of the same length, we see that $\alpha+\beta\in \Delta(\g) $,
(resp., $\alpha-\beta\in \Delta(\g) $) iff $(\alpha,\beta)< 0$,
(resp., $(\alpha,\beta)>0$).

%By its definition,  $\Delta_\Gamma$ has the following property:
%\begin{itemize}
%\item[(A)] $\Delta_\Gamma$ coincides with the set of all roots  $\alpha\in \Delta(\g)$
%such that $s_\alpha$ is contained in the subgroup of $W(\g)$
%generated by elements of $\Delta_\Gamma$.
%\end{itemize}

We need to check that $\alpha\in \Delta_{\Gamma}$ implies
$-\alpha\in \Delta_\Gamma$ and that $\alpha,\beta\in
\Delta_\Gamma,\alpha+\beta\in \Delta(\g) $ imply $\alpha+\beta\in
\Delta_\Gamma$. The first implication follows directly from the definition of $\Delta_\Gamma$. To
prove the second one we note that $\alpha+\beta=s_{\alpha}\beta$
whenever $\alpha,\beta,\alpha+\beta\in\Delta(\g)$, while $s_{s_\alpha\beta}=s_\alpha s_\beta s_\alpha\in \Gamma$
provided $\alpha,\beta\in \Delta_\Gamma$.
\end{proof}
\begin{proof}[Proof of Proposition~\ref{Prop:5.2.6}]
The subgroup $\Gamma\subset W(\g)$ is large iff
\begin{itemize}
\item[(A)]  $\{\alpha,\beta,\alpha+\beta\}\cap \Delta_\Gamma\neq\varnothing$
for all  $\alpha,\beta\in\Delta(\g)$ such that
$\alpha+\beta\in\Delta(\g)$.
\end{itemize}

One checks directly that a maximal root subsystem
$\Delta_\Gamma\subset \Delta(\g)$ of characteristic  1 or 2
satisfies (A).

Now let  $\Delta_\Gamma$ be a root subsystem of $\Delta(\g)$
satisfying (A). At first, assume that $\Delta_\Gamma$ is not
maximal. Let $\Delta_1$ be a maximal proper root subsystem of
$\Delta(\g) $ containing $\Delta_\Gamma$. Choose
$\alpha\in\Delta_1\setminus \Delta_\Gamma$. We see that
$\alpha+\beta\not\in\Delta(\g) $ for all $\beta\not\in\Delta_1$.
Otherwise
$\{\alpha,\beta,\alpha+\beta\}\cap\Delta_\Gamma=\varnothing$.
Analogously, $\alpha-\beta\not\in\Delta(\g)$. Therefore $\alpha\perp
\Delta\setminus\Delta_1$. Since the root system $\Delta$ is
irreducible, there is $\gamma\in \Delta$ such that
$(\alpha,\gamma)\neq 0, \gamma\not\perp \Delta\setminus\Delta_1$. By the above, any such $\gamma$ necessarily lies in $\Delta_\Gamma$.
Without loss of generality, we may assume that $(\alpha,\gamma)=-1$
whence $\alpha+\gamma\in \Delta$. Then, automatically,
$\alpha+\gamma\in \Delta_1\setminus\Delta_\Gamma$. It follows that
$\alpha+\gamma\perp\Delta\setminus\Delta_1$, which contradicts the
choice of $\gamma$.

It remains to show that the characteristic of $\Delta_\Gamma$ is
less than 3. Assume that
$\Delta_\Gamma=\Span_\Z\{\alpha_0,\ldots,\alpha_{i-1},\alpha_{i+1},\ldots,\alpha_r\}\cap\Delta(\g)
$, where $n_i>2$. Let $\pi_i^\vee$ denote the dual fundamental
weight corresponding to $\alpha_i$. The subset $\Delta_\Gamma\subset
\Delta(\g) $ coincides with the set of all $\alpha\in\Delta$ such
that $n_i$ divides $\pi_i^\vee(\alpha)$. So it is enough to check
that there are $\alpha,\beta\in \Delta(\g) $ such that
$\langle\pi_i^{\vee},\alpha\rangle=\langle\pi_i^{\vee},\beta\rangle=1$,
and $\alpha+\beta\in\Delta(\g) $. There is  $\gamma\in \Delta(\g) $
with $\langle\pi_i^{\vee},\gamma\rangle=2$. Choose such an element
$\gamma=\sum_{j=1}^r m_j\alpha_j$ such that $\sum m_j$ is minimal
possible. One sets $\alpha:=\alpha_i,\beta:=\gamma-\alpha_i\in
\Delta(\g) $.
\end{proof}

\begin{Cor}\label{Cor:5.2.7}
Suppose $\g$ is a simple classical Lie algebra. Then $\Gamma\subset
W(\g)$ is large iff $\Delta_\Gamma$ is listed in
Table~\ref{Tbl:5.2.8}.
\end{Cor}
 \begin{longtable}{|c|l|}\caption{Subsets $\Delta_\Gamma$
 for large subgroups
 $\Gamma\subset W(\g)$  when $\g$ is classical}\label{Tbl:5.2.8}\\\hline $\g$&$\Delta_\Gamma$\\\hline
$A_l,l\geqslant 2$& $\{\varepsilon_i-\varepsilon_j|i\neq j, i,j\in I\text{
or }i,j\not\in I\}, I\subsetneq \{1,\ldots,n+1\},
I\neq\varnothing$\\\hline $B_l,l\geqslant 3$&(a) $\{\pm
\varepsilon_i\pm\varepsilon_j|i\neq j, i,j\in I \text{ or }i,j\not\in
I\}\cup \{\pm\varepsilon_i| i\in I\},
I\subsetneq\{1,\ldots,n\}$\\&(b) $\{\pm
\varepsilon_i\pm\varepsilon_j|i\neq j, i,j\in I \text{ or }i,j\not\in
I\}\cup \{\pm\varepsilon_i| i\in \{1,2,\ldots,n\}\},
I\subsetneq\{1,\ldots,n\},I\neq\varnothing$\\&(c)
$\{\varepsilon_i-\varepsilon_j|i\neq j, i,j\in I \text{ or } i,j\not\in
I\}\cup \{\pm(\varepsilon_i+\varepsilon_j), i\in I,j\not\in I)\},
I\subset \{1,\ldots,n\}$
\\\hline $C_l,l\geqslant 2$&(a) $\{\pm \varepsilon_i\pm\varepsilon_j|i\neq j, i,j\in
I \text{ or }i,j\not\in I\}\cup \{\pm 2\varepsilon_i| i\in I\},
I\subsetneq\{1,\ldots,n\}$\\&(b) $\{\pm
\varepsilon_i\pm\varepsilon_j|i\neq j, i,j\in I \text{ or }i,j\not\in
I\}\cup \{\pm 2\varepsilon_i| i\in \{1,2,\ldots,n\}\},
I\subsetneq\{1,\ldots,n\}, I\neq\varnothing$\\&(c)
$\{\varepsilon_i-\varepsilon_j|i\neq j, i,j\in I \text{ or } i,j\not\in
I\}\cup \{\pm(\varepsilon_i+\varepsilon_j), i\in I,j\not\in I)\},
I\subset \{1,\ldots,n\}$ \\\hline $D_l,l\geqslant 3$& (a) $\{\pm
\varepsilon_i\pm\varepsilon_j|i\neq j, i,j\in I \text{ or }i,j\not\in I\},
I\neq\{1,\ldots,n\},\varnothing$\\&(b)
$\{\varepsilon_i-\varepsilon_j|i\neq j, i,j\in I \text{ or } i,j\not\in
I\}\cup \{\pm(\varepsilon_i+\varepsilon_j), i\in I,j\not\in I)\},
I\subset \{1,\ldots,n\}$
\\\hline
\end{longtable}
Note that some subsets $\Delta_\Gamma$ appear in
Table~\ref{Tbl:5.2.8} more than once.
\begin{proof}
For $\g$  of type $A_l$ or $D_l$ the required assertion stems
directly from Proposition~\ref{Prop:5.2.6}.

Suppose $\g\cong \sp_{2l},l\geqslant 2$. If $l=2$, then
$\Gamma\subset W(\g)$ is large iff $\Delta_\Gamma\neq \varnothing$.
All nonempty subsets $\Delta_\Gamma\subset\Delta(\sp_4) $ do appear
in Table~\ref{Tbl:5.2.8}. Now suppose $l>2$. Let  $\Delta_0$ denote
the subset of all short roots in $\Delta(\g) $ and $W_0$  the
subgroup of $W(\g)$ generated by $s_\alpha,\alpha\in\Delta_0$. Note
that $W_0$ is the Weyl group of the root system $D_l$. By the
definition of a large subgroup, the subgroup $\Gamma_0$ generated by
$s_\alpha,\alpha\in \Delta_\Gamma\cap \Delta_0,$ is large in $W_0$.
If $\Delta_0\cap\Delta_\Gamma$ is of type (a) (see Table
\ref{Tbl:5.2.8}), then $\Gamma$ is large in $W(\g)$ iff
$\Delta_\Gamma$ contains a long root. If $\Delta_0\cap
\Delta_\Gamma$ is of type (b) or $\Delta_0\subset\Delta_\Gamma$,
then $\Gamma_0$ is large in $W(\g)$. Since $\Gamma\subset
N_{W(\g)}(\Gamma_0)$, we see that large subgroups in $W(\g)$ are
precisely those presented in Table \ref{Tbl:5.2.8}.

The proof for   $\g\cong\so_{2l+1},l>2,$ follows easily from the
duality between the root systems $B_l,C_l$.
\end{proof}

\subsection{Examples of computation of Weyl
groups}\label{SUBSECTION_Weyl_computation} In this subsection we
classify  pairs $(G,V)$, where $G$ is a simple algebraic group, and
$V$ is a symplectic $G$-module such that $\defe_G(V)=\rank G$,
$W^{(\cdot)}_{G,V}\neq W(\g)$. The computation for $V\cong U\oplus
U^*$ (and, more generally, $X=T^*(G*_HV)$) is made in \cite{Weyl},
Section 5, so here we consider only the case $X\not\cong U^*\oplus
U$.

\begin{Lem}\label{Lem:5.9}
Let $G$ be a simple group, $X:=M_G(H,\eta,V)$, where $\eta$ is
nilpotent, satisfy the equivalent conditions of Lemma
\ref{Lem:2.3.5}. If $s_\alpha\not\in W_{G,X}^{(\cdot)}$ for some
$\alpha\in \Delta(\g)$, then there exist a subalgebra $\s\subset \h$
such that $\s\sim_G \g^{(\alpha)}$ and
\begin{equation}\label{eq:5.9:1}\frac{\tr_{U\oplus
V}(h^2)}{\tr_\h(h^2)}=1-\frac{4}{\tr_\h(h^2)}.\end{equation} Here
$U:=\z_\g(\eta)/\h$ and $h$ is a coroot in $\s$.
\end{Lem}
\begin{proof}
By Corollary \ref{Cor:5.2.3}, $S^{(\alpha)}\rightsquigarrow_\g X$.
Equivalently, there is a subalgebra $\s\subset \h$ such that
$\s\sim_G\g^{(\alpha)}$ and $(\s,\C^2\oplus \C^2)\rightsquigarrow_\h
U\oplus V$. The last condition implies that the $\s$-modules
$\h/\s\oplus (\C^2)^{\oplus 2}$ and $U\oplus V$ differ by a trivial
summand. Comparing the traces of $h^2$ on these two modules, we get
the claim.
\end{proof}

Here is the main result of this subsection.

\begin{Prop}\label{Prop:5.10}
Let $G$ be a simple group and $V$ a symplectic $G$-module satisfying
the equivalent conditions of Lemma \ref{Lem:2.3.5} such that
$V\not\cong U\oplus U^*$ for any $G$-module $U$. Then
$W_{G,V}^{(\cdot)}\neq W(\g)$ iff $V$ is contained in Table
\ref{Tbl:5.11}. The group $W_{G,V}^{(\cdot)}$ is presented in the forth
column of the table.
\end{Prop}

\begin{longtable}{|c|c|c|c|}
\caption{$G$-modules $V$ such that $W_{G,V}^{(\cdot)}\neq
W(\g)$}\label{Tbl:5.11}\\\hline
N&$\g$&$V$&$W_{G,V}^{(\cdot)}$\\\endfirsthead\hline
N&$\g$&$V$&$W_{G,V}^{(\cdot)}$\\\endhead\hline 1&$\g=\sl_6$&$
V=V(\pi_3)\oplus V(\pi_1)^{\oplus 2}\oplus V(\pi_5)^{\oplus
2}$&$A_1\times A_3$\\\hline 2& $\g=\sp_4$&$ V=V(\pi_1)\oplus
V(\pi_2)^{\oplus 2}$&$C_1\times C_1$\\\hline 3&$\g=\sp_6$&$
V=V(\pi_3)\oplus V(\pi_1)^{\oplus 2}$&$C_1\times C_2$\\\hline 4&
$\g=\so_{11}$&$ V=V(\pi_5)\oplus V(\pi_1)^{\oplus 4}$&$B_1\times
B_4$\\\hline 5& $\g=\so_{13}$&$ V=V(\pi_6)\oplus V(\pi_1)^{\oplus
2}$&$B_2\times B_4$\\\hline
\end{longtable}

In the fourth column we indicate the type of a root subsystem in
$\Delta(\g)$ such that the reflections corresponding to its roots
generate $W^{(\cdot)}_{G,V}$. By $B_1$ (resp., $C_1$) we mean a root
subsystem containing two opposite short (resp., long) roots in $B_n$
(resp., $C_n$). Root subsystems indicated in column 4 are determined
uniquely up to $W(\g)$-conjugacy.

\begin{proof}[Proof of Proposition \ref{Prop:5.10}]
By Corollary \ref{Cor:5.2.3}, $S^{(\alpha)}\rightsquigarrow_\g V$
for some $\alpha\in\Delta(\g)$. There is an $\SL_2$-stable prime
divisor $D'$ on $\C^2\oplus \C^2$ such that $m_{\SL_2}(D')=2$.
Applying the Luna slice theorem, we see that there is a prime
$G$-stable divisor $D$ on $V$
 such that $m_G(D)<\dim G$. All $G$-modules $V$ with $m_G(V)=\dim G$  possessing
 such a divisor $D$ were classified by Knop and Littelmann,
\cite{KL}. All such symplectic modules $V$ such that $V\not\cong
U\oplus U^*$ are presented in Table \ref{Tbl:5.11}. Let us show that
for these modules the inequality $W_{G,V}^{(\cdot)}\neq W(\g)$ does
hold.

{\it Case 1. $\g=\sl_6, V=V(\pi_3)\oplus V(\pi_1)^{\oplus 2}\oplus
V(\pi_5)^{\oplus 2}$.} We can consider $V$ as a symplectic
$\widetilde{G}:=\SL_6\times\SL_2$-module, where $\SL_2$ acts on
$V(\pi_1)^{\oplus 2}\oplus V(\pi_5)^{\oplus 2}$ as on $\C^2\otimes
(V(\pi_1)\oplus V(\pi_5))$). This module has a finite stabilizer in
general position and is coisotropic, see \cite{Knop7},\cite{cois}.
The Weyl group $W_{\widetilde{G},V}^{(\cdot)}$ was computed in
\cite{Knop7}, Table 12, it corresponds to the root system $A_1\times
A_1\times A_3$. By the discussion preceding Lemma \ref{Lem:1.9},
$W^{(\cdot)}_{\widetilde{G},V}=W_{\SL_2,V}^{(\cdot)}\times
W_{G,V}^{(\cdot)}$. It follows that $W_{G,V}^{(\cdot)}=A_1\times
A_3$.

{\it Case 2. $\g=\sp_4, V=V(\pi_1)\oplus V(\pi_2)^{\oplus 2}$.} We
can consider $V$ as a symplectic
$\widetilde{G}:=\Sp_4\times\C^\times$-module. Here $\C^\times$ acts
trivially on $V(\pi_1)$ and as $\SO_2$ on $V(\pi_2)^{\oplus 2}\cong
\C^2\otimes V(\pi_2)$. Again, this module is coisotropic and has a
finite stabilizer in general position. Using tables obtained in
\cite{Knop7}, we see that $W_{\widetilde{G},V}^{(\cdot)}\cong
A_1\oplus A_1$. But
$W_{G,V}^{(\cdot)}=W_{\widetilde{G},V}^{(\cdot)}$, see assertion 3
of Lemma \ref{Lem:1.9}. Using Lemma \ref{Lem:5.9}, we see that
$s_\alpha\in W_{G,V}^{(\cdot)}$ for all long roots $\alpha$.

{\it Case 3. $\g=\sp_6, V=V(\pi_3)\oplus V(\pi_1)^{\oplus 2}$.} One
argues exactly as in the previous case.

Before proceeding to the remaining two cases let us make some
remarks.

Firstly,  $s_\alpha\in W_{G,V}^{(\cdot)}$ for all short roots
$\alpha$. One checks this using Lemma \ref{Lem:5.9} (the fraction in
the l.h.s. of (\ref{eq:5.9:1}) (the index of the $G$-module $V$) can
be computed using  Table 1 of \cite{AEV}).

Since $s_\alpha\in W_{G,V}^{(\cdot)}$ for any short root $\alpha$,
it follows from  Corollary~\ref{Cor:5.2.5}, Proposition
\ref{Prop:5.2.6} that $W_{G,V}^{(\cdot)}$ is either the whole Weyl
group $W(\g)$ or is maximal among all proper subgroups generated by
reflections. The latter holds iff $\C[C_{G,V}]\cong \z(\C[V]^G)$
(the center of the Poisson algebra $\C[V]^G$) contains two linearly
independent elements of degree 4.

{\it Case 4. $\g=\so_{11}, V=V(\pi_5)\oplus V(\pi_1)^{\oplus 4}$.}

Note that $\Sp(V)^G\cong \Sp_4$. We consider the $\Z^2$-grading on
$\C[V]$ induced by the degrees with respect  to $V(\pi_5)$ and
$V(\pi_1)^{\oplus 4}$. By \cite{Schwarz}, $\C[V]^G$ is freely
generated by a 21-dimensional subspace $U\subset \C[V]^G$ such that
\begin{enumerate}
\item $U$ is an $\Sp_4$-submodule in $\C[V]^G$.
\item There is the decomposition $U=U_1\oplus U_2\oplus U_3\oplus
U_4$, where $U_1\cong S^2\C^4, U_2\cong \bigwedge^2\C^4, U_3\cong
\C^4, U_4\cong\C$ (isomorphisms of $\Sp_4$-modules).
\item $U_i, i=\overline{1,4},$ is homogeneous with respect to the
$\Z^2$-grading on $\C[V]$. The degrees of $U_1,U_2,U_3,U_4$ are
$(2,0),(2,2),(1,2),(0,4)$, respectively.
\end{enumerate}
Let $P_1,P_2$ denote the Poisson bivectors on $V(\pi_5)$ and
$V(\pi_1)^{\oplus 4}$ respectively. Now let $f_1,f_2$ be homogeneous
elements of $\C[V]$ of bidegrees, say $(d_1,d_1'),(d_2,d_2')$. Then
$\{f_1,f_2\}=\langle P_1,df_1\wedge df_2\rangle+\langle
P_2,df_1\wedge df_2\rangle$. The bidegrees of the first and the
second summand are, respectively,
$(d_1+d_2-2,d_1'+d_2'),(d_1+d_2,d_1'+d_2'-2)$. For instance,  if
$f_1\in \C[V(\pi_5)]$, then $\{f_1,f_2\}$ is homogeneous of bidegree
$(d_1+d_2-2,d_2')$.

Let us check that $U_4\subset \z(\C[V]^G)$. Since $V(\pi_5)$ and
$V(\pi_1)^{\oplus 4}$ are skew-orthogonal, we get $\{U_4,U_1\}=0$.
Suppose $\{U_4,U_2\}\neq \{0\}$. Since $U_4\subset
\C[V]^{G\times\Sp_4}$, we see that $\{U_4,U_2\}$ is isomorphic (as %REVIEWER'S REMARK: EXPLAIN
an $\Sp_4$-module) to a submodule of $\bigwedge^2\C^4$. But
$\{U_4,U_2\}$ consists of homogeneous elements of degree (2,4)
whence $\{U_4,U_2\}\subset U_1U_4+U_3^2$. Both summands are
isomorphic to $S^2\C^4$. This contradicts $\{U_4,U_2\}\neq \{0\}$.
Finally, the degree of $\{U_3,U_4\}$ equals $(1,4)$ whence
$\{U_3,U_4\}=0$.

Let $q$ denote a homogeneous element in $\C[\g]^G$ corresponding to
an invariant nondegenerate form. Then $\mu_{G,V}^*(q)$ is a
homogeneous element of $\C[V]^G$ of degree 4. It remains to check
that $\mu_{G,V}^*(q)\not\in U_4$. One checks easily that
$\im\mu_{G,V(\pi_1)^{\oplus 4}}$ contains an element $\xi$ such that
$(\xi,\xi)\neq 0$. If $v\in V(\pi_1)^{\oplus 4}\hookrightarrow V$ is
such that $\mu_{G,V}(v)=\xi$, then $q(\mu_{G,V}(v))\neq 0, f(v)=0$
for any $f\in U_4$. The last equality holds because $U_4\subset \C[V(\pi_5)]$.

By Corollary \ref{Cor:5.2.7}, $W_{G,V}^{(\cdot)}$ corresponds
either to $B_1\oplus B_4$ or to $ B_2\oplus B_3$. Thanks to
Corollary \ref{Cor:5.2.3}, it remains to prove that
$S^{(A)}\not\rightsquigarrow_\g V$, where
$A=\{\varepsilon_1-\varepsilon_2,\varepsilon_3-\varepsilon_4\}$. If
$S^{(A)}\rightsquigarrow_\g V$, then
\begin{equation}\label{eq:5.9:2}\dim
V^{\g^{(A)}}+\dim\g-\dim\n_\g(\g^{(A)})= \dim V-8.\end{equation} But
$\dim V(\pi_1)^{\g^{(A)}}=3, \dim V(\pi_5)^{\g^{(A)}}=8$ (recall
that the weight system of $V(\pi_5)$ consists of all weights of the
form $\pm\frac{1}{2}(\varepsilon_1\pm\ldots\pm\varepsilon_5)$
without multiplicities;
 $V(\pi_5)_\lambda\subset V(\pi_5)^{\g^{(A)}}$ iff $(\lambda,\varepsilon_1-\varepsilon_2)=
 (\lambda,\varepsilon_3-\varepsilon_4)=0$),
 $\dim\n_\g(\g^{(A)})=15$. So (\ref{eq:5.9:2}) does not hold.

{\it Case 5. $\g=\so_{13}, V=V(\pi_6)\oplus V(\pi_1)^{\oplus 2}$.}
By \cite{Schwarz}, the algebra $\C[V]^G$ is freely generated by 12
elements $f_{(2,0,0)},f_{(1,1,0)},f_{(0,2,0)}, f_{(0,0,4)},
f_{(0,0,8)},f_{(1,0,4)},f_{(0,1,4)},f_{(2,0,4)},f_{(1,1,4)},f_{(0,2,4)},
f_{(1,1,2)}, f_{(1,1,6)}$, where the lower index indicates the
grading with respect to the decomposition $V=V(\pi_1)\oplus
V(\pi_1)\oplus V(\pi_6)$. Note that $\Sp(V)^G\cong \SL_2$. The
elements $f_{(0,0,4)},f_{(0,0,8)},f_{(1,1,2)},f_{(1,1,6)}$ are
$\SL_2$-invariant, $\Span_\C (f_{(2,0,0)},f_{(1,1,0)},f_{(0,2,0)}),
\Span_\C( f_{(2,0,4)},f_{(1,1,4)},f_{(0,2,4)})\cong S^2\C^2,
\Span(f_{(1,0,4)},f_{(0,1,4)})\cong \C^2$.

Analogously to the previous case (i.e., using the grading and the
$\SL_2$-module structure of $\C[V]^G$), we check  that $f_{(0,0,4)},
f_{(0,0,8)}\in \z(\C[V]^G)$. Let $\mu_1,\mu_2,\mu$ denote the moment
maps for the actions $G:V(\pi_6),G:V(\pi_1)^{\oplus 2},V$. Clearly,
$\mu=\mu_1+\mu_2$. Further, put $f_2(\xi)=\tr(\xi^2),
f_4(\xi)=\tr(\xi^4), \xi\in\g$ (the traces are taken in the
tautological $\so_{13}$-module). We have shown that
$\mu_1^*(f_2),\mu_1^*(f_4)\in \z(\C[V]^G)$. On the other hand,
$\mu^*(f_2),\mu^*(f_4)\in \z(\C[V]^G)$. Let us check that
$\mu_1^*(f_2),\mu_1^*(f_4),\mu^*(f_2),\mu^*(f_4)$ are algebraically
independent. Analogously to the previous case one checks that
$\mu_1^*(f_2),\mu^*(f_2)$ are independent. It remains to check that
the equality
\begin{equation}\label{eq:5.9:3}
\begin{split}
&a\tr(\xi_1^4)+b\tr((\xi_1+\xi_2)^4)+c\tr(\xi_1^2)^2+d\tr(\xi_1^2)\tr((\xi_1+\xi_2)^2)+
e\tr((\xi_1+\xi_2)^2)^2=0,\\& \forall \xi_i\in \overline{\im\mu_i},
i=1,2,
\end{split}
\end{equation}
implies $a=b=0$. The isotropy subalgebras in general positions for
$V(\pi_6),V(\pi_1)^{\oplus 2}$ are
$\g^{(\alpha_1,\alpha_2,\alpha_4,\alpha_5)},\g^{(\alpha_2,\ldots,\alpha_6)}$,
respectively. By Lemma \ref{Lem:2.3.11},
$\overline{\im\mu_1}=\overline{G\Span_\C(\varepsilon_1+\varepsilon_2+\varepsilon_3,
\varepsilon_4+\varepsilon_5+\varepsilon_6)},$ $
\overline{\im\mu_2}=\overline{G\Span_\C(\varepsilon_1)}$. Since
$\tr(\xi_1^4)$ and $\tr(\xi_1^2)^2$ are not proportional, we get
$a+b=0$. Writing down the terms of (\ref{eq:5.9:3}) of bidegree
(3,1) with respect to $(\xi_1,\xi_2)$, we get
$$4b\tr(\xi_1^3\xi_2)+2d\tr(\xi_1^2)\tr(\xi_1\xi_2)+4e\tr(\xi_1^2)\tr(\xi_1\xi_2)=0.$$
Putting
$\xi_1=\varepsilon_1+\varepsilon_2+\varepsilon_3+i(\varepsilon_4+\varepsilon_5+\varepsilon_6),
\xi_2=\varepsilon_1$, we see that $b=0$.

To prove that the group $W_{G,V}^{(\cdot)}$ has the form indicated
in Table \ref{Tbl:5.11} it is enough to check  that
$S^{(A)}\not\rightsquigarrow_\g V$ for
$A=\{\varepsilon_1-\varepsilon_2,\varepsilon_3-\varepsilon_4,\varepsilon_5-\varepsilon_6\}$.
This is done analogously to the previous case.
\end{proof}

\section{Fibers of $\widehat{\psi}_{G,X}$ and untwisted
varieties}\label{SECTION_untwisted} Throughout this section $G,X$
are as in the previous one.

The goal of this section is to prove Theorem \ref{Thm:0.5} and
establish some examples of untwisted varieties. Subsection
\ref{SUBSECTION_Utw1} contains some technical results used in the
proof of Theorem \ref{Thm:0.5}. The proof itself is given in
Subsection \ref{SUBSECTION_Utw2}. In Subsection
\ref{SUBSECTION_Utw3} we describe some classes of untwisted
Hamiltonian varieties. We state a result by Knop that the cotangent
bundle of an affine variety is untwisted and show that any
symplectic module is untwisted. Finally, in Subsection
\ref{SUBSECTION_counterexamples} we give two counterexamples: of a
Hamiltonian variety not satisfying (Irr) and of a conical
coisotropic model variety not satisfying (Utw2). The former
counterexample is due to F. Knop.

\subsection{Reducedness of fibers of $\widehat{\psi}_{G,X}$}\label{SUBSECTION_Utw1}
\begin{Prop}\label{Prop:3.1}
 Let $\xi\in\a_{G,X}^{(\cdot)}$ be such that
$(W_{G,X}^{(\cdot)})_\xi=\{1\}$. Then the fiber
$\widehat{\psi}_{G,X}^{-1}(\pi_{W_{G,X}^{(\cdot)},\a_{G,X}^{(\cdot)}}(\xi))$
is reduced.
\end{Prop}
\begin{proof}
We preserve the notation of Proposition \ref{Prop:4.6.1} and Remark
\ref{Rem:2.8} and put $\widehat{G}=M$. The image of $\xi$ in
$\a_{G,X}^{(\cdot)}/W_{G,X}^{(\cdot)}$ is a smooth point. Thanks to
Theorem \ref{Thm:4.0.1}, the schematic fiber in interest is a local
complete intersection. So to verify that this fiber is reduced it is
enough to prove that it is generically reduced (see, for example,
\cite{Eisenbud}, Propositions 18.13,18.15). In other words, we need
to show that $\widehat{\psi}_{G,X}$ is smooth at any point $x\in X$
satisfying conditions (a)-(e) of Proposition \ref{Prop:4.6.1} for
any irreducible component $Z$. By Remark \ref{Rem:2.8},
$W_{M,X_M}^{(\cdot)}=\{1\}$. Using the commutative diagram of Remark
\ref{Rem:2.8}, we see that it is enough to prove the proposition in
the case when $\xi=0$ and $X=M_G(H,\eta,V)$, where $\eta$ is
nilpotent, $\cork_G(X)=0$, and $W_{G,X}^{(\cdot)}=\{1\}$. In this
case $\a_{G,X}^{(\cdot)}/W^{(\cdot)}_{G,X}\cong X\quo G$. Let
$X_L,L,G_0,X_0,L_0$ be as in Proposition \ref{Prop:4.2.2}. An
element of $N_G(L_0)$ acting trivially on $\a_{G,X}^{(X_L)}$ lies in
$Z_G(\a_{G,X}^{(X_L)})=L$. Thus (see the discussion preceding Lemma
\ref{Lem:1.9})
$W_{G_0,X_0}^{(\cdot)}/W_{G^\circ_0,X_0}^{(\cdot)}\cong
G_0/G_0^\circ$. So $G_0$ is connected. There is a point $y\in X_0$
of the form $[g,0]$, we may assume that $g=e$. It follows directly
from Example \ref{Ex:2.1.6} that $X_0\cong
M_{G_0}(N_H(L_0)/L_0,\eta,V^{L_0})$. Note that $\eta$ is a nilpotent
element of $ \g^{L_0}$  so  $\eta\in \g^{L_0}\cap\lfr_0^\perp\cong
\g_0$. By assertion 1 of Proposition \ref{Prop:4.2.2}, it is enough
to check that the fiber $\pi_{G_0,X_0}^{-1}(0)$ is generically
reduced. Replacing $(G,X)$ with $(G_0,X_0)$ we may assume, in
addition, that $X$ satisfies the equivalent conditions of Lemma
\ref{Lem:2.3.5}. It follows that $X$ satisfies condition (*) of
Proposition \ref{Prop:5.2.1}.

We may replace $G$ with a covering and assume that $G=T_0\times
H^\circ$, where $T_0$ is a torus. Further, by assertion 4 of Lemma
\ref{Lem:1.9}, $W_{G,\widetilde{X}}^{(\cdot)}=\{1\}$, where
$\widetilde{X}:=M_G(H^\circ,\eta,V)$. So we may replace $X$ with
$\widetilde{X}$ and assume that $H$ is connected. Since in this case
$X=T^*(T_0)\times V$, we reduce to the case $H=G,X=V$. Changing $G$
by a covering again, we may assume that $G\cong (G,G)\times Z$,
where $Z$ is a torus.

Recall that $\a_{G,V}\cong V\quo G$. The required claim will follow
if we show that the zero fibers of the morphisms
$\pi_{(G,G),V},\pi_{Z,V\quo (G,G)}$ are reduced. For the former
morphism this stems easily from the decomposition $V\cong
V^{(G,G)}\oplus\bigoplus_{i=1}^k V_i$. Put, for brevity,
$U_1:=V^{(G,G)},U_2:=\bigoplus_{i=1}^k V_i$. Note that
$\Sp(U_2)^{(G,G)}$ is a torus of dimension $k$ acting trivially on
$U_2\quo (G,G)$. So it remains to prove that $\pi_{Z,U_1}^{-1}(0)$
are reduced. Since $\cork_G(V)=0$, we have $\dim V=4k+2\dim Z$. It
follows that $\dim U_2=2\dim Z$. Further, by the above, $k+\dim
Z=\dim V\quo G=\dim U_1\quo Z+\dim U_2\quo (G,G)$, whence $\dim
U_1\quo Z=\dim Z$. Since $\C(U_1)^Z=\Quot(\C[U_1]^Z)$ (see
\cite{alg_hamil}, Theorem 1.2.9, for the proof in the general case),
it follows that $Z$ acts on $U_1$ locally effectively. Thus the
weight system  of $Z$ in $U_1$ coincides with
$\lambda_1,\ldots,\lambda_r,-\lambda_1,\ldots,-\lambda_r$, where
$\lambda_1,\ldots,\lambda_r$ form a basis in $\z^*$. Now the claim
is easy.
\end{proof}

\begin{Prop}\label{Prop:3.3}
Suppose $0\in\im\psi_{G,X}$, so that $0\in\a_{G,X}^{(X_L)}$. Let
$s\in W_{G,X}^{(\cdot)}$ be a reflection. Put
$\Gamma_s:=(\a_{G,X}^{(\cdot)})^s$ (the fixed point hyperplane of
$s$), $D_s=\pi_{W_{G,X}^{(\cdot)},\a_{G,X}^{(\cdot)}}(\Gamma_s)$.
Let $\widetilde{Z}$ be an irreducible component of
$\widehat{\psi}_{G,X}^{-1}(D_s)$.  Let $\xi$ be a general point in
$\Gamma_s$, $M:=Z_G(\xi)$, and $x\in \widetilde{Z}$ satisfy
conditions (a)-(f) of Propositions \ref{Prop:4.6.1},
\ref{Prop:4.6.5}. In the notation of those propositions put
$\widehat{G}=T_0(M,M)$ and let $\widehat{X}$ be as in (e) of
Proposition \ref{Prop:4.6.1}. Then the multiplicity of
$\widetilde{Z}$ in $\widehat{\psi}_{G,X}^{-1}(D_s)$ equals 1 or 2
and the latter holds iff
$W_{\widehat{G},\widehat{X}}^{(\cdot)}=\{1\}$.
\end{Prop}
\begin{proof}
We use the notation of Propositions
\ref{Prop:4.6.1},\ref{Prop:4.6.5} and Remark \ref{Rem:2.8}.  From
the choice of $M$ it follows that $s\in
W_{\widehat{G},X_M}^{(X_L)}$. Put
$\Gamma'_s=(\a_{\widehat{G},X_M}^{(X_L)})^s,$ and let $D_s'$ be the
image of $\Gamma_s'$ in
$\a_{\widehat{G},\widehat{X}}^{(\widehat{X}_{\widehat{L}})}/
W_{\widehat{G},\widehat{X}}^{(\widehat{X}_{\widehat{L}})}$. Let
$\widetilde{Z}_M$ be an irreducible component of $\widetilde{Z}\cap
X_M$ containing $x$. Also $\widetilde{Z}_M$ is  an irreducible
component of
$\widehat{\psi}_{\widehat{G},X_M}^{-1}(\pi_{W_{\widehat{G},X_M}^{(X_L)},\a_{\widehat{G},X_M}^{(X_L)}}(\Gamma_s'))$.
Let $\widehat{Z}'$ denote an irreducible component of
$\widehat{\psi}^{-1}_{\widehat{G},\widehat{X}'}(D_s')$ containing a
connected component of $O\cap \widetilde{Z}_M$, , where $O$ is as
in the proof of Proposition \ref{Prop:4.6.3}. Clearly,
$\widehat{Z}'=\widehat{Z}\times V^H$ for some subvariety
$\widehat{Z}\subset \widehat{X}$, which  is an irreducible component
of $\widehat{\psi}_{\widehat{G},\widehat{X}}^{-1}(D_s')$. From the
commutative diagram of Remark \ref{Rem:2.8} one deduces that
precisely one of the following possibilities holds:
\begin{enumerate}
\item $W_{\widehat{G},\widehat{X}}^{(\cdot)}\cong\Z_2$.
The multiplicity of $\widetilde{Z}$ in
$\widehat{\psi}_{G,X}^{-1}(D_s)$ equals the multiplicity of
$\widehat{Z}$ in
$\widehat{\psi}_{\widehat{G},\widehat{X}}^{-1}(D_s')$.
\item $W_{\widehat{G},\widehat{X}}^{(\cdot)}=\{1\}$. By Proposition
\ref{Prop:3.1}, the multiplicity of $\widehat{Z}$ in
$\widehat{\psi}_{\widehat{G},\widehat{X}}^{-1}(D_s')$ is one.
Further, the multiplicity of $\widetilde{Z}$ in
$\widehat{\psi}_{G,X}^{-1}(D_s)$ is 2.
\end{enumerate}

It remains to consider the first possibility. Note that, by
definition of $M$, one gets  $\Gamma_s\subset \z(\m)$ whence
$\Gamma_s'\subset \z(\widehat{\g})$. Since
$W_{\widehat{G},\widehat{X}}^{(\cdot)}\neq \{1\}$, it follows that
$\dim\a_{\widehat{G},\widehat{X}}^{(\cdot)}\cap
[\widehat{\g},\widehat{\g}]=1$. Further, note that $T_0$ acts
trivially on $\widehat{X}_{\widehat{L}}$. Since
$\a_{\widehat{G},\widehat{X}}^{(\widehat{X}_{\widehat{L}})}=\z(\widehat{\lfr})\cap\t_0^\perp$
and $\t_0$ projects surjectively to $\z(\widehat{\g})$ (recall that
$\widehat{\g}=\t_0+[\m,\m]$), we see that
$\a_{\widehat{G},\widehat{X}}^{(\widehat{X}_{\widehat{L}})}\cap\z(\widehat{\g})=\{0\}$.
Therefore
$\a_{\widehat{G},\widehat{X}}^{(\widehat{X}_{\widehat{L}})}\subset
[\g,\g],\dim\a_{\widehat{G},\widehat{X}}^{(\cdot)}=1,
Z(\widehat{G})^\circ\subset T_0$. From the last inclusion it follows
that $Z(\widehat{G})^\circ$ acts trivially on $\widehat{X}$.

Replacing $(G,X)$ with $((\widehat{G},\widehat{G}),\widehat{X})$, we
reduce the problem to the proof of the following claim.
\begin{itemize}
\item[(**)] Suppose that $X=M_G(H,\eta,V)$, $G$ is semisimple, $\eta$ is nilpotent,
$\cork_G(X)=0$, $\dim\a_{G,X}^{(\cdot)}=1,
W_{G,X}^{(\cdot)}=\{1,s\}$, where $s$ is a reflection.
 Then the fiber
$\widehat{\psi}_{G,X}^{-1}(0)$ is reduced.
\end{itemize}

As in the proof of Proposition \ref{Prop:3.1}, we see that
$\widehat{\psi}_{G,X}=\pi_{G,X}$. So it is enough to check that
$\pi_{H,U\oplus V}^{-1}(0)$ is reduced, where
$U=(\z_\g(\eta)/\h)^*$.  Note that $\C[U\oplus V]^H$ is generated by
an element of degree 4. Recall that $\eta\in (U^*)^H\subset
\C[U\oplus V]^H$ is of degree 4. So if $\eta\neq 0$, we are done. If
$\g\neq\h$, then there is an element  $q\in\C[\g/\h]^H$
corresponding to an $H$-invariant quadratic form on $\g/\h\cong
\h^{\perp}$. It has degree 4.  Any fiber of $q$ is reduced, since
$\dim\g/\h>1$. Therefore it remains to consider the case $H=G,X=V$.
The reducedness of fibers in this case follows from the observation
that a homogeneous generator of $\C[V]^G$ is irreducible.
\end{proof}

\begin{Prop}\label{Prop:3.4}
Again, we keep the notation of Proposition \ref{Prop:4.6.1} and
Remark \ref{Rem:2.8}. Suppose $X$ is  untwisted and $0\in
\im\widehat{\psi}_{G,X}$. Let $x$ be a point satisfying conditions
(a)-(d) of Proposition \ref{Prop:4.6.1} for a point $\xi_0\in
\a_{G,X}^{(X_L)}$. Then the Hamiltonian $\widehat{G}$-variety
$\widehat{X}$ is untwisted and
$W_{\widehat{G},\widehat{X}}^{(\widehat{X}_{\widehat{L}})}=(W_{G,X}^{(X_L)})_{\xi_0}$.
\end{Prop}
\begin{proof}
At first, we consider the case when $\widehat{G}=M$. Let $s$ be a
reflection lying in $(W_{G,X}^{(X_L)})_{\xi_0}$. Let $\Gamma_s,D_s$
have the same meaning as in Proposition \ref{Prop:3.3} and $D_s'$
denote the image of $\Gamma_s$ in
$\a_{M,\widehat{X}}^{(\widehat{X}_L)}/W_{M,\widehat{X}}^{(\widehat{X}_L)}$
(we write $\widehat{X}_L$ instead of $\widehat{X}_{\widehat{L}}$
because $L=\widehat{L}$). By Proposition \ref{Prop:3.1}, it is
enough to show that $\widehat{\psi}_{M,\widehat{X}}^{-1}(D'_s)$ is
reduced. Note that the last scheme is non-empty because
$\widehat{\psi}_{M,\widehat{X}}(x)\in D_s'$. Choose a component
$\widehat{Z}$ of $\widehat{\psi}_{M,\widehat{X}}^{-1}(D'_s)$. Note
that $\widehat{\psi}_{G,X}(\widehat{Z}\cap O)\subset D_s$. Choose a
component $Z$ of $\widehat{\psi}_{G,X}^{-1}(D_s)$ containing a
connected component of $\widehat{Z}\cap O$. Since $\widehat{G}=M$,
we see that the morphism
$\a_{M,X_M}^{(X_L)}/W_{M,X_M}^{(X_L)}\rightarrow
\a_{\widehat{G},X_M}^{(X_L)}/W_{\widehat{G},X_M}^{(X_L)}$ from the
commutative diagram of Remark \ref{Rem:2.8} can be inverted. So one
can consider the morphism $\widehat{\psi}:\widehat{X}\rightarrow
\a_{G,X}^{(X_L)}/W_{G,X}^{(X_L)}$ from this diagram. It follows from
the  diagram that the multiplicities of $\widehat{Z}$ in
$\widehat{\psi}^{-1}(D_s)$ and of $Z$ in
$\widehat{\psi}_{G,X}^{-1}(D_s)$ coincide. By (Utw2) the latter is
$1$. It follows that the multiplicity of $\widehat{Z}$ in
$\widehat{\psi}_{M,\widehat{X}}^{-1}(D_s')$ is 1 and that the
morphism
$\a_{M,\widehat{X}}^{(\widehat{X}_L)}/W_{M,\widehat{X}}^{(\widehat{X}_L)}\rightarrow
\a_{G,X}^{(X_L)}/W_{G,X}^{(X_L)}$ is unramified along $D_s'$. The
latter is equivalent to $s\in W_{M,\widehat{X}}^{(\widehat{X}_L)}$.
Since $W_{G,X}^{(X_L)}$ is generated by reflections, so is
$(W_{G,X}^{(X_L)})_{\xi_0}$. This completes the proof when
$\widehat{G}=M$.

We proceed to the general case.

\begin{Lem}\label{Lem:3.6}
Let $X_1,X_2$ be Hamiltonian $G$-varieties and
$\varphi:X_1\rightarrow X_2$ a Hamiltonian $G$-morphism. If $X_2$ is
untwisted, then so is $X_1$ and a natural morphism
$\varphi_0:C_{G,X_1}\rightarrow C_{G,X_2}$ induced by $\varphi$ is
\'{e}tale.
\end{Lem}
\begin{proof}[Proof of Lemma \ref{Lem:3.6}]
The morphism $\varphi_0$ is finite and dominant. The morphism
$\widetilde{\psi}_{G,X_2}\circ\varphi$ is smooth in codimension 1.
Therefore $\varphi_0$ is \'{e}tale in codimension 1. Since
$C_{G,X_2}$ is smooth, we can apply the Zariski-Nagata theorem on
the purity of branch locus. We see that $C_{G,X_1}$ is smooth and
$\varphi_0$ is \'{e}tale.
\end{proof}

Changing $M$ by a covering, we may assume that $M=\widehat{G}\times
T_0$, where $T_0$ is a torus. There is a natural Hamiltonian
morphism  $\rho:T^*(T_0)\times \widehat{X}\twoheadrightarrow
M_M(H,\eta,V)$. By Lemma \ref{Lem:3.6}, $T^*(T_0)\times \widehat{X}$
is untwisted and the natural morphism $C_{M, T^*(T_0)\times
\widehat{X}}\rightarrow C_{M, M_M(H,\eta,V)}$ is \'{e}tale. From
Proposition \ref{Thm:2.2} it follows that the Weyl groups of
$M_M(H,\eta,V)$ and $T^*(T_0)\times \widehat{X}$ coincide. This
implies all required claims.
\end{proof}

\subsection{Proof of Theorem \ref{Thm:0.5}}\label{SUBSECTION_Utw2}
\begin{Lem}\label{Lem:4.1}
If $X$ is an affine Hamiltonian variety satisfying (Utw1), then all schematic
fibers of $\widetilde{\psi}_{G,X}\quo G$ are Cohen-Macaulay.
\end{Lem}
\begin{proof}
By the Hochster-Roberts theorem (see, for instance, \cite{VP},
Theorem 3.19), $X\quo G$ is Cohen-Macaulay. Since $C_{G,X}$ is
smooth  and $\widetilde{\psi}_{G,X}\quo G$ is equidimensional (from
Proposition \ref{Lem:4.4.1}), we see that any fiber of
$\widetilde{\psi}_{G,X}\quo G$ is a locally complete intersection in
a Cohen-Macaulay scheme whence Cohen-Macaulay (\cite{Eisenbud},
Proposition 18.13).
\end{proof}

\begin{proof}[Proof of assertion 1]
Recall that $C_{G,X}\cong \a_{G,X}^{(\cdot)}/W_{G,X}^{(\cdot)}$,
Proposition \ref{Thm:2.2}.
 Thanks to Lemma \ref{Lem:4.1}, it remains to prove that any fiber of
 $\widehat{\psi}_{G,X}\quo G$ is smooth  in codimension 1. Since $X$
 is conical (of degree, say, $k$), there are actions $\C^\times:X\quo G,\a_{G,X}^{(\cdot)}/W_{G,X}^{(\cdot)}$
 such that the  former extends to a morphism $\C\times X\quo G:\rightarrow X\quo
 G$, the latter is induced by $\C^\times\times\a_{G,X}\rightarrow \a_{G,X}, (t,v)\mapsto t^k
 v$, and the morphism $\widehat{\psi}_{G,X}\quo G$ is
 $\C^\times$-equivariant. Applying a standard argument, we see
 that it is enough to prove that $\widehat{\psi}_{G,X}\quo
 G^{-1}(0)$ is smooth in codimension 1.

Let us use the notation of Corollary \ref{Cor:4.5.2}. Put
$\lambda:=0$ and choose $z\in Z_0$ and $x\in \pi_{G,X}^{-1}(z)$ with
closed $G$-orbit. Put $\widehat{X}:=M_{G}(H,\eta,V/V^H)$. By
Proposition \ref{Prop:3.4},
$\a_{G,\widehat{X}}^{(\cdot)}/W_{G,\widehat{X}}^{(\cdot)}\cong
\a_{G,X}^{(\cdot)}/W_{G,X}^{(\cdot)}$. Since
$\cork_{\widehat{G}}(\widehat{X})=0$, we have $\widehat{X}\quo
G\cong \a_{G,X}^{(\cdot)}/W_{G,X}^{(\cdot)}$. Taking quotients in
the commutative diagram of Remark \ref{Rem:2.8}, we get the
following commutative diagram

\begin{picture}(80,30)
\put(2,22){$V^H\times \widehat{X}\quo G$}\put(40,22){$O\quo
G$}\put(66,22){$X\quo
G$}\put(35,4){$\a_{G,X}^{(\cdot)}/W_{G,X}^{(\cdot)}$}
\put(38,23){\vector(-1,0){13}}\put(52,23){\vector(1,0){13}}
\put(22,20){\vector(1,-1){13}} \put(67,20){\vector(-1,-1){13}}
\put(30,24){\footnotesize $\hookleftarrow$}\put(56,24){\footnotesize
$\hookrightarrow$}\put(23,12){\footnotesize $\operatorname{pr}_2$}
\put(64,12){\footnotesize $\widehat{\psi}_{G,X}\quo G$}
\end{picture}

It follows that $\widehat{\psi}_{G,X}\quo G$ is smooth at $z$. Since
one may take an arbitrary point of $Z_0$ for $z$ we are done by
Proposition \ref{Prop:4.2.2}.
\end{proof}

\begin{proof}[Proof of assertion 2]
Choose a nonzero point $\xi\in\a_{G,X}^{(\cdot)}$ and put
$\a_0:=\C\xi$. Further, put $Y:=(X\quo
G)\times_{\a_{G,X}^{(\cdot)}/W_{G,X}^{(\cdot)}}\a_{G,X}^{(\cdot)},
Y_0:=(X\quo G)\times_{\a_{G,X}^{(\cdot)}/W_{G,X}^{(\cdot)}}\a_0$
(here $\a_0$ maps to $\a_{G,X}^{(\cdot)}/W_{G,X}^{(\cdot)}$ via the
composition $\a_0\hookrightarrow \a_{G,X}^{(\cdot)}\rightarrow
\a_{G,X}^{(\cdot)}/W_{G,X}^{(\cdot)}$).

Let us check that $Y$ is normal (as a scheme) and Cohen-Macaulay.
Indeed, the morphism $\a_{G,X}^{(\cdot)}\rightarrow
\a_{G,X}^{(\cdot)}/W_{G,X}^{(\cdot)}$ is flat, since $X$ satisfies
(Utw1). Therefore the morphism $Y\rightarrow X\quo G$ is flat. But,
as we have already remarked, $X\quo G$ is Cohen-Macaulay. By
Corollary 18.17 from \cite{Eisenbud}, $Y$ is Cohen-Macaulay. Note
that $X\quo G$ is smooth in codimension 1 over
$\a_{G,X}^{(\cdot)}/W_{G,X}^{(\cdot)}$. Hence $Y$ is smooth in
codimension 1 over $\a_{G,X}^{(\cdot)}$ hence normal.

Again, being a complete intersection in a Cohen-Macaulay variety,
$Y_0$ is Cohen-Macaulay. Similarly to the previous paragraph, $Y_0$
is normal.

Let us show that $\C[\a_0]$ is integrally closed in $\C[Y_0]$. Let
$\widetilde{\a}_0$ denote the spectrum of the integral closure of
$\C[\a_0]$ in $\C[Y]$. There is an action of $\C^\times$ on $Y_0$
lifted from $\C^\times:X\quo G$. The morphism $Y_0\rightarrow \a_0$
is $\C^\times$-equivariant. Therefore there is an action
$\C^\times:\widetilde{\a}_0$ contracting $\widetilde{\a}_0$ to the unique point  over $0\in\a_0$. It follows that
$\widetilde{\a}_0\cong\C$. Since the zero fiber of the morphism
$Y\rightarrow \a_0$ is normal, we see that the morphism
$\widetilde{\a}_0\rightarrow \a_0$ is \'{e}tale in 0. From the
$\C^\times$-equivariance it follows that it is an isomorphism.

Thus a general fiber of the morphism $Y_0\rightarrow \a_0$ is
irreducible. Thanks to the presence of $\C^\times$-action, the same
is true for any fiber but the zero one.  It follows easily from
(Con1),(Con2) that $(\widehat{\psi}_{G,X}\quo G)^{-1}(0)$ is
connected. Since $(\widehat{\psi}_{G,X}\quo G)^{-1}(0)$ is normal, it
is irreducible.
\end{proof}

\begin{proof}[Proof of assertion 3]
By Proposition \ref{Thm:2.2}, $X$ satisfies (Utw1). Assume that $X$
does not satisfy (Utw2).  By Proposition \ref{Prop:3.1}, there is
$s\in W_{G,X}^{(\cdot)}$ such that some irreducible component
$\widetilde{Z}\subset \widehat{\psi}_{G,X}^{-1}(D_s)$ (where, as
above, $D_s$ denotes the image of $(\a_{G,X}^{(\cdot)})^s$ in
$\a_{G,X}^{(\cdot)}/W_{G,X}^{(\cdot)}$) is of multiplicity 2. Put
$Y=\pi_{G,X}(\widetilde{Z})$. By Proposition \ref{Lem:4.4.1}, $Y$ is
an irreducible component of $(\widehat{\psi}_{G,X}\quo
G)^{-1}(D_s)$. It follows from (Irr) that
$Y=(\widehat{\psi}_{G,X}\quo G)^{-1}(D_s)$. Thanks to Theorem
\ref{Thm:4.0.1}, the set of closed orbits of any two components
$\widetilde{Z}_1,\widetilde{Z}_2\subset
\widehat{\psi}_{G,X}^{-1}(D_s)$ is the same. By Proposition
\ref{Prop:3.3}, the multiplicity of any component $\widetilde{Z}_1$
in $\widehat{\psi}_{G,X}^{-1}(D_s)$ is 2. Let $f\in
\C[\a_{G,X}^{(\cdot)}]^{W_{G,X}^{(\cdot)}}$ be such that $(f)=D_s$.
Let us remark that $f$ is not a square in $\C(X)$. Assume the
converse, let $f=f_1^2, f_1\in \C(X)$. Then $f_1\in
\C[C_{G,X}]=\C[\a_{G,X}^{(\cdot)}]^{W_{G,X}^{(\cdot)}}$ which is
absurd.

Put $\widetilde{A}:=\C[X][t]/(t^2-f)$. There is a natural morphism
of schemes $\rho:\Spec(\widetilde{A})\rightarrow X$. This morphism
is unramified over $X\setminus \widehat{\psi}_{G,X}^{-1}(D_s)$. Note
also that the group $\Z_2$ acts on $\Spec(\widetilde{A})$ so that
$\rho$ is the quotient for this action. Hence the restriction of
$\rho$ to any irreducible component of $\Spec(\widetilde{A})$ is
dominant. Since $f$ is not a square in $\C[X]$, we see that
$\Spec(\widetilde{A})$ is irreducible. Recall that  $(f)=2D$ for
some $D$. It follows that $\rho$ is unramified over
$\widehat{\psi}_{G,X}^{-1}(D_s)$. So $\rho$ is \'{e}tale in
codimension
1. %Let $\widetilde{X}$ denote the normalization of the reduced
%scheme associated with $\Spec(\widetilde{A})$. There is the natural
%morphism $\widetilde{\rho}:\widetilde{X}\rightarrow X$. As we have
%shown above, $\widetilde{\rho}$ is finite and \'{e}tale in codimension
%1.
Let $\widetilde{X}$ denote the normalization of
$\Spec(\widetilde{A})$ and $\widetilde{\rho}$ is the natural
morphism $\widetilde{X}\rightarrow X$. Since $X$ is smooth, we see
that $\widetilde{\rho}$ is also \'{e}tale in codimension 1. Besides,
$\widetilde{\rho}$ is finite.  Applying the Zariski-Nagata theorem
to $\widetilde{\rho}$, we obtain that $\widetilde{\rho}$ is
\'{e}tale. However, by our assumptions, $X$ is simply connected
whence $\widetilde{\rho}$ is an isomorphism. It follows that the
image of $t$ in $\C[\widetilde{X}]$ lies in $\C(X)$. This
contradicts the condition that $f$ is not a square in $\C(X)$.
\end{proof}

\begin{Rem}\label{Rem:0.6}
In fact, assertion 3 can be generalized to non simply connected
varieties. Namely, suppose $X$ satisfies (Irr). Then there exists an
untwisted conical Hamiltonian $G$-variety $\widetilde{X}$ and a free
action of a finite group $\Gamma$  on $\widetilde{X}$ by Hamiltonian
automorphisms (see Definition \ref{defi:3.5}) such that $X\cong
\widetilde{X}/\Gamma$ and $\pi_{\Gamma,X}:\widetilde{X}\rightarrow
X$ is a Hamiltonian morphism. The proof of this claim is similar to
that of assertion 3.
\end{Rem}

\subsection{Some classes of untwisted
varieties}\label{SUBSECTION_Utw3}
\begin{Prop}\label{Prop:5.0}
Let $X$ be  coisotropic, simply connected and conical  (e.g. a
symplectic vector space). Then $X$ is untwisted.
\end{Prop}
\begin{proof}
Thanks to Proposition \ref{Thm:2.2}, $X$ satisfies (Utw1).
Furthermore, $X$ obviously satisfies (Irr). Applying assertion 3 of
Theorem \ref{Thm:0.5}, we complete the proof.
\end{proof}

\begin{Thm}\label{Thm:5.1}
Let $X_0$ be a smooth irreducible affine $G$-variety. Then
$X:=T^*X_0$ is an untwisted Hamiltonian $G$-variety.
\end{Thm}
\begin{proof}
(Utw1) is checked in Satz 6.6 of \cite{Knop1}. (Utw2) follows from
\cite{Knop2}, Corollary 7.6.
\end{proof}

As the preprint \cite{Knop2} is not published, below we present
alternative proofs of Theorem \ref{Thm:5.1}.

\begin{Thm}\label{Thm:5.2}
Let $V$ be a symplectic $G$-module. Then $V$ is an untwisted
Hamiltonian $G$-variety.
\end{Thm}

We will prove this theorem after some auxiliary considerations.

\begin{Prop}\label{Prop:5.3}
Let $X$ be a conical Hamiltonian $G$-variety  and $G_0,X_0$ be such
as in the discussion preceding Proposition \ref{Prop:4.2.2}. If the
Hamiltonian $G_0^\circ$-variety $X_0$ satisfies (Irr), then so does
$X$.
\end{Prop}
\begin{proof}
The action $\C^\times:X$ preserves $X_0$ and so gives rise to the
structure of a conical Hamiltonian $G_0$-variety on $X_0$.  By
Proposition \ref{Prop:4.2.2}, the following  diagram, where the
horizontal arrows are quotient morphisms for the actions
$G_0/G_0^\circ$ on $X_0\quo G_0^\circ,
\a_{G,X}^{(\cdot)}/W_{G_0^\circ,X_0}^{(\cdot)}$, is commutative.

\begin{picture}(60,30)
\put(4,22){$X_0\quo G_0^\circ$}\put(45,22){$X\quo G$}
\put(2,4){$\a_{G,X}^{(\cdot)}/W_{G_0^\circ,X_0}^{(\cdot)}$}\put(42,4){$\a_{G,X}^{(\cdot)}/W_{G,X}^{(\cdot)}$}
\put(19,24){\vector(1,0){24}}\put(25,6){\vector(1,0){15}}
\put(9,20){\vector(0,-1){12}}\put(50,20){\vector(0,-1){12}}
\put(10,13){\footnotesize $\widehat{\psi}_{G_0^\circ,X_0}$}
\put(51,13){\footnotesize $\widehat{\psi}_{G,X}$}
\end{picture}

Choose $\xi\in\a_{G,X}^{(\cdot)}/W_{G,X}^{(\cdot)}$ and a point
$\xi'\in\a_{G_0,X_0}^{(\cdot)}/W_{G_0^\circ,X_0}^{(\cdot)}$ mapping
to $\xi$. By the previous commutative diagram,
$(\widehat{\psi}_{G,X}\quo G)^{-1}(\xi)$ is the quotient of
$\widehat{\psi}_{G_0^\circ,X_0}\quo G_0^{\circ-1}(\xi')$ by some
finite group. In particular, $(\widehat{\psi}_{G,X}\quo G)^{-1}(\xi)$
is irreducible.
\end{proof}

\begin{Prop}\label{Prop:5.4}
Let $X$ be an irreducible conical affine Hamiltonian $G$-variety
such that $\dim \a_{G,X}^{(\cdot)}=\rank G$. Let $G=Z(G)^\circ
G_1\ldots G_k$ be the decomposition of $G$ into the locally direct
product of simple normal subgroups and the unit component of the
center. If $X$ satisfies (Utw1) and is untwisted as a Hamiltonian
$G_i$-variety for any $i$, then $X$ is untwisted as a Hamiltonian
$G$-variety. Conversely, if $X$ is untwisted as a Hamiltonian
$G$-variety, then so is it as a Hamiltonian $G_i$-variety.
\end{Prop}
\begin{proof}
By Proposition \ref{Prop:3.1},  if $\widehat{\psi}_{G,X}^{-1}(D)$ is
not reduced for some divisor $D$ of
$\a_{G,X}^{(\cdot)}/W_{G,X}^{(\cdot)}$, then (in the notation of
Proposition \ref{Prop:3.3}) $D=D_s$ for some reflection $s\in
W_{G,X}^{(\cdot)}$. By assertion 5 of Lemma \ref{Lem:1.9}, if $s$ is
a reflection in $W_{G,X}^{(\cdot)}$, then $s$ is also a reflection
in $W_{G_i,X}^{(\cdot)}$ for some $i$ and vice versa. If
$W_{G,X}^{(\cdot)}$ is generated by reflections, then
$\widehat{\psi}_{G,X}$ is identified with the product of the
morphisms $\widehat{\psi}_{G_i,X}$. Now the proof is
straightforward.
\end{proof}

\begin{proof}[Proof of Theorem \ref{Thm:5.2}]
Applying Proposition \ref{Prop:5.3},  we reduce the proof to the
case when $V$ satisfies the equivalent conditions of Lemma
\ref{Lem:2.3.5}. Further, thanks to Proposition \ref{Prop:5.4}, we
may (and will) assume that $G$ is simple.

Suppose $V$ is not untwisted. By Proposition \ref{Thm:2.2},
$\widehat{\psi}_{G,V}$ is not smooth in codimension 1. By
Proposition \ref{Prop:3.3}, in the notation of that proposition, for
some $s\in W_{G,V}^{(\cdot)}$there is a point $x\in
\widehat{\psi}_{G,V}^{-1}(D_s)$ satisfying conditions (a)-(f) with
$W_{\widehat{G},\widehat{X}}^{(\cdot)}=\{1\}$. From  Remark
\ref{Rem:4.6.4} it follows that one can take $(M,M)$ for
$\widehat{G}$. Note that $\widehat{\g}\cong \sl_2$. By Proposition
\ref{Prop:5.2.1}, $\g_x=\sl_2, V/(\g_*x+V^{\g_x})\cong
\C^2\oplus\C^2$ (here $\C^2$ denotes the irreducible two-dimensional
$\sl_2$-module). In the proof of Proposition \ref{Prop:5.10} we have
seen that \begin{itemize}\item[(A)] all modules $V$ containing such
a point $x$ are presented in Table \ref{Tbl:5.11}; \item[(B)] if
$\alpha\in\Delta(\g)$ is such that $S^{(\alpha)}\rightsquigarrow_\g
V$, then $s_{w\alpha}\not\in W_{G,V}^{(\cdot)}$ for some $w\in
W(\g)$.
\end{itemize}

Let us choose a point $x\in \widehat{\psi}_{G,V}^{-1}(0)$ satisfying
conditions (a)-(e) of Proposition \ref{Prop:4.6.1}. Let
$\widehat{X}$ be the corresponding model $G$-variety. Let us check
that $(\widehat{\psi}_{G,V}\quo G)^{-1}(0)$ is irreducible. Clearly,
$(\widehat{\psi}_{G,V}\quo G)^{-1}(0)$ is connected. From
Proposition \ref{Prop:4.4.2} it follows that
$(\widehat{\psi}_{G,V}\quo G)^{-1}(0)$ is smooth in codimension 2
(as a variety). Applying the Hartshorne  connectedness theorem, we
see that $(\widehat{\psi}_{G,V}\quo G)^{-1}(0)$ is irreducible.

Therefore $\widehat{X}$ does not depend on the choice of $x$. If
$S^{(\alpha)}\rightsquigarrow_\g \widehat{X}$ for some $\alpha\in
\Delta(\g)$, then $S^{(\alpha)}\rightsquigarrow V$.  By Proposition
\ref{Prop:4.6.3}, $W_{G,\widehat{X}}^{(\cdot)}$ is $W(\g)$-conjugate
to a subgroup in $W_{G,V}^{(\cdot)}$. Let us check that these two
groups are, in fact, $W(\g)$-conjugate. By (B) and Corollary
\ref{Cor:5.2.3}, if $s_{w\alpha}\in W_{G,V}^{(\cdot)}$ for all $w\in
W(\g)$, the the same inclusions hold for $\widehat{X}$. But both
$W_{G,V}^{(\cdot)},W_{G,\widehat{X}}^{(\cdot)}$ are contained in
Table \ref{Tbl:5.2.8}. Therefore they are $W(\g)$-conjugate.

 It follows from the
commutative diagram of Remark \ref{Rem:2.8} that the fiber
$(\widehat{\psi}_{G,V}\quo G)^{-1}(0)$ is smooth in $x$. By
Proposition \ref{Prop:4.2.2}, this fiber is smooth in codimension 1.
Proceeding as in the proof of assertion 1 of Theorem \ref{Thm:0.5},
we see that any fiber of $\widehat{\psi}_{G,V}\quo G$ is normal.
Applying assertions 2,3 of Theorem \ref{Thm:0.5}, we see that $V$ is
untwisted.
\end{proof}

\begin{proof}[First alternative proof  of Theorem \ref{Thm:5.1}]
Suppose that (Utw2) does not hold. Let $Y$ be a prime divisor in $X$
consisting of singular points of $\widehat{\psi}_{G,X}$.

{\it Step 1.}
 Note that $Y$ is $\C^\times$-stable. Therefore there is a point
$x\in X_0\cap Y$ with closed $G$-orbit. Set $H:=G_x,
V_0:=T_xX_0/\g_*x, X'_0:=G*_HV_0, X':=T^*X_0$. As we noted in
\cite{slice}, $X'\cong M_G(H,0,V_0\oplus V_0^*)$. By \cite{Knop1},
Satz 6.5, $W_{G,X}^{(\cdot)}$ is conjugate to $W_{G,X'}^{(\cdot)}$.
Using Remark \ref{Rem:2.8}, we see that $X'$ does not satisfy
(Utw2). So we may assume that $X_0\cong G*_HV_0$. Also \cite{Knop1},
Satz 6.5, implies that $T^*\widetilde{X}_0$, where
$\widetilde{X}_0:=G*_{H^\circ}V_0$, is not untwisted. So we may
assume that $X$ is simply connected.

{\it Step 2.} Analogously to the proof of Theorem \ref{Thm:5.2}, we
may assume that $m_G(X)=\dim G$ and $G$ is simple. In \cite{Weyl},
Section 5 (see especially Lemma 5.4.1), we checked that condition
(B) of the proof of Theorem \ref{Thm:5.2} (with $V$ replaced by $X$)
holds. Now we can proceed as in that proof. Note that \cite{Weyl}
uses results of Section \ref{SECTION_Weyl} but not of Section
\ref{SECTION_untwisted}.
\end{proof}

Let us also present a proof that does not use the classification
results of \cite{Weyl}.

\begin{proof}[Second alternative proof  of Theorem \ref{Thm:5.1}]
Again, we may assume that $X$ is simply connected. Recall the notation $Y$ from
the previous proof.

 Set
$D:=\overline{\widetilde{\psi}_{G,X}(Y)}$. Let us check that any
irreducible component of $\widetilde{\psi}_{G,X}^{-1}(D)$ is of
multiplicity 2. We may assume that $Y\neq
\widetilde{\psi}_{G,X}^{-1}(D)$.  It follows from the connectedness theorem
of \cite{Knop12}, that $\widetilde{\psi}_{G,X}^{-1}(z)$ is connected
for any $z\in D$.  Applying  the Hartshorne connectedness theorem to
the Cohen-Macaulay scheme $(\widetilde{\psi}_{G,X}\quo G)^{-1}(D)$
and using the Knop connectedness theorem, we may assume that (probably, after
replacing $Y$ with another component of multiplicity 2) there is an
irreducible component $Y_1$ of $\widetilde{\psi}_{G,X}^{-1}(D)$ such
that
\begin{itemize}
\item
$\codim_{Y\quo G} (Y_1\quo G)\cap (Y\quo G)=1$,
\item
$\widetilde{\psi}_{G,X}(Y_1\cap Y)$ is dense in  $D$.
\item $Y_1$ is of multiplicity 1.
\end{itemize}
Choose a general point $\alpha\in D$. It follows from Proposition
\ref{Prop:4.4.2} that there is a point $x\in Y\cap Y_1$ satisfying
the conditions (a)-(f) of Propositions
\ref{Prop:4.6.1},\ref{Prop:4.6.5}. By Proposition \ref{Prop:3.3},
$Y_1$ has multiplicity 2, contradiction.

So any component of $\widetilde{\psi}_{G,X}^{-1}(D)$ has
multiplicity 2. Recall that by our assumptions $X$ is simply
connected. Arguing as in the last paragraph of the proof of
assertion 3 of Theorem \ref{Thm:0.5}, we get a contradiction.
\end{proof}

%\begin{Rem}\label{Rem:5.12}
%This remark contains a sketch of an alternative proof of Theorem
%\ref{Thm:5.1}. Let $X=T^*X_0$, where $X_0$ is a smooth affine
%$G$-variety. Applying Propositions \ref{Prop:5.3}, \ref{Prop:5.4} we
%see that it is enough to assume that $G$ is simple and
%$\a_{G,X}=\t$. Let $x_0\in X_0$ be a point with  closed $G$-orbit.
%Put $X_0':=G*_{G_{x_0}}(T_xX_0/\g_*x_0), X':=T^*X_0$. By Satz 6.5
%from \cite{Knop1},
%$\a_{G,X}^{(\cdot)}=\a_{G,X'}^{(\cdot)},W_{G,X}^{(\cdot)}=W_{G,X'}^{(\cdot)}$.
%Applying the Luna slice theorem, we see that the statements  of
%Theorem \ref{Thm:5.1} for $X,X'$ are equivalent. So we may assume
%that $X_0$ is a homogeneous vector bundle. It follows from Theorem
%7.1.2 in \cite{Weyl} that $S^{(\alpha)}\rightsquigarrow_\g X$ iff
%$s_{w\alpha}\not\in W_{G,X}^{(\cdot)}$ for some $w\in W$. Now the
%proof is analogous to the final part of the proof of Theorem
%\ref{Thm:5.2}.
%\end{Rem}

\subsection{Some counterexamples}\label{SUBSECTION_counterexamples}
The following example of a Hamiltonian variety not satisfying (Irr)
is due to F. Knop.

\begin{Ex}\label{Ex:5.13}
Put $G=\C^\times, X=\C^\times\times\C^\times\times\C\times\C$.
Choose coordinates $x_1,\ldots,x_4$ on $X$ corresponding to the
above decomposition. Define the action $G:X$ by
$t(x_1,x_2,x_3,x_4)=(tx_1,t^{-1}x_2,x_3,x_4)$. Put
$\alpha:=(x_1x_2^2-x_1^{-1}x_3^2)dx_1+x_4dx_3$. Clearly, $\alpha$ is
$G$-invariant. Put $$\omega:=-d\alpha=2x_1x_2dx_1\wedge
dx_2-2x_1^{-1}x_3 dx_1\wedge dx_3-dx_3\wedge  dx_4.$$ One checks
directly that $\omega$ is nondegenerate. The action $G:X$ is
Hamiltonian with $\mu_{G,X}(x)=\langle\alpha,
\frac{\partial}{\partial t}\rangle_x=x_1^2x_2^2-x_3^2$. It is clear
that $\mu_{G,X}^{-1}(a)$ is irreducible whenever $a\neq \{0\}$. It
follows that $C_{G,X}\cong \g$. On the other hand,
$\mu_{G,X}^{-1}(0)$ has two connected components.
\end{Ex}

We remark that the Hamiltonian variety in Example \ref{Ex:5.13} is
the smallest one in the sense that both group and variety have the
smallest possible dimensions.

Now let us present an example of a coisotropic conical  model
variety $X=M_G(H,\eta,V)$
 that is twisted. An example, where the group
$W_{G,X}^{(\cdot)}$ is not generated by reflections, can be found in
\cite{alg_hamil}, Subsection 5.10. In the following example
$W_{G,X}^{(\cdot)}$ is generated by reflections but (Utw2) does not
hold. Note that this example is very similar to that from
\cite{alg_hamil}.

\begin{Ex}\label{Ex:5.14}
Put $G=\SL_2\times \C^\times, X:=M_G(\Z_2\times\SL_2,0,\C^2\oplus
\C^2)$, where $\C^2$ denotes the two-dimensional irreducible
$\SL_2$-module with the symplectic form given by $(u,v)\mapsto
\det(u,v)$ and the nontrivial element $\sigma\in \Z_2\subset
\C^\times$ acts on $\C^2\oplus \C^2$ as follows:
$\sigma(v_1,v_2)=(v_2,-v_1)$. One easily checks that
$W_{G,X}^{(\cdot)}=N_G(T)/T\cong \Z_2$ and (Utw2) does not hold.
\end{Ex}

\section{Some open problems}\label{SECTION_open}
Firstly, we state two conjectures concerning property (Irr). Below
$G$ is a connected reductive group.

\begin{Conj}\label{Conj:1}
Any conical irreducible Hamiltonian $G$-variety $X$ satisfies (Irr).
\end{Conj}

The following conjecture is a weaker version of the first one.
\begin{Conj}\label{Conj:2}
$X=M_G(H,\eta,V)$, where $\eta$ is nilpotent, satisfies (Irr).
\end{Conj}

In virtue of the local cross-section and symplectic slice theorems
(Propositions \ref{Prop:1.1}, \ref{Prop:4.3.3}) one can deduce from
Conjecture \ref{Conj:2}  that any fiber of $\psi_{G,X}$ is normal
(as a variety).

Unlike the first conjecture, the second one can be reduced to some
case-by-case consideration. Let us sketch the scheme of this
reduction.

At first, one reduces the problem to the case when $X$ satisfies the
equivalent conditions of Lemma \ref{Lem:2.3.5} and then to the case
when $X$ is algebraically simply connected. Here one should check
that $X$ satisfies (Utw2). This will follow if one verifies the
following assertion:

\begin{itemize}
\item[(*)]
for any $\alpha\in \Delta(\g)$ such that
$S^{(\alpha)}\rightsquigarrow_\g X$ there is $w\in W(\g)$ such that
$s_{w\alpha}\not\in W_{G,X}^{(\cdot)}$.
\end{itemize}

Finally, it is enough to check (*) only for some special quadruples
$(G,H,\eta,V)$. By analogy with Section 7 of \cite{Weyl}, we call
such quadruples {\it quasiessential}. By definition, a quadruple
$(G,H,\eta,V)$ is quasiessential if $M_G(H,\eta,V)$ satisfies the
equivalent conditions of Lemma \ref{Lem:2.3.5} and for any ideal
$\h_1\subset \h$ there is $\alpha\in \Delta(\g)$ such that
$S^{(\alpha)}\rightsquigarrow_\g M_G(H,\eta,V)$ but
$S^{(\alpha)}\not\rightsquigarrow_\g M_G(H_1,\eta,V)$, where $H_1$
is a subgroup of $H$ corresponding to $\h_1$. It is not very
difficult to show that if $(G,H,\eta,V)$ is quasiessential, then $G$
is simple and $H$ is semisimple.

The next conjecture strengthens assertion 1 of Theorem
\ref{Thm:0.5}.

Note that any fiber of $\widetilde{\psi}_{G,X}\quo G$ has the
natural structure of a Poisson variety. The open stratum (in the
sense of Subsection \ref{SUBSECTION_affham4}) is symplectic.

\begin{Conj}\label{Conj:3}
Let $X=M_G(H,\eta,V)$, where $\eta$ is nilpotent, be untwisted. Then
any fiber $Y$ of $\widetilde{\psi}_{G,X}\quo G$ has symplectic
singularities. This means that there is a resolution of
singularities $\widetilde{Y}\rightarrow Y$ such that the symplectic
form on the smooth part of $Y$  is extended to some regular form on
$\widetilde{Y}$.
\end{Conj}

Finally, we would like to propose a conjecture giving an estimate on
dimensions of fibers of $\mu_{G,X}$.

\begin{Conj}\label{Conj:4}
Let $X$ be an irreducible affine Hamiltonian $G$-variety. Then
$\dim\mu_{G,X}^{-1}(\eta)\leqslant \dim X-(m_G(X)+\defe_G(X)+\dim
G\eta)/2$.
\end{Conj}

If $X$ is the cotangent bundle of a (not necessarily affine)
$G$-variety $X_0$ this conjecture follows from Vinberg's theorem on
the modality of the action of a Borel subgroup of $G$ on $X_0$, see
\cite{Vinberg_complexity}.

\section{Notation and conventions}\label{SECTION_Notation}
For an algebraic group denoted by a capital Latin letter we denote
its Lie algebra by the  corresponding small German letter. For roots
and weights of semisimple Lie algebras we use the notation of
\cite{VO}.
\begin{longtable}{p{5.5cm} p{10.5cm}}
\\$\sim_G$& the equivalence relation induced by an action of a group
$G$.
%\\$A^{(B)}$& the subset of all $B$-semiinvariant functions in a $G$-algebra $A$.\\
%$A^\times$& the group of all invertible elements of an algebra $A$.
%\\ $\Aut(\g)$ (resp., $\Int(\g)$)& the group of all (resp., inner)
%automorphisms of a Lie algebra $\g$.
\\ $C_{G,X}$& the spectrum of the integral closure of $\psi_{G,X}^*(\C[\g]^G)$ in
$\C[X]^G$.
\\ $\cork_G(X)$& the corank of a Hamiltonian $G$-variety $X$.
\\ $\defe_G(X)$& the defect of a Hamiltonian $G$-variety $X$.
\\$e_\alpha$& a nonzero element of the root subspace $\g^\alpha\subset \g$.
\\ $(f)$& the zero divisor of a rational function $f$.
\\ $(G,G)$ (resp., $[\g,\g]$)& the commutant of a group $G$
(resp., of a Lie algebra $\g$)
\\ $G^{\circ}$& the connected component of unit of an algebraic group $G$.
\\ $G*_HV$& the homogeneous bundle over $G/H$ with a fiber $V$.
\\ $[g,v]$& the equivalence class of $(g,v)$ in $G*_HV$.
\\ $G_x$& the stabilizer of $x\in X$ under an action
$G:X$.
\\ $\g^{\alpha}$& the root subspace of $\g$ corresponding to a root  $\alpha$.
\\ $\g^{(A)}$ (resp., $G^{(A)}$)& the subalgebra $\g$ generated by $\g^{\alpha}, \alpha\in A\cup -A$ (resp.,
the corresponding connected subgroup of $G$).
%\\$\Gr(V,d)$& the Grassmainian of $d$-dimensional subspaces of a
%vector space $V$.
\\ $m_G(X)$&$=\max_{x\in X}\dim Gx$.
\\ $N_G(H)$, (resp., $N_G(\h),\n_\g(\h)$)& the normalizer of an algebraic subgroup $H$ in
an algebraic group $G$ (resp., of a subalgebra $\h\subset \g$ in
$G$, of a subalgebra $\h\subset \g$ in $\g$).
\\ $N_G(H,Y)$& the stabilizer of $Y$ under the action of $N_G(H)$.
\\ $\Quot(A)$& the fraction field of $A$.
\\ $\rank(G)$& the rank of an algebraic group $G$.
%\\ $R(\lambda)$& the irreducible representation of a reductive algebraic group (or a
%reductive Lie algebra) corresponding to a highest weight $\lambda$.
%\\ $\Rad_u(H)$ (resp., $\Rad_u(\h)$)& the unipotent radical of an
%algebraic group $H$ (resp., of an algebraic Lie algebra $\h$).
\\
 $s_\alpha$& the reflection in a Euclidian space corresponding to a vector $\alpha$.
 \\ $\Span_{F}(A)$& the linear span of a subset $A$ of a module
over a field $F$.
\\ $\td A$& the transcendence degree of an algebra $A$.
\\ $U^{\skewperp}$& the skew-orthogonal complement to a subspace $U\subset V$ of
a symplectic vector space $V$.
\\ $V^\g$& $=\{v\in V| \g v=0\}$, where $\g$ is a Lie algebra and
$V$ is a $\g$-module.
%\\ $V_{(\lambda)}$& ��������� ���������� $\g$-������
%$V$ ($\g$~--- ����������� ������� ��), ��������������� �������� ����
%$\lambda$.
\\ $V(\lambda)$& the irreducible module of the highest weight $\lambda$
over a reductive algebraic group or a reductive Lie algebra.\\
$W(\g)$& the Weyl group of a reductive Lie algebra $\g$.\\
%$\X(G)$& the character lattice of an algebraic group $G$.\\
%$\X_{G}$& the weight lattice of a reductive algebraic group $G$.\\
$X^G$& the fixed point set for an action $G:X$.\\
$X\quo G$& the categorical quotient for an action $G:X$, where $G$
is a reductive group and $X$ is an affine  $G$-variety.
%$X/G$& a rational quotient for an action of an algebraic group $G$
%on a
%variety $X$.\\
%$X^{reg}$& the subvariety consisting of all smooth points of a variety $X$.\\
%$\#X$& the number of elements in a set $X$.
\\ $Z(G)$(resp., $\z(\g))$& the center of an algebraic group $G$ (resp., a Lie algebra $\g$).
\\   $Z_G(H)$, (resp., $Z_G(\h),\z_\g(\h)$)& the centralizer of a subgroup $H$
in an algebraic group $G$ (resp., of a subalgebra $\h\subset \g$ in
$G$, of a subalgebra $\h\subset \g$ in a Lie algebra  $\g$).
\\ $\alpha^\vee$& the dual root of $\alpha\in\Delta(\g)$.
\\  $\Delta(\g)$& the root system of a reductive Lie algebra $\g$.
%\\ $\lambda^*$& the dual highest weight to $\lambda$.
%\\ $\Lambda(\g)$& the root lattice of a reductive Lie algebra $\g$.
\\ $\mu_{G,X}$& the moment map for a Hamiltonian
$G$-variety $X$.
\\ $\xi_s$ (resp., $\xi_n$)& the semisimple (the nilpotent) part of
an element $\xi$ of an algebraic Lie algebra.
\\ $\xi_*$& the velocity vector field associated with $\xi\in\g$.
\\ $\Pi(\g)$& the system of simple roots for a reductive Lie algebra $\g$.
\\$\pi_{G,X}$& the (categorical) quotient morphism $X\rightarrow X\quo G$.
\\$\varphi\quo G$& the morphism of (categorical) quotients induced by a
$G$-equivariant morphism $\varphi$.
\\ $\varphi^*$& the homomorphism $\C[X_2]\rightarrow \C[X_1]$
induced  by a morphism $\varphi:X_1\rightarrow X_2$.
\\ $\psi_{G,X}$& $:=\pi_{G,\g}\circ\mu_{G,X}$.
\\ $\widetilde{\psi}_{G,X}$& the natural morphism $X\rightarrow
C_{G,X}$.
\\ $\widehat{\psi}_{G,X}$& the natural morphism $X\rightarrow
\a_{G,X}^{(\cdot)}/W_{G,X}^{(\cdot)}$.
\end{longtable}

\bigskip
{\Small Department of Mathematics, Massachusetts Institute of
Technology, 77 Massachusetts Avenue, Cambridge, MA 02139, USA.

\noindent E-mail address: ivanlosev@math.mit.edu}
\end{document}